\documentclass[11pt,leqno]{article}

\usepackage{latexsym}
\usepackage[all]{xy}
\usepackage{amsmath,amssymb,theorem}
\usepackage{amsfonts}

\newtheorem{thm}{Theorem}[section]
\newtheorem{prop}[thm]{Proposition}
\newtheorem{lem}[thm]{Lemma}
\newtheorem{df}[thm]{Definition}
\newtheorem{cor}[thm]{Corollary}

{\theorembodyfont{\rmfamily} \newtheorem{rmk}[thm]{Remark}}

\newcommand{\opi}[1]{\operatorname{{\it #1}}}
\newcommand{\op}[1]{\operatorname{#1}}
\newcommand{\sff}[1]{{\sf{#1}}}
\newcommand{\No}{\boldsymbol{N}}

\newcommand{\sch}[2]{\mbox{$({#1}\:\!\otimes\:\! {#2})^{\operatorname{sch}}$}}
\newcommand{\bB}{\boldsymbol{B}}

\DeclareFontFamily{U}{rsf}{}
\DeclareFontShape{U}{rsf}{m}{n}{
  <5> <6> rsfs5 <7> <8> <9> rsfs7 <10->  rsfs10}{}
\DeclareMathAlphabet{\mathscr}{U}{rsf}{m}{n}
\newcommand{\mycal}[1]{\mathscr{#1}}

\begin{document}

\title{Algebraic and topological aspects of the schematization functor}
\author{L. Katzarkov, T. Pantev, B. To\"en}
\date{}
\maketitle

\def\abstractname{Abstract}

\begin{abstract}
We study some basic properties of schematic homotopy types and the
schematization functor. We describe two different algebraic models for
schematic homotopy types: cosimplicial Hopf alegbras and equivariant
cosimplicial algebras, and provide explicit constructions of the
schematization functor for each of these models.  We also investigate
some standard properties of the schematization functor helpful for the
description of the schematization of smooth projective complex
varieties. In a companion paper these results are used in the
construction of a non-abelian Hodge structure on the schematic
homotopy type of a smooth projective variety.
\end{abstract}


\tableofcontents

\newpage

\section{Introduction}

The schematization functor is a device which converts a pair $(X,k)$,
consisting of topological space $X$ and field $k$, into an
algebraic-geometric object $\sch{X}{k}$. The characteristic property
of $\sch{X}{k}$ is that it encodes the $k$-linear part of the homotopy
type of $X$. Namely $\sch{X}{k}$ captures all the information about
local systems of $k$-vector spaces on $X$ and their cohomology. For a
simply connected $X$ and $k = {\mathbb Q}$ (respectively $k = {\mathbb
F}_{p}$) the object $\sch{X}{k}$ is a model for the rational
(repectively $p$-adic) homotopy type of $X$. An important advantage of
$\sch{X}{k}$ is that it makes sense for non-simply connected $X$ and
that it detects non-nilpotent information.

The object $\sch{X}{k}$ belongs to a special class of
algebraic $\infty$-stacks over $k$, called schematic homotopy types
\cite{t1}. The existence and functoriality of the schematization
$\sch{X}{k}$ are proven in \cite{t1} but the construction is somewhat
abstract and unwieldy. In this paper we supplement \cite{t1} by
describing explicit algebraic models for $\sch{X}{k}$. We also study
in detail some of the basic properties of the schematization - the Van
Kampen theorem, the schematization of homotopy fibers, de Rham models,
etc.  These results are used in an essential way in \cite{kpt} in
which we construct mixed Hodge structures on schematizations of smooth
complex projective varieties.

The paper is organized in three parts. In section \ref{sec:review} we
briefly review the definition of schematic homotopy types and the
existence results for the schematization functor from \cite{t1}. In
section \ref{sec:models} we present two different algebraic models for
the schematization of a space, namely equivariant cosimplicial
$k$-algebras and cosimplicial Hopf $k$-algebras. These generalize two
well known ways for modelling rational homotopy types - via dg
algebras and via nilpotent dg Hopf algebras \cite{ta}. Each of the two
models utilizes a different facet of the homotopy theory of a space
$X$. The equivariant cosimplicial $k$-algebras codify the cohomology
of $X$ with $k$-local system coefficients together with their
cup-product structure, whereas the cosimplicial Hopf $k$-algebra is
the algebra of representative functions on the simplicial loop group
associated with $X$ via Kan's construction \cite[Section~V.5]{gj}. 
The two models have different ranges of aplicability. For
instance the cosimplicial Hopf algebra model is needed to construct
the weight tower for the mixed Hodge structure (MHS) on the
schematization of a smooth projective variety \cite{kpt}. On the other
hand the Hodge decomposition on the schematic homotopy type of a
smooth projective variety is defined in terms of the equivariant
cosimplicial model.

Section \ref{sec:properties} gathers some useful facts about the
behavior of the schematization functor. As an application of the
equivariant cosimplicial algebra model, we describe the schematization
of differentiable manifolds in terms of de Rham complexes of flat
connections. This description generalizes a theorem of Sullivan's
\cite{su} expressing the real homotopy theory of a manifold in terms
of its de Rham complex.  To facilitate computations we prove a
schematic analogue of the van Kampen theorem which allows us to build
schematizations by gluing schematizations of local pieces.  Since in
the algebraic-geometric setting we can not use contractible
neighborhoods as the building blocks we are forced to study the
schematizations of Artin neighborhoods and more generally of
$K(\pi,1)$'s. This leads to the notion of $k$-algebraically good
groups which are precisely the groups $\Gamma$ with property that the
schematization of $K(\Gamma,1)$ has no higher homotopy groups. We give
various examples of such groups and prove that the fundamental groups
of Artin neighborhoods are algebraically good. Note that the analogous
statement in rational homotopy theory is unknown and probably
false. Finally we prove two exactness properties of the schematization
functor. First we establish a Lefschetz type right exactness property
of schematizations, useful for understanding homotopy types of
hyperplane sections. We also give sufficient conditions under which
the schematization commutes with taking homotopy fibers. This
criterion is used in the construction of new examples \cite{kpt}
of non-K\"{a}hler homotopy types.

\

\bigskip

\noindent
{\bf Acknowledgements:} We would like to thank C.Simpson for his
constant encouragement and for many stimualting conversations on the
subject of this paper. We would also like to thank A.Beilinson and
P.Deligne for their help with the discussion about good groups. We
thank M.Olsson for pointing out a mistake in a preliminary
version of the text and for several helpful comments. We want
to thank MSRI for providing excellent working conditions during the
semester on `Intersection theory, algebraic stacks and non-abelian
Hodge theory' when most of this work was done.

Finally, we are grateful to the referee for his useful comments which helped
improving the readability of this paper.

\

\bigskip

\noindent
\textit{Conventions:} We will fix two universes $\mathbb{U}$ and
$\mathbb{V}$, with $\mathbb{U} \in \mathbb{V}$, and we assume
$\mathbb{N} \in \mathbb{U}$.

We denote by $k$ a base field in $\mathbb{U}$.  We consider
$\sff{Aff}$, the category of affine schemes over $\op{Spec}\, k$
belonging to $\mathbb{U}$. The category $\sff{Aff}$ is a
$\mathbb{V}$-small category. We endow it with the faithfully flat and
quasi-compact topology, and consider the model category $\sff{SPr}(k)$
of presheaves of $\mathbb{V}$-simplicial sets on the site
$(\sff{Aff},\op{fpqc})$. We will use the local projective model structure on
simplicial presheaves described in \cite{bl,t1} (note that as the site
$\sff{Aff}$ is $\mathbb{V}$-small, the model category $\sff{SPr}(k)$
exists).  We denote by $\sff{SPr}(k)_{*}$ the model category of
pointed objects in $\sff{SPr}(k)$. The expression
\emph{stacks} will always refer to 
objects in $\op{Ho}(\sff{SPr}(k))$. In the same way, \emph{morphisms
of stacks} refers to morphisms in $\op{Ho}(\sff{SPr}(k))$. \\

\section{Review of the schematization functor} \label{sec:review}

In this section, we review the theory of affine stacks and
schematic homotopy types introduced in
\cite{t1}. The main goal is to recall the theory and fix
the notations and the terminology.

We will denote by $\sff{Alg}^{\Delta}$ the category of cosimplicial
(commutative) $k$-algebras that belong to the universe $\mathbb{V}$.
The category $\sff{Alg}^{\Delta}$ is endowed with a simplicial closed
model category structure for which the fibrations are the epimorphisms
and the equivalences are the quasi-isomorphisms. This model category
is known to be cofibrantly generated, and even finitely generated
\cite[Theorem~2.1.2]{t1}.
\

\bigskip

\noindent
There is a natural spectrum functor
\[
\op{Spec} : (\sff{Alg}^{\Delta})^{\op{op}} \longrightarrow
\sff{SPr}(k),
\]
defined by the formula
\[
\xymatrix@R=1pt@C=9pt{
\op{Spec} A \; : & \sff{Aff}^{\op{op}} \ar[r] &
\sff{SSet} \\
&  \op{Spec} B   \ar@{|->}[r] & \underline{\op{Hom}}(A,B),
}
\]
As usual $\underline{\op{Hom}}(A,B)$ denotes the simplicial set of
morphisms from the cosimplicial algebra $A$ to the algebra $B$.
Explicitly, if $A$ is given by a cosimplicial object $[n] \mapsto
A_{n}$, then the presheaf of $n$-simplices
of $\op{Spec} A$ is given by $(\op{Spec} B) \mapsto \op{Hom}(A_{n},B)$.

The functor $\op{Spec}$ is a right Quillen functor, and its
right derived functor is denoted by
\[
\mathbb{R}\op{Spec} : \op{Ho}(\sff{Alg}^{\Delta})^{\op{op}} \longrightarrow
\op{Ho}(\sff{SPr}(k)).
\]
The restriction of $\mathbb{R}\op{Spec}$ to the full sub-category of
$\op{Ho}(\sff{Alg}^{\Delta})$ consisting of objects isomorphic to a
cosimplicial algebra in $\mathbb{U}$ is fully faithful (see
\cite[Corollary~2.2.3]{t1}). By defintion, an affine stack is an object
$F \in \op{Ho}(\sff{SPr}(k))$ isomorphic to an object of the form
$\mathbb{R}\op{Spec} A$, for some cosimplicial algebra $A$ in
$\mathbb{U}$. Moreover, by \cite[Theorem~2.4.1,2.4.5]{t1} the following
conditions are equivalent for a given pointed stack $F$:

\begin{enumerate}
\item The pointed stack $F$ is affine and connected.
\item The pointed stack $F$ is connected and for all $i>0$ the sheaf
$\pi_{i}(F,*)$ is represented by an affine unipotent group scheme.
\item There exist a cohomologically connected cosimplicial algebra
$A$ (i.e. $H^{0}(A)\simeq k$), which
belongs to $\mathbb{U}$, and such that
$F\simeq \mathbb{R}\opi{Spec} A$.
\end{enumerate}

Recall next that for a pointed simplicial presheaf $F$, one can
define its simplicial presheaf of loops
$\Omega_{*}F$. The functor $\Omega_{*} : \sff{SPr}_{*}(k)
\longrightarrow \sff{SPr}(k)$ is right Quillen,
and can be derived to a functor defined on the level of homotopy categories
\[
\mathbb{R}\Omega_{*}F : \op{Ho}(\sff{SPr}_{*}(k)) \longrightarrow
\op{Ho}(\sff{SPr}(k)).
\]

A pointed and connected stack $F \in \op{Ho}(\sff{SPr}_{*}(k))$ is
called a {\em pointed affine $\infty$-gerbe} if
the loop stack $\mathbb{R}\Omega_{*}F \in \op{Ho}(\sff{SPr}(k))$ is affine.
A pointed schematic homotopy type is a pointed affine
$\infty$-gerbe which in addition satisfies
a cohomological condition (see
\cite[Def. 3.1.2]{t1} for details).

\

\bigskip

The main result on affine stacks that we need is the existence theorem of
\cite{t1}. Embed the category $\sff{SSet}$ into the category
$\sff{SPr}(k)$ by viewing a simplicial set $X$ as a constant simplicial
presheaf on $(\sff{Aff},\op{ffqc})$.  With this convention we have the
following important definition:

\begin{df}{(\cite[Definition $3.3.1$]{t1})}
Let $X$ be a pointed and connected simplicial set in $\mathbb{U}$. The
schematization of $X$ over $k$ is
a pointed schematic homotopy type $\sch{X}{k}$,
together with a morphism 
\[
u : X \longrightarrow \sch{X}{k}
\]
in $\op{Ho}(\sff{SPr}_{*}(k))$
which is a universal for morphisms from $X$ to pointed schematic
homotopy types (in the category $\op{Ho}(\sff{SPr}_{*}(k))$.
\end{df}

\

\smallskip

We have stated the above definition only for simplicial sets in order
to simplify the exposition. However, by using the {\em singular
functor} $\op{Sing}$, attaching to each topological space $T$ the
simplicial set of singular chains in $T$ (see e.g. \cite{ho} for
details), one can define the schematization of a pointed connected
topological space.  In what follows we will always assume implicitly
that the functor $\op{Sing}$ has been applied when necessary and we
will generally not distinguish between topological spaces and
simplicial sets when considering the schematization functor.

\

Finally, recall the main existence theorem.

\begin{thm}{(\cite[Theorem  3.3.4]{t1})}
Any pointed and connected simplicial set $(X,x)$ in $\mathbb{U}$
possesses a schematization over $k$.  Furthermore, for any $i>0$ the
sheaf $\pi_{i}(\sch{X}{k},x)$ is represented by an affine group
scheme, which is commutative and unipotent for $i>1$.
\end{thm}

\

Let $(X,x)$ be a pointed connected simplicial set in $\mathbb{U}$, and let
$\sch{X}{k}$ be its schematization. Then, one has:

\begin{enumerate}

\item The affine group scheme $\pi_{1}(\sch{X}{k},x)$
is naturally isomorphic to the pro-algebraic completion
of the discrete group $\pi_{1}(X,x)$ over $k$.

\item Let $V$ be a local system of finite dimensional $k$-vector spaces
  on $X$. In particular $V$  corresponds to a linear representation of
  $\pi_{1}(\sch{X}{k},x)$ and gives rise to a local system
  $\mathcal{V}$ on $\sch{X}{k}$. Then there is a natural isomorphism
\[
H^{\bullet}(X,V)\simeq H^{\bullet}(\sch{X}{k},\mathcal{V}),
\]

\item If $X$ is simply connected and of finite type (i.e. the homotopy
type of a simply connected finite $CW$ complex), then for any
$i>1$, the group scheme $\pi_{i}(\sch{X}{k},x)$ is naturally
isomorphic to the pro-unipotent completion of the discrete groups
$\pi_{i}(X,x)$. In other words, for any $i>1$
\[
\begin{split}
\pi_{i}(\sch{X}{k},x) & \simeq
\pi_{i}(X,x)\otimes_{\mathbb{Z}}\mathbb{G}_{a} \qquad \text{if }
\op{char}(k)=0,  \\
\pi_{i}(\sch{X}{k},x) & \simeq
\pi_{i}(X,x)\otimes_{\mathbb{Z}}\mathbb{Z}_{p} \qquad \text{if }
\op{char}(k)=p>0.
\end{split}
\]
Here the groups $\pi_{i}(X,x)$ appearing in the right hand side are
thought of as constant group schemes over $k$.
\end{enumerate}

\section{Algebraic models} \label{sec:models}

In this section we discuss two different algebraic models for pointed
schematic homotopy types: cosimplicial commutative Hopf algebras and equivariant
cosimplicial commutative algebras. These two models give rise to two
different explicit formulas for the schematization of a space, having
their own advantages and disavantages according to the situation. Our
first model will allow us to complete the proof of
\cite[Theorem~3.2.9]{t1}.

We will start by introducing an intermediate model category structure on the
category of simplicial affine $k$-group schemes. This model category 
structure will then be localized in order to get the right homotopy category of 
of cosimplicial Hopf algebras suited for the setting of schematic
homotopy type. We will present this intermediate model category in a first separate section
as we think that it might have some independant interests. 

\subsection{Simplicial affine group schemes}

By a \textit{Hopf algebra} we will mean a \textit{unital and co-unital
commutative Hopf $k$-algebra}. The category of Hopf algebras will be
denoted by $\sff{Hopf}$, which is therefore equivalent to the opposite of
the category $\sff{GAff}$ of affine $k$-group schemes.  Recall that every
Hopf algebra is equal to the colimit of its Hopf subalgebras of finite
type (see \cite[III, \S 3, No. 7]{dg}). In particular the category $\sff{Hopf}$ is the category of ind-objects
in the category of of Hopf algebras of finite type and thus is
complete and co-complete \cite[Expos\'e~1,Proposition~8.9.1(b)]{sga4.1}.
We consider the category of cosimplicial Hopf algebras
$\sff{Hopf}^{\Delta}$, dual to the category of simplicial affine group
schemes $\sff{sGAff}$. When we will need to specify universes we will write
$\sff{Hopf}_{\mathbb{U}}$, $\sff{Hopf}^{\Delta}_{\mathbb{U}}$,
$\sff{GAff}_{\mathbb{U}}$, $\sff{sGAff}_{\mathbb{U}}$ \dots.

The category $\sff{GAff}_{\mathbb{U}}$ has all $\mathbb{U}$-small
limits and colimits.  In particular, the category of simplicial
objects $\sff{sGAff}_{\mathbb{U}}$ is naturally endowed with tensor
and co-tensor 
structures over the category $\sff{SSet}_{\mathbb{U}}$ of
$\mathbb{U}$-small simplicial sets. By duality, the category
$\sff{Hopf}^{\Delta}_{\mathbb{U}}$ has also a natural tensor and
co-tensor structure over $\sff{SSet}_{\mathbb{U}}$. For $X \in
\sff{SSet}_{\mathbb{U}}$ and $B_{*} \in \sff{Hopf}^{\Delta}$,
$G_{*}=\op{Spec}\, B_{*}$, we will use the standard notations
\cite[Theorem~2.5]{gj}:
\[
X\otimes B_{*}, \qquad B_{*}^{X}, \qquad X\otimes G_{*}, \qquad
G_{*}^{X}.
\] 
Explicitly for a simplicial
affine group scheme $G_{*}$ and a simplicial set $X$ we have 
$(X\otimes G)_{n} = \coprod_{X_{n}} G_{n}$, where the
coproduct is taken in the category of affine group schemes. 
To describe $G^{X}$ we
first define 
$(G^{X})_{0}$ as the equalizer 
\[
\xymatrix@1{
{\displaystyle (G^{X})_{0}} \ar[r] & {\displaystyle \prod_{n} G_{n}^{X_{n}}}
  \ar@<.5ex>[r]^-{a}\ar@<-.5ex>[r]_-{b} & {\displaystyle \prod_{p
      \stackrel{u}{\to} q} 
  G_{q}^{X^{p}}.}
}
\]
Here $a$ is the composition $\prod_{n} G_{n}^{X_{n}} \twoheadrightarrow
 G_{p}^{X_{p}} \to G_{q}^{X_{p}}$ where the second map is induced from
 $u$. Similarly $b$ is the composition $\prod_{n} G_{n}^{X_{n}}
 \twoheadrightarrow 
 G_{q}^{X_{q}} \to G_{q}^{X_{p}}$ where the second map is induced from
 $u$. 
With this definition we can now set $(G^{X})_{n} = (G^{X\times
 \Delta_{n}})_{0}$. The definitions of  $X\otimes B_{*}$ and
 $B_{*}^{X}$ are analogous. Note that
\[
\op{Spec}\, (X\otimes B_{*})\simeq G_{*}^{X} \qquad \text{and} \qquad
\op{Spec}\, (B_{*}^{X}) \simeq X\otimes B_{*}.
\]

For any simplicial set $K$ and any simplicial affine group scheme $G_{*}$, we will use the following notation
$$Map(K,G_ {*}):=(G_{*}^{K})_{0}.$$
Note that $Map(\Delta^{n},G_{*})\simeq G_{n}$.

Let  $n\geq 0$. We consider the simplicial sphere $S^{n}:=\partial \Delta^{n+1}$, 
pointed by the vertex $0\in \Delta^{n+1}$. There is a natural morphism
of affine group schemes
$$Map(S^{n},G_{*}) \longrightarrow Map(*,G_{*})\simeq G_{0}.$$
The kernel of this morphism will be denoted by $Map_{*}(S^{n},G_{*})$. 
In the same way, the kernel of the morphism
$$Map(\Delta^{n+1},G_{*}) \longrightarrow Map(*,G_{*})=G_{0}$$
will be denoted by $Map_{*}(\Delta^{n+1},G_{*})$. The inclusion
$S^{n} \subset \Delta^{n+1}$ induces a morphism
$$Map_{*}(\Delta^{n+1},G_{*}) \longrightarrow Map_{*}(S^{n},G_{*}).$$
The cokernel of this morphism, taken in $\sff{sGAff}$, will be denoted by 
$\Pi_{n}(G_{*})$.

For any $n$ and any $0\leq k\leq n$, we denote by 
$\Lambda^{n,k}$ the $k$-th horn of $\Delta^{n}$, which by definition obtained from
$\partial \Delta^{n}$ by removing its $k$-th face. For a morphism $G_{*} \longrightarrow H_{*}$
in $\sff{sGAff}$ we then have a morphism of affine group schemes
$$Map(\Delta^{n},G_{*}) \longrightarrow Map(\Lambda^{n,k},G_{*})\times_{Map(\Lambda^{n,k},G_{*})}
Map(\Delta^{n},H_{*}),$$
or equivalently 
$$G_{n} \longrightarrow Map(\Lambda^{n,k},G_{*})\times_{Map(\Lambda^{n,k},G_{*})}H_{n}.$$

\begin{df}\label{dd1}
Let $f : G_{*} \longrightarrow H_{*}$ be a morphism in $\sff{sGAff}$. 
\begin{enumerate}
\item The morphism $f$ is an \emph{equivalence} if for all $n$  the induced morphism
$$\Pi_{n}(G_{*}) \longrightarrow \Pi_{n}(H_{*})$$
is an isomorphism. 
\item The morphism $f$ is a \emph{fibration} if for all $n$ and all $0\leq k \leq n$
the induced morphism
$$G_{n} \longrightarrow Map(\Lambda^{n,k},G_{*})\times_{Map(\Lambda^{n,k},G_{*})}H_{n}$$
is a faithfully flat morphism of schemes.
\item The morphism $f$ is a \emph{cofibration} if it has the left lifting property
with respect to all morphisms which are fibrations and equivalences.
\end{enumerate}
\end{df}

It is useful to notice here that a morphism of affine group schemes
$G\longrightarrow H$ is faithfully flat if and only if the induced morphism
of Hopf algebras $\mathcal{O}(H) \longrightarrow \mathcal{O}(G)$ is injective (see
\cite[III, \S 3, No. 7]{dg}).

\begin{thm}\label{tt1}
The above definition makes $\sff{sGAff}_{\mathbb{U}}$ into a 
model category. 
\end{thm}

\textit{Proof:} We are going to apply the theorem \cite[Thm. 2.1.19]{ho}
to the opposite category $\sff{sGAff}_{\mathbb{U}}^{op}=\sff{Hopf}^{\Delta}_{\mathbb{U}}$. 
For this, we let $I$ be a set of representative of all morphisms in $\sff{sGAff}_{\mathbb{U}}$
which are fibrations $G_{*} \longrightarrow H_{*}$ such that for any $n\geq 0$ 
the group schemes $G_{n}$ and $H_{n}$ are of finite type (as schemes over $Spec\, k$).
In the same way, we let $J\subset $ be the subset corresponding to morphisms
in $I$ which are also equivalences. We need to prove that the
category $\sff{sGAff}_{\mathbb{U}}^{op}=\sff{Hopf}^{\Delta}_{\mathbb{U}}$, the class
$W$ of equivalences, and the sets $I$ and $J$ satisfy the following seven conditions.

\begin{enumerate}

\item The subcategory $W$ has the two out of the three property and is closed
under retracts.
\item The domains and codomains of $I$ are small (relative to $I$-cell).
\item The domains and codomains of $J$ are small (relative to $J$-cell).
\item $J-cell \subset W\cap I-cof$.
\item $I-inj\subset W\cap J-inj$.
\item $W\cap I-cell \subset J-cof$.
\end{enumerate}

These properties will prove the existence of a model category structure
on $\sff{Hopf}^{\Delta}_{\mathbb{U}}$ whose equivalences are the one
defined in \ref{dd1}, and whose cofibrations are generated by the set $I$. 
Before going futher in the proof of the properties above we check that the cofibrations generated by the set 
$I$ are precisely the one given in definition \ref{dd1} (note that 
by duality fibrations in $\sff{Hopf}^{\Delta}_{\mathbb{U}}$ correspond to 
cofibrations in $\sff{sGAff}_{\mathbb{U}}$ and conversly). 

\begin{lem}\label{ll0}
\begin{enumerate}
\item A morphism $G_{*} \longrightarrow H_{*}$ in $\sff{sGAff}_{\mathbb{U}}$ such that 
for each $n$, the morphism $G_{n} \longrightarrow H_{n}$ is faithfully flat is a fibration. 
\item 
A morphism in $\sff{sGAff}$ has the left lifting property with respect to $I$ if and only if
it has the left lifting property with respect every fibration.
\end{enumerate}
\end{lem}

\textit{Proof of lemma \ref{ll0}:} $(1)$ Let $K_{*}$ be the kernel of the morphism
$G_{*} \longrightarrow H_{*}$. Let $0\leq k\leq n$ and set 
$$L_{n}:=H_{n}\times_{Map(\Lambda^{n,k},H_{*})}Map(\Lambda^{n,k},G_{*}).$$
We have a commutative diagram of affine group schemes
$$\xymatrix{
K_{n} \ar[r] \ar[d] & G_{n} \ar[r] \ar[d] & H_{n} \ar[d]^-{id} \\
Map(\Lambda^{n,k},K_{*}) \ar[r] & L_{n} \ar[r] & H_{n}.}$$
Using a version of the five lemma we see that it is enough to prove that 
$H_{n} \longrightarrow Map(\Lambda^{n,k},K_{*})$ is faithfully flat. In other words, 
we can assume that $H_{*}=\{e\}$.

The simplicial presheaf $h_{G_{*}}$ represented by $G_{*}$ is a 
presheaf in simplicial groups on the site of all affine schemes. It is 
a globally fibrant simplicial presheaf (see \cite[Thm. 17.1]{may}).
As the functor $G_{*} \mapsto h_{G_{*}}$ commutes with exponential by simplicial sets (because
it commutes with arbitrary  limits) we see that this implies that for $0\leq n\leq n$ 
the morphism
$$G_{n} \longrightarrow Map(\Lambda^{n,k},G_{*})$$
induces a surjective morphism on the associated presheaves. As 
$Map(\Lambda^{n,k},G_{*})$ is an affine scheme this implies that this morphism
has in fact a section, and thus is faithfully flat by \cite[III. \S 3, No. 7]{dg}. \\

$(2)$ Let $G_{*} \longrightarrow H_{*}$ be a morphism having the left lifting property with respect 
to $I$, and let $K_{*} \longrightarrow L_{*}$ be a fibration together with a commutative
diagram
$$\xymatrix{
G_{*} \ar[r] \ar[d] & K_{*} \ar[d] \\
H_{*} \ar[r] & L_{*}.}$$
We consider the factorisation $K_{*} \longrightarrow K_{*}' \longrightarrow L_{*}$
into an faithfully flat followed by an injective morphism. The morphism 
$K'_{*} \longrightarrow L_{*}$ stays a fibration, as shown by the following 
commutative square with faithfully flat rows
$$\xymatrix{
K_{n} \ar[r] \ar[d] & K'_{n} \ar[d] \\
Map(\Lambda^{n,k},K_{*})\times_{Map(\Lambda^{n,k},L_{*})}L_{n} \ar[r] & 
Map(\Lambda^{n,k},K'_{*})\times_{Map(\Lambda^{n,k},L_{*})}L_{n}.}$$
Moreover, the morphism $K_{*} \longrightarrow K_{*}'$ also is a fibration because of part $(1)$ of the lemma
\ref{ll0}.
This implies that we are reduced to treat two cases, either $K_{*} \longrightarrow L_{*}$
is faithfully flat, or it is injective. 

We start to treat the case where $K_{*} \longrightarrow L_{*}$ is injective. In this case, the induced
morphism of simplicial sheaves
$$h_{K_{*}} \longrightarrow h_{L_{*}}$$
is a monomorphism and a local fibration. The local lifting property
with respect to the inclusion $*\hookrightarrow \Delta^{1}$ and using that the morphism is mono
implies that 
the induced morphism 
$\pi_{0}(h_{K_{*}}) \longrightarrow \pi_{0}(h_{L_{*}})$
is a monomorphism of sheaves. As monomorphisms and local fibrations are stable
by exponentiation by a finite simplicial set we also see that the induced morphism
$\pi_{i}(h_{K_{*}}) \longrightarrow \pi_{i}(h_{L_{*}})$ is a monomorphism of sheaves for 
all $i\geq 0$. Finally, the local lifting property 
for the inclusion $*\hookrightarrow \Delta^{n}$ and the fact the morphism is mono also implies that 
the induced morphism 
$\pi_{n}(h_{K_{*}}) \longrightarrow \pi_{n}(h_{L_{*}})$ is surjective for all $n>0$. 
In other words, 
the following square
$$\xymatrix{h_{K_{*}} \ar[r] \ar[d] & h_{L_{*}} \ar[d] \\
\pi_{0}(h_{K_{*}}) \ar[r] & \pi_{0}(h_{L_{*}}) }$$
is cartesian.
Equivalently, the following square 
$$\xymatrix{K_{*} \ar[r] \ar[d] & L_{*} \ar[d] \\
\Pi_{0}(K_{*}) \ar[r] & \Pi_{0}(L_{*}) }$$
is also cartesian. Therefore, in order to prove the existence of a lifting
$H_{*} \longrightarrow K_{*}$ we can replace $K_{*}$ by $\Pi_{0}(K_{*}) $
and $L_{*}$ by $\Pi_{0}(L_{*})$. We are therefore reduced to prove the following 
fact about affine group schemes: if $p : G\longrightarrow H$ is a morphism of affine group 
schemes having the left lifting property with respect to every morphism between affine
group schemes of finite type, then $p$ is an isomorphism. This last assertion follows easily 
from the fact that the category of affine group schemes is the category of pro-objects in the
category of affine group schemes of finite type.
 
We now assume that $K_{*} \longrightarrow L_{*}$ is faithfully flat. 
We consider a certain set $X$.  Its elements are
cosimplicial Hopf sub-algebras $A_{*}\subset \mathcal{O}(K_{*})$ such that:
\begin{itemize}
\item the morphism $\mathcal{O}(L_{*}) \longrightarrow \mathcal{O}(K_{*})$
factors throught $A_{*}$, 

\item there exists a morphism $H_{*} \longrightarrow Spec\, A_{*}$ making the diagram
$$\xymatrix{G_{*} \ar[r] \ar[d] & Spec\, A_{*} \ar[d] \\
H_{*} \ar[r] \ar[ru] & L_{*}}$$
commutative. 
\end{itemize}

This set is non-empty as the image of $\mathcal{O}(L_{*}) \hookrightarrow \mathcal{O}(K_{*})$
is an element of $X$.
We next order $X$ by the order induced by inclusion of cosimplicial Hopf sub-algebras.
The ordered set $X$  is inductive, and we let $A_{*} \subset \mathcal{O}(K_{*})$ be a maximal
element. Assume that $A_{*}\neq \mathcal{O}(K_{*})$. We chose a lift 
$H_{*} \longrightarrow Spec\, A_{*}$, and we consider the following commutative diagram
$$
\xymatrix{G_{*} \ar[r] \ar[d] & K_{*} \ar[d] \\
H_{*} \ar[r]  & Spec\, A_{*}.}$$
As $A_{*}\neq \mathcal{O}(K_{*})$, there exist a cosimplicial Hopf subalgebra
$D_{*} \subset \mathcal{O}(K_{*})$, such that for any $n$ the algebra $D_{n}$ is of finite type,
and with  $A_{*} \nsubseteq D_{*}$. We let $D'_{*}:=A_{*}\cap D_{*}$. Finally, we let 
$B_{*}$ be the cosimplicial Hopf sub-algebra of $\mathcal{O}(K_{*})$ which is generated by 
$D_{*}$ and $A_{*}$. There exists a commutative diagram
$$\xymatrix{
G_{*} \ar[r] \ar[d] & K_{*} \ar[d] \ar[r] & Spec\, B_{*} \ar[dl] \ar[r] & Spec\, D_{*} \ar[d] \\
H_{*} \ar[r] & Spec\, A_{*} \ar[rr] & & Spec\, D_{*}',}$$
and the square on the right hand side is furthermore cartesian. Finally, 
as $D_{*}' \longrightarrow D_{*}$ is injective, $Spec\, D_{*} \longrightarrow Spec\, D_{*}'$
is a fibration (because of part $(1)$ of lemma \ref{ll0}) between simplicial affine group schemes of finite type.
By assumption a lift $H_{*} \longrightarrow Spec\, B_{*}$ exists. But this contredicts the 
maximality of $A_{*}$. Therefore, $A_{*}=\mathcal{O}(K_{*})$ and a lift
$H_{*} \longrightarrow K_{*}$ exists.
\hfill $\Box$ \\

Let us now check that the properties $(1)$ to $(5)$ above are satisfied.
The property $(1)$ is true as equivalences are the morphisms inducing 
isomorphisms on $\Pi_{n}$ (and isomorphisms of group schemes do
have the two out of the three property and are closed under retracts). The points $(2)$ and $(3)$ 
are true as any object in the category of Hopf algebras is small (with respect to some cardinal 
depending of this object) with respect to the whole
category, therefore the domains and codomains of $I$ and $J$ are also small with respect
to the whole category $\sff{Hopf}^{\Delta}_{\mathbb{U}}$. \\

Let us now prove $(4)$. 
We have $J\subset I\subset I-cof$, and therefore $J-cell \subset I-cof$ because 
$I-cof$ is stable by push-outs and transfinite compositions. In order to prove that 
$J-cell \subset W$ it is enough to prove the following two properties:
\begin{enumerate}
\item[(a)] The trivial fibrations in $\sff{sGAff}_{\mathbb{U}}$ are stable by base-change.
\item[(b)] The trivial fibrations in $\sff{sGAff}_{\mathbb{U}}$ are stable by 
(${\mathbb{U}}$-small) filtered limits.
\end{enumerate}

Let 
$$\xymatrix{
G_{*}'\ar[r] \ar[d]_-{f'} & G_{*} \ar[d]^-{f} \\
H_{*}' \ar[r] & H_{*}}$$
be a cartesian diagram in $\sff{sGAff}_{\mathbb{U}}$. For any simplicial set $K$, 
the diagram 
$$\xymatrix{
Map(K,G_{*}')\ar[r] \ar[d] & Map(K,G_{*}) \ar[d]^-{f} \\
Map(K,H_{*}') \ar[r] & Map(K,H_{*})}$$
is a cartesian diagram of affine group schemes. Therefore, for
any $n$ and $0\leq k\leq n$ we have a
cartesian diagram of affine group schemes
$$\xymatrix{
G_{n}' \ar[r] \ar[d] & G_{n} \ar[d] \\
Map(\Lambda^{n,k},G'_{*})\times_{Map(\Lambda^{n,k},G'_{*})}H'_{n} \ar[r] & 
Map(\Lambda^{n,k},G_{*})\times_{Map(\Lambda^{n,k},G_{*})}H_{n}.}$$
As the faithfully flat morphisms are stable by base change this implies that 
$f'$ is a fibration if $f$ is so. 

We consider the functor $h : \sff{sGAff}_{\mathbb{U}} \longrightarrow \sff{SPr}_{*}(k)$, 
sending a simplicial affine group scheme $G_{*}$ to the simplicial presheaf
$X\mapsto Hom(X,G_{*})$, pointed at the unit of $G_{0}$. As the functor sending 
an affine group scheme to the sheaf of groups it represents commutes with
limits and quotients (see \cite[III, \S 3, No. 7]{dg}), we see that the homotopy sheaves 
of $\pi_{n}(h_ {G_{*}})$ are representable by the group schemes $\Pi_{n}(G_{*})$. Moreover, 
a morphism $f : G_{*} \longrightarrow H_{*}$ is a fibration our sense if and only if
the induced morphism $h_{G_{*}} \longrightarrow h_{H_{*}}$ is a local fibration in the
sense of \cite[\S 1]{ja} (i.e. satisfies the local right lifting property with respcect to 
$\Lambda^{n,k}\subset \Delta^{n}$). In particular, this morphism
of simplicial presheaves provides a long exact sequence on homotopy sheaves when $f$ is a fibration 
(see \cite[Lemma 1.15]{ja}). 
We deduce from this that if $K_{*}$ denotes the kernel of the morphism $f$, and if $f$ is a fibration, 
then there exists a long exact sequence of affine group schemes
$$\xymatrix{
\dots \Pi_{n+1}(H_{*}) \ar[r] & \Pi_{n}(K_{*}) \ar[r] & \Pi_{n}(G_{*}) \ar[r] & \Pi_{n}(H_{*}) \ar[r] & 
\Pi_{n-1}(K_{*}) \ar[r] & \dots}.$$
Using this fact we deduce that a fibration in $\sff{sGAff}_{\mathbb{U}}$ is 
an equivalence if and only if its kernel $K_{*}$ is acyclic (i.e; $\Pi_{n}(K_{*})=0$). As kernels
are stable by base change, we see that this implies that trivial fibrations are 
also stable by base change. This proves the property $(a)$ above. 

In order to prove the property $(b)$, let $G_{*}=lim_{\alpha}G_{*}^{(\alpha)}$ be 
a filtered limit of objects in $\sff{sGAff}_{\mathbb{U}}$.

\begin{lem}\label{ll1} 
For any $n\geq 0$, the natural morphism 
$$\Pi_{n}(G_{*}) \longrightarrow lim_{\alpha}\Pi_{n}(G_{*}^{\alpha})$$
is an isomorphism. 
\end{lem}

\textit{Proof of the lemma \ref{ll1}:}
For any $n$, we have a cocartesian 
square of affine group schemes
$$
\xymatrix{
Map_{*}(\partial\Delta^{n+1},G_{*}) \ar[r] \ar[d] & Map_{*}(\Delta^{n+1},G_{*}) \ar[d] \\
\{e\} \ar[r] & \Pi_{n}(G_{*}).}$$
As the functor $Map(K,G_{*})$ commutes with limits, we are reduced to show that 
filtered limits of affine group schemes preserves cocartesian squares. Let 
$$\xymatrix{
E^{(\alpha)} \ar[r] \ar[d] & F^{(\alpha)} \ar[d] \\
G^{(\alpha)} \ar[r] & H^{(\alpha)}}$$
be a filtered diagram of cocartesian squares in $\sff{GAff}$, et let 
$$\xymatrix{
E \ar[r] \ar[d] & F \ar[d] \\
G \ar[r] & H}$$
be the limit diagram. We need to prove that for any $G_{0} \in \sff{sGAff}$ 
the induced diagram of sets
$$\xymatrix{
Hom(H,G_{0}) \ar[r] \ar[d] & Hom(F,G_{0}) \ar[d] \\
Hom(G,G_{0}) \ar[r] & Hom(E,G_{0})}$$
is cartesian. As any affine group scheme $G_{0}$ is the projective limit 
of its quotients of finite type we can assume that $G_{0}$ is of finite type. 
But then, $Hom(-,G_{0})$ sends filtered limits to filtered colimits, as the
Hopf algebra corresponding to $G_{0}$ is of finite type as an algebra and thus
is a compact object in the category of Hopf algebras (i.e. $Hom(\mathcal{O}(G_{0}),-)$
commutes with filtered colimits). Therefore, the above diagram of sets is isomorphic to 
$$\xymatrix{
colim_{\alpha}Hom(H^{(\alpha)},G_{0}) \ar[r] \ar[d] & colim_{\alpha}Hom(F^{(\alpha)},G_{0}) \ar[d] \\
colim_{\alpha}Hom(G^{(\alpha)},G_{0}) \ar[r] & colim_{\alpha}Hom(E^{(\alpha)},G_{0}).}$$
This last diagram is cartesian because filtered colimits in sets preserve 
finite limits.  \hfill $\Box$  \\

The previous lemma implies that equivalences in $\sff{sGAff}_{\mathbb{U}}$ are stable
by filtered limits, and in particular by transfinite composition on the left. 
Moreover, faithfully flat morphisms of affine group schemes are also 
stable by filtered limits, because injective morphisms of Hopf algebras are
stable by filtered colimits. By definition of fibrations this implies that 
the fibrations are also stable by filtered limits. This finishes the proof of the
point $(b)$ above. The property $(4)$ 
is thus proven. \\

Before proving the last two properties $(5)$ and $(6)$ we need some 
additional lemmas.

\begin{lem}\label{ll2}
A morphism $G_{*} \longrightarrow H_{*}$ in $\sff{sGAff}_{\mathbb{U}}$
is a fibration if and only if for any trivial cofibration $A\subset B$ of $\mathbb{U}$-small simplicial sets
the natural morphism
$$Map(B,G_{*}) \longrightarrow Map(A,G_{*})\times_{Map(A,H_{*})}Map(B,H_{*})$$
is faithfully flat.
\end{lem}

\textit{Proof of the lemma \ref{ll2}:} This follows from the definition, from the fact that
the morphisms $\Lambda^{n,k} \subset \Delta^{n}$ generate the trivial cofibrations 
of simplicial sets (see \cite[Definition 3.2.1]{ho}), and from the fact that faithfully flat morphisms 
of affine group schemes are stable by pull-backs and filtered compositions (because 
injective morphisms of Hopf algebras are stable by filtered colimits).
\hfill $\Box$ \\

\begin{lem}\label{ll3}
Let $G_{*} \in \sff{sGAff}_{\mathbb{U}}$ and $A\subset B$ be a cofibration 
of $\mathbb{U}$-small simplicial sets. Then 
morphism 
$G_{*}^{B} \longrightarrow G_{*}^{A}$ is 
a fibration.
\end{lem}

\textit{Proof of the lemma \ref{ll3}:} This follows from the formula
$$Map(K,G_{*}^{A})\simeq Map(A\times K,G_{*}),$$
from the lemma \ref{ll2}, and from the fact that for any $n$ and any $k$ the natural
morphism
$$(A\times \Delta^{n})\coprod_{A\times \Lambda^{n,k}} B\times \Lambda^{n,k}
\longrightarrow B\times \Delta^{n}$$
is a trivial cofibration of simplicial sets.
\hfill $\Box$ \\

We have $J\subset I$ and thus $I-inj\subset J-inj$. In order to prove $(5)$ we need to
prove that $I-inj\subset W$. Let $i : G_{*} \longrightarrow H_{*}$ be a morphism in 
$\sff{sGAff}_{\mathbb{U}}$ having the 
left lifting property with respect to $I$. By the lemma \ref{ll0} we know that it
also has the left lifting property with respect to all fibrations. Moreover, lemma \ref{ll0} $(1)$ implies that 
$G_{*}$ itself is 
fibrant, and the lifting property for the diagram
$$\xymatrix{
G_{*} \ar[r]\ar[d] & G_{*} \ar[d] \\
H_{*} \ar[r] & \{e\}}$$
implies the existence of $r : H_{*} \longrightarrow G_{*}$ such that $r\circ i=id$.
The lemma \ref{ll3} implies that $H^{\Delta^{1}}_{*} \longrightarrow  H^{\partial \Delta^{1}}_{*}=H^{\Delta^{1}}\times H^{\Delta^{1}}$
is fibration. Therefore, the lifting property for the diagram
$$\xymatrix{
G_{*} \ar[r]^-{i}\ar[d] & H^{\Delta^{1}}_{*} \ar[d] \\
H_{*} \ar[r]_-{id,i\circ r} & H_{*}\times H_{*},}$$
implies the existence of a morphism
$$h : H_{*} \longrightarrow H^{\Delta^{1}}_{*}$$
which is a homotopy between the identity and $i\circ r$. This implies that 
$i : h_{G_{*}} \longrightarrow h_{H_{*}}$ is a homotopy equivalence of simplicial presheaves, and thus 
induces isomorphisms on homotopy sheaves. As we have already seen the homotopy sheaves of 
$h_{G_{*}}$ and $h_{H_{*}}$ are representable by the homotopy group schemes of $G_{*}$ and $H_{*}$. Therefore, 
$i$ itself is an equivalence in $\sff{sGAff}_{\mathbb{U}}$. \\

It only remains to prove that the property $(5)$ is satisfied. But as we have already seen during the proof of point $(4)$
the morphism in $I-cell$ are fibrations in $\sff{sGAff}$ (because fibrations are stable by base change and filtered limits).
In particular any morphism in $I-cell\cap W$ is a trivial fibration. Because of lemma \ref{ll0} this implies that 
$I-cell\cap W\subset J-cof$. This finishes the proof of the theorem \ref{tt1}.
\hfill $\Box$ \\

\noindent
Before continuing further we need to make an important 
remark concerning the dependance of the universe $\mathbb{U}$
for the model category  $\sff{Hopf}_{\mathbb{U}}^{\Delta}$.
If $\mathbb{U} \in \mathbb{V}$ are two universes, one gets a natural
inclusion functor $\sff{Hopf}_{\mathbb{U}}^{\Delta} \longrightarrow
\sff{Hopf}^{\Delta}_{\mathbb{V}}$.  An important fact about this inclusion is that
the sets $I$ and $J$ defined in the proof of the previous proposition are
independant of the choice of the universe. Therefore, the sets $I$ and $J$ are
generating set of cofibrations and trivial cofibrations for both model categories 
$\sff{Hopf}_{\mathbb{U}}^{\Delta}$ and 
$\sff{Hopf}^{\Delta}_{\mathbb{V}}$. As we know that 
a morphism is a cofibration 
if and only if it is a retract of a relative $I$-cell complex (see \cite[Def. 2.1.9]{ho}
and \cite[Prop. 2.1.18 (b)]{ho}), we
see that a morphism in $\sff{Hopf}_{\mathbb{U}}^{\Delta}$ is 
a cofibration if and only if it is a cofibration in 
$\sff{Hopf}^{\Delta}_{\mathbb{V}}$. In the same way, fibrations are the morphisms
having the right lifting property with respect to $J$, and thus
a morphism in $\sff{Hopf}_{\mathbb{U}}^{\Delta}$ is 
a fibration if and only if it is a fibration in 
$\sff{Hopf}^{\Delta}_{\mathbb{V}}$. Moreover, using the functorial factorization
of cosimplicial algebras via the small object argument with respect
to $I$ and $J$ (as described in the proof above) one sees that the
full sub-category $\sff{Hopf}_{\mathbb{U}}^{\Delta}$ is stable by the
functorial factorizations in $\sff{Hopf}^{\Delta}_{\mathbb{V}}$.  In particular, 
a path object (see \cite[Def. 1.2.4 (2)]{ho}) in $\sff{Hopf}_{\mathbb{U}}^{\Delta}$ is also a path
object in $\sff{Hopf}^{\Delta}_{\mathbb{V}}$. As a consequence, we see that 
two morphisms in $\sff{Hopf}_{\mathbb{U}}^{\Delta}$ between fibrant and
cofibrant objects are homotopic (see \cite[Def. 1.2.4 (5)]{ho}) if and only if they are homotopic
in $\sff{Hopf}^{\Delta}_{\mathbb{V}}$.

A consequence of these remarks and of the fact that 
the homotopy category of a model category is equivalent to the
category of fibrant and cofibrant objects and homotopy classes of objects between them
(see \cite[Thm. 1.2.10]{ho}) is that the induced functor
\[
\op{Ho}(\sff{Hopf}^{\Delta}_{\mathbb{U}}) \longrightarrow
\op{Ho}(\sff{Hopf}^{\Delta}_{\mathbb{V}}) 
\]
is fully faithful. The essential image of this functor consists of all
cosimplicial Hopf algebras in $\mathbb{V}$ which are equivalent to 
a cosimplicial Hopf algebra in $\mathbb{U}$. 

The conclusion of this short discussion is that increasing the size of the ambient universe
does not destroy the homotopy theory of cosimplicial Hopf algebras, and its
only effect is to add more objects to the corresponding homotopy category.
\\

An important additional property of the model category $\sff{sGAff}_{\mathbb{U}}$ is its right properness.

\begin{cor}\label{cc1}
The model category $\sff{sGAff}_{\mathbb{U}}$ is right proper (i.e. 
equivalences are stable by  pull-backs along fibrations).
\end{cor}

\textit{Proof:} This follows from lemma \ref{ll0} $(1)$, as it implies in particular that 
any object in $\sff{sGAff}_{\mathbb{U}}$ is fibrant (see \cite[Corollary 13.1.3]{hir}).
This could also have been checked directly using the long exact sequence in homotopy for fibrations (see
the proof of point $(4)$ of the theorem \ref{tt1}).  \hfill $\Box$ \\

Another interesting property is the compatibility between the model structure and
the simplicial enrichement.

\begin{cor}\label{cc2}
Together with its natural simplicial enrichement, the model category $\sff{sGAff}_{\mathbb{U}}$
is a simplicial model category. 
\end{cor}

\textit{Proof:} This is a consequence of lemma \ref{ll2}. Indeed, we need to verify the axiom
M7 of \cite[Definition 9.1.6]{hir}.  For this, let $i : A\subset B$ be a cofibration of $\mathbb{U}$-small simplicial sets,
and $f : G_{*} \longrightarrow H_{*}$ a fibration in $\sff{sGAff}_{\mathbb{U}}$. The fact that 
$$G_{*}^{B} \longrightarrow G_{*}^{A}\times_{H_{*}^{A}}H_{*}^{B}$$
is a again fibration follows from the definition, from lemma \ref{ll2} and from the fact that the
morphism
$$B\times \Lambda^{n,k}\coprod_{A\times \Lambda^{n,k}}B\times \Delta^{n} \longrightarrow
B\times \Delta^{n}$$
is a trivial cofibration. 

If moreover, $i$ or $f$ is an equivalence then, because of \cite[Corollary 1.5]{ja}
and \cite[Lemma 1.1.5]{ja}, the morphism
$$h_{G_{*}} \longrightarrow h_{G_{*}^{A}\times_{H_{*}^{A}}H_{*}^{B}}\simeq
h_{G_{*}}^{A}\times_{h_{H_{*}}^{A}}h_{H_{*}}^{B}$$
is a local equivalence. This implies that 
$$G_{*}^{B} \longrightarrow G_{*}^{A}\times_{H_{*}^{A}}H_{*}^{B}$$
is an equivalence. \hfill $\Box$ \\

The simplicial Homs of the simplicial category $\sff{SGAff}_{\mathbb{U}}$ 
will be denoted by $\underline{Hom}$. By proposition \ref{cc2} these simplicial
Homs possess derived version (see \cite[\S 4.3]{ho})
$$\mathbb{R}\underline{Hom}(-,-) : 
\opi{Ho}(\sff{sGAff}_{\mathbb{U}})^{op} \times \opi{Ho}(\sff{sGAff}_{\mathbb{U}})
\longrightarrow \opi{Ho}(\sff{SSet}_{\mathbb{U}}).$$

\subsection{The P-localization}

Let $K$ be an affine group scheme of finite type over $k$, and $V$ be a finite
dimensional linear representation of $K$. We consider $V$ as an affine group scheme
(its group law being the addition), as well as the simplicial affine
group scheme $K(V,n)$ for $n\geq 0$. By definition $K(V,n)$ is the classifying 
space of $K(V,n-1)$, as defined for instance in \cite[\S 1.3]{t1}, and we set 
$K(V,0)=V$. The group scheme $K$ acts on $V$ and therefore 
acts on $K(V,n)$. We will denote by $K(K,V,n)$ the simplicial 
group scheme which is the semi-direct product of $K$ (considered
as a constant simplicial object) by $K(V,n)$. We therefore have 
a split exact sequence of simplicial affine group schemes
$$\xymatrix{1 \ar[r] & K(V,n) \ar[r] & K(K,V,n) \ar[r] & K \ar[r] & 1.}$$

Ideally, we would like now to construct the left Bousfield localization of $\sff{sGAff}_{\mathbb{U}}$ with 
respect to the set of objects 
$K(K,V,n)$, for all $K$, $V$ and $n$. Dually, this would correspond to 
perform a right Bousfield localization of $\sff{Hopf}^{\Delta}_{\mathbb{U}}$
with respect to the corresponding set of objects. The only general result insuring the
existence of a right Bousfield localization  we are aware about  
is the theorem \cite[Theorem 5.1.1]{hir} which requires the model category 
to be cellular. Unfortunately the model category $\sff{Hopf}^{\Delta}_{\mathbb{U}}$ is not 
cellular, as cofibrations are simply not monomorphisms. It is therefore unclear that the
localized model structure exists (we think it does). In this section we will
show that the existence of a localization functor
$$(-)^{P} : \opi{Ho}(\sff{sGAff}_{\mathbb{U}}) \longrightarrow
\opi{Ho}^{P}(\sff{sGAff}_{\mathbb{U}}),$$
which will be enough to prove the equivalence between 
pointed schematic homotopy types and cosimplicial Hopf
algebras up to P-equivalences. \\

We start by  defining our new equivalences in $\sff{sGAff}_{\mathbb{U}}$. 
We will see later that they are precisely the quasi-isomorphisms of
cosimplicial Hopf algebras (see Corollary \ref{cc2'}).

\begin{df}\label{d1}
A morphism $f : G_{*} \longrightarrow H_{*}$ is a \emph{P-equivalence} if for 
any affine group scheme of finite type $K$, any linear representation of finite dimensional
$V$ of $K$, and any $n\geq 1$ the induced morphism
$$f^{*} : \mathbb{R}\underline{Hom}(H_{*},K(K,V,n)) \longrightarrow 
\mathbb{R}\underline{Hom}(G_{*},K(K,V,n))$$
is an isomorphism in $\opi{Ho}(\sff{SSet}_{\mathbb{U}})$.
\end{df}

The definition above gives the new class of equivalences
on $\sff{sGAff}_{\mathbb{U}}$. The localization of the category
$\sff{sGAff}_{\mathbb{U}}$ with respect to 
$P$-equivalences will be denoted by 
$$\opi{Ho}^{P}(\sff{sGAff}_{\mathbb{U}}).$$
As an equivalence is a P-equivalence we have a natural functor
$$l : \opi{Ho}(\sff{sGAff}_{\mathbb{U}}) \longrightarrow \opi{Ho}^{P}(\sff{sGAff}_{\mathbb{U}}).$$

\begin{prop}\label{pp2}
The above functor
$$l : \opi{Ho}(\sff{sGAff}_{\mathbb{U}}) \longrightarrow \opi{Ho}^{P}(\sff{sGAff}_{\mathbb{U}})$$
possesses a right adjoint
$$j : \opi{Ho}^{P}(\sff{sGAff}_{\mathbb{U}}) \longrightarrow \opi{Ho}(\sff{sGAff}_{\mathbb{U}})$$
which is fully faithful. The essential image of $j$ is of the smallest
full sub-category of $\opi{Ho}(\sff{sGAff}_{\mathbb{U}})$ containing the objects
$K(K,V,n)$ and which is stable by homotopy limits. 
\end{prop}

\textit{Proof:} This an application of the existence of a cellularization functor
applied to the model category $\sff{Hopf}^{\Delta}_{\mathbb{U}}$
and to set of objects $\{K(K,V,n)\}$ (see \cite[Proposition 5.2.3, 5.2.4]{hir}). For the convenience of the
reader we reproduce the argument here. 

We let $X$ be a set of representative for the morphisms
$$K(K,V,n)^{\Delta^{m}} \longrightarrow K(K,V,n)^{\partial \Delta^{m}},$$
for all affine group scheme of finite type $K$, all finite dimensional linear 
representation $V$ of $K$, and all integers $n\geq 1$ and $m\geq 0$.
Because of Corollary \ref{cc2} all the morphisms in $X$ are fibrations.
For a given cofibrant object $G_{*}\in \sff{sGAff}_{\mathbb{U}}$  
we construct a tower of cofibrant objects in $G_{*}/\sff{sGAff}_{\mathbb{U}}$
\[
\xymatrix{G_{*} \ar[r] &  \dots G_{*}^{(i)} \ar[r] & G_{*}^{(i-1)} \ar[r] & \dots
  \ar[r] & G_{*}^{(0)}=*,} 
\]
defined inductively in the following way. We let $I_{i}$ be the
set of all commutative squares in $\sff{sGAff}_{\mathbb{U}}$
\[
\xymatrix{
G_{*} \ar[r] \ar[d] & H_{*} \ar[d]^-{u} \\
G_{*}^{(i-1)} \ar[r] & H'_{*},}
\]
where $u \in X$. The set $I_{i}$ is then $\mathbb{U}$-small.

We now define an object $G_{*}(i)$ by the pull-back square
\[
\xymatrix{
G_{*}(i) \ar[r] \ar[d] & \prod_{j\in I_{i}}H_{*} \ar[d] \\
G_{*}^{(i-1)} \ar[r] & \prod_{j\in I_{i}}H_{*}'.}
\]
We let 
$G_{*}\rightarrow G_{*}^{(i)}=Q(G_{*}(i))$ be the cofibrant replacement of $F
\rightarrow G_{*}(i)$ in $G_{*}/\sff{sGAff}_{\mathbb{U}}$.
This defines the tower inductively on $i$.  Finally, we consider the morphism
\[
\alpha : G_{*} \longrightarrow \widetilde{G_{*}}:= \op{Lim}_{i}G_{*}^{(i)}.
\]
We first claim that the object $\widetilde{G_{*}}$ 
is P-local, in the sense that 
for any P-equivalence $H_{*} \longrightarrow H_{*}'$ the induced morphism
$$\mathbb{R}\underline{Hom}(H_{*}',G_{*}) \longrightarrow  \mathbb{R}\underline{Hom}(H_{*},G_{*})$$
is an isomorphism in $\opi{Ho}(\sff{SSet})$.
Indeed, by construction it is a $\mathbb{U}$-small homotopy limit
of P-local objects.  

It then only remains to
see that the morphism $\alpha$ is a $P$-equivalence, as this would imply 
formally that $G_{*} \longrightarrow \widetilde{G_{*}}$ is a P-localization
(i.e. un universal P-equivalence with a P-local object). This in turn would imply the result as 
the functor $G_{*} \mapsto \widetilde{G_{*}}$ would then identify the
localization $\opi{Ho}^{P}(\sff{sGAff}_{\mathbb{U}})$ with the full sub-category 
of P-local objects in $\opi{Ho}(\sff{sGAff}_{\mathbb{U}})$. Moreover, the construction of
$\widetilde{G_{*}}$ shows that the P-local objects are obtained by successive
homotopy limits of objects of the form $K(K,V,n)$. 

Let us then  consider $K(K,V,n)$ for a given 
affine group scheme of finite type $K$, a finite dimensional linear 
representation $V$ of $K$ and an integer $n\geq 1$. 
Using that the simplicial affine group scheme 
$K(K,V,n)$ is levelwise of finite type and is $n$-truncated we see that 
it is a $\omega$-cosmall object in $\sff{sGAff}_{\mathbb{U}}$. 
Moreover, as $\widetilde{G_{*}}$, $G_{*}^{(i)}$ and $G_{*}$ are all cofibrant
 the morphism 
\[
\alpha^{*} : \mathbb{R}\underline{\op{Hom}}(\widetilde{G_{*}},K(K,V,n)) 
\longrightarrow
\mathbb{R}\underline{\op{Hom}}(\widetilde{G_{*}},K(K,V,n))
\]
is isomorphic in $\op{Ho}(\sff{SSet})$ to the natural morphism
\[
\alpha^{*} : \op{Colim}_{i}\underline{\op{Hom}}(G_{*}^{(i)},K(K,V,n))
\longrightarrow 
\underline{\op{Hom}}(G_{*},K(K,V,n)).
\]
Thus, by the inductive construction of the tower, we deduce from this that
for any $m$, 
the morphism 
\[
\mathbb{R}\underline{\op{Hom}}(\widetilde{G_{*}},K(K,V,n))^{\partial
  \Delta^{m}} 
\longrightarrow
\mathbb{R}\underline{\op{Hom}}(\widetilde{G_{*}},K(K,V,n))^{\Delta^{m}}
\times^{h}_{\mathbb{R}\underline{\op{Hom}} 
(G_{*},K(K,V,n))^{\Delta^{m}}}
\mathbb{R}\underline{\op{Hom}}(G_{*},K(K,V,n))^{\partial \Delta^{m}}
\] 
is
surjective on connected components. This implies that $\alpha^{*}$
is an isomorphism, and therefore that $G_{*} \longrightarrow \widetilde{G_{*}}$
is a $P$-equivalence as required.
\hfill $\Box$ \\

Using the previous proposition, we will always implicitely identify
the category $\opi{Ho}^{P}(\sff{sGAff}_{\mathbb{U}})$ with 
the full sub-category of $\opi{Ho}(\sff{sGAff}_{\mathbb{U}})$
consisting of P-local objects, and also with the smallest full sub-category 
of $\opi{Ho}(\sff{sGAff}_{\mathbb{U}})$ containing 
the $K(K,V,n)$ and which is stable by homotopy limits.
The left adjoint of the inclusion functor
$$\opi{Ho}^{P}(\sff{sGAff}_{\mathbb{U}}) \hookrightarrow 
\opi{Ho}(\sff{sGAff}_{\mathbb{U}})$$
is isomorphic to the localization functor, and can be identified with
the construction $G_{*} \mapsto \widetilde{G_{*}}$ given during the proof of the
theorem. We will more often denote this functor by 
$$G_{*} \mapsto G_{*}^{P}.$$

We will give now a more explicit description of the
P-equivalences related to Hochschild cohomology of affine group schemes
with coefficients in linear representations. 
For a Hopf algebra $B$ and a $B$-comodule $V$  (most of the time assumed to be 
of finite dimension but this is not needed),  one
can consider the cosimplicial $k$-vector space
\[
\begin{array}{cccc}
C^{*}(B,V) : & \Delta & \longrightarrow & \mbox{$k$-\sff{Vect}} \\
 & [n] & \mapsto & C^{n}(B,V):=V\otimes B^{\otimes n},
\end{array}
\]
where the transitions morphisms $V\otimes B^{\otimes n}
\longrightarrow V\otimes B^{\otimes m}$ are given by the 
co-action and co-unit morphisms. From a dual point of view, $V$
corresponds to a linear representation of the affine group scheme $G=Spec\, B$,
and can be considered as a quasi-coherent sheaf $\mathbb{V}$ on the
simplicial affine scheme $BG$.  The cosimplicial space $C^{*}(B,V)$
is by definition the cosimplicial space of sections
$\Gamma(BG,\mathbb{V})$ of this sheaf on $BG$.

The cosimplicial vector space $C^{*}(B,V)$ has an associated total
complex, whose cohomology groups will be denoted by
\[
H^{i}(B,V):=H^{i}(\op{Tot}(C^{*}(B,V))).
\]
These are the Hochschild cohomology groups of $B$ with coefficients in
$V$.  From a dual point of view, the complex $\op{Tot}(C^{*}(B,V))$
also computes the cohomology of the affine group scheme $G$ with
coefficients in the linear representation $V$. We will also use the
notations

$$C^{*}(G,V):=C^{*}(B,V) \qquad 
H^{*}(G,V):=H^{*}(B,V). $$

Let now $B_{*} \in \sff{Hopf}_{\mathbb{U}}^{\Delta}$ be a cosimplicial Hopf
algebra, and let us consider the Hopf algebra $H^{0}(B_{*})$ of $0$-th
cohomology of $B_{*}$. We have 
$Spec\, H^{0}(B_{*})\simeq \Pi_{0}(G_{*})$, where $G_{*}:=\op{Spec}\, B_{*}$.  
By construction $H^{0}(B_{*})$ is the limit
(in the category of Hopf algebras) of the cosimplicial diagram $[n]
\mapsto B_{n}$, and therefore comes equipped with a natural
co-augmentation $H^{0}(B_{*}) \rightarrow B_{*}$.  In particular, if
$V$ is any $H^{0}(B_{*})$-comodule, $V$ can also be considered
naturally as comodule over each $B_{n}$. We get in this way a
cosimplicial object in the category of cosimplicial vector spaces
(i.e. a bi-cosimplicial vector space)
\[
[n] \mapsto C^{*}(B_{n},V).
\]
We will denote the diagonal associated to this bi-cosimplicial space  by
\[
C^{*}(B_{*},V):= \op{Diag}\left(
[n] \mapsto C^{*}(B_{n},V)\right):=\left([n] \mapsto C^{n}(B_{n},V)\right).
\]
If we denote by $G_{*}=\op{Spec}\, B_{*}$ the associated simplicial affine
group scheme, we will also use the notation
\[
C^{*}(G_{*},V):=C^{*}(B_{*},V).
\]
The cohomology of the total complex associated to $C^{*}(B_{*},V)$
is called the \textit{Hochschild cohomology} of the cosimplicial
algebra $B_{*}$ (or equivalently of the simplicial affine group scheme
$G_{*}=\op{Spec}\, B_{*}$)  
with coefficients in the comodule $V$, and is denoted by
\[
H^{n}(B_{*},V):=H^{n}(\op{Tot}(C^{*}(B_{*},V))) \qquad
H^{n}(G_{*},V):=H^{n}(\op{Tot}(C^{*}(G_{*},V))).
\]

We see that $C^{*}(G_{*},V)$ is the cosimplicial space of sections of $V$ considered as 
a quasi-coherent sheaf on $h_{BG_{*}}$, represented by the simplicial affine scheme $BG_{*}$.
As the sheaf $V$ is quasi-coherent we have natural isomorphisms
$$H^{i}(Tot(C^{*}(G_{*},V)))\simeq H^{i}(h_{BG_{*}},V):=\pi_{0}(Map_{SPr(k)/h_{BG_{*}}}(h_{BG_{*}},F(V,i))),$$
where $F(V,i) \longrightarrow h_{BG_{*}}$ is the relative Eleinberg-MacLane construction 
on the sheaf of abelian groups $V$. This can be easily deduced from the
special case of a non-simplicial affine group scheme treated in \cite[\S 1.3, \S 1.5]{t1}, simply  by
noticing that $h_{BG_{*}}$ is naturally equivalent to the homotopy colimit of the
$h_{BG_{n}}$ when $n$ varies in $\Delta^{op}$. 

\begin{prop}\label{ll4}
A morphism $f : G_{*} \longrightarrow H_{*}$ is a P-equivalence if and only if 
it satisfies the following two properties.
\begin{enumerate}
\item For any finite dimensional linear representation $V$ of $\Pi_{0}(H_{*})$ the induced morphism
$$f^{*} : H^{*}(H_{*},V) \longrightarrow H^{*}(G_{*},V)$$
is an isomorphism.
\item The induced morphism $\Pi_{0}(G_{*}) \longrightarrow \Pi_{0}(H_{*})$ is an isomorphism.
\end{enumerate}
\end{prop}

\textit{Proof:} Let $G_{*} \in \sff{sGAff}_{\mathbb{U}}$, $K$ an 
affine group scheme of finite type, and $V$ a finite dimensional linear representation of $K$.
We assume that $G_{*}$ is a cofibrant object. Then, there exists a natural
morphism of simplicial sets
$$\underline{Hom}(G_{*},K(K,V,n)) \longrightarrow \underline{Hom}_{*}(Bh_{G_{*}},Bh_{K(K,V,n)}),$$
where $\underline{Hom}_{*}$ denotes the simplicial set of morphisms of the
category $\sff{SPr}_{*}(k)$. Composing with the fibrant and cofibrant
replacement functors in $\sff{SPr}_{*}(k)$ we get this way a natural morphism 
in $\opi{Ho}(\sff{SSet_{*}})$
$$\mathbb{R}\underline{Hom}(G_{*},K(K,V,n)) \longrightarrow 
\mathbb{R}\underline{Hom}_{*}(Bh_{G_{*}},Bh_{K(K,V,n)}).$$
This morphism comes equiped with a natural projection to the set 
$Hom(\Pi_{0}(G_{*}),K)$
$$\xymatrix{
\mathbb{R}\underline{Hom}(G_{*},K(K,V,n)) \ar[rr] \ar[rd] & & 
\mathbb{R}\underline{Hom}_{*}(Bh_{G_{*}},Bh_{K(K,V,n)}) \ar[dl] \\
 & Hom(\Pi_{0}(G_{*}),K). &}$$
Moreover, the homotopy fiber $F$ of the right hand side morphism 
is such that 
$$\pi_{i}(F)\simeq H^{n-i}(G_{*},V).$$
This shows that in order to prove the lemma it is enough to prove that 
the horizontal morphism (recall that $G_{*}$ is cofibrant, and that 
$K(K,V,n)$ is always fibrant)
$$\underline{Hom}(G_{*},K(K,V,n))\simeq \mathbb{R}\underline{Hom}(G_{*},K(K,V,n)) \longrightarrow 
\mathbb{R}\underline{Hom}_{*}(Bh_{G_{*}},Bh_{K(K,V,n)})$$
is an isomorphism.

For this, we use the fact that $\underline{Hom}(G_{*},K(K,V,n))$
is naturally isomorphic to the total space (see \cite[Definition 18.6.3]{hir}) of the
cosimplicial space $m \mapsto \underline{Hom}(G_{m},K(K,V,n))$
$$\mathbb{R}\underline{Hom}(G_{*},K(K,V,n))\simeq
\underline{Hom}(G_{*},K(K,V,n)) \simeq Tot(m \mapsto \underline{Hom}(G_{m},K(K,V,n))).$$
In the same way, the simplicial presheaf $Bh_{G_{*}}$ is equivalent to 
the homotopy colimit of the diagram $n \mapsto Bh_{G_{n}}$, and we thus have
$$\mathbb{R}\underline{Hom}_{*}(Bh_{G_{*}},Bh_{K(K,V,n)})\simeq
Holim_{n}\mathbb{R}\underline{Hom}_{*}(Bh_{G_{m}},Bh_{K(K,V,n)}).$$
Therefore, in order to prove the lemma we need to prove the following two 
separate statements. 
\begin{enumerate}
\item[(a)] The natural morphism (see \cite[Definition 18.7.3]{hir})
$$Tot \left( m \mapsto \underline{Hom}(G_{m},K(K,V,n)) \right) \longrightarrow
Holim \left( m \mapsto \underline{Hom}(G_{m},K(K,V,n)) \right)$$
is an isomorphism in $\opi{Ho}(\sff{SSet})$.

\item[(b)] For any $m$, the natural morphism
$$\underline{Hom}(G_{m},K(K,V,n)) \longrightarrow \mathbb{R}\underline{Hom}_{*}(Bh_{G_{m}},
Bh_{K(K,V,n)})$$
is an isomorphism in $\opi{Ho}(\sff{SSet})$.
\end{enumerate}

For property $(a)$ we use that $K(K,V,n)$ are abelian group objects in 
the category of simplicial affine group schemes over $K$. This implies that the morphism
$$Tot \left( m \mapsto \underline{Hom}(G_{m},K(K,V,n)) \right) \longrightarrow
Holim \left( m \mapsto \underline{Hom}(G_{m},K(K,V,n)) \right)$$
is a morphism of abelian group objects in the category of simplicial sets over the set $Hom(G_{m};K)$, 
or in other words is a morphism between disjoint union of simplicial abelian groups. 
The property $(a)$ above then follows from the fact that for any cosimplicial 
object $X_{*}$ in the category of simplicial abelian groups, the natural morphism
$$Tot(X_{*}) \longrightarrow Holim (X_{*})$$
is a weak equivalence (see e.g. \cite[X \S 4.9, XI \S 4.4]{bk}). 

The property $(b)$ is clear for $n=0$, as we simply find a bijection of sets
$$Hom(G_{m},K\rtimes V) \simeq Hom_{*}(Bh_{G_{m}},Bh_{K\rtimes V}).$$
Let us assume that $n>0$. Because we have 
$$\pi_{i}(\underline{Hom}(G_{m},K(K,V,n)))\simeq \pi_{i-1}(\underline{Hom}(G_{m},K(K,V,n-1)))$$
$$\pi_{i}(\mathbb{R}\underline{Hom}_{*}(Bh_{G_{m}},
Bh_{K(K,V,n)})) \simeq \pi_{i-1}(\mathbb{R}\underline{Hom}_{*}(Bh_{G_{m}},
Bh_{K(K,V,n-1)}))$$
it is enough by induction to prove that the natural morphism
$$\pi_{0}(\underline{Hom}(G_{m},K(K,V,n))) \longrightarrow \pi_{0}(\mathbb{R}\underline{Hom}_{*}(Bh_{G_{m}},
Bh_{K(K,V,n)}))$$
is bijective. To prove this it is enough to prove that both morphisms
$$\pi_{0}(\underline{Hom}(G_{m},K(K,V,n))) \longrightarrow Hom(G_{m},K)$$
$$\pi_{0}(\mathbb{R}\underline{Hom}_{*}(Bh_{G_{m}},
Bh_{K(K,V,n)})) \longrightarrow Hom(G_{m},K)$$
are isomorphisms. As the natural projection $K(K,V,n) \longrightarrow K$ posses
a section these two morphisms are surjective. It remains to show that these morphisms
are injective. This in turns will follow from the following lemma.

\begin{lem}\label{sl}
Let $G_{*}$ be a cofibrant simplicial affine group scheme. Then, for any $m\geq 0$ and
any faithfully flat morphism of affine group schemes $H \longrightarrow H'$ the induced
morphism
$$Hom(G_{m},H') \longrightarrow Hom(G_{m},H)$$
is surjective.
\end{lem}

\textit{Proof of the lemma:} We consider the evaluation functor
$H_{*} \mapsto H_{m}$, and its right adjoint
$$i_{m}^{*} : \sff{GAff}_{\mathbb{U}} \longrightarrow \sff{sGAff}_{\mathbb{U}}.$$
We have
$$i_{m}^{*}(H)_{p}\simeq H^{Hom([m],[p])},$$
showing that $i_{m}^{*}$ preserves faithfully flat morphisms. Moreover, we have for any 
simplicial set $A$
$$Map(A,i_{m}^{*}(H))\simeq H^{A_{m}},$$
which easily implies that $\Pi_{i}(i_{m}^{*}(H))=0$ for any $i$. In particular, if $H\longrightarrow H'$
is a faithfully flat morphism we see from lemma \ref{ll0} $(12)$ that the induced morphism
$$i_{m}^{*}(H) \longrightarrow i_{m}^{*}(H')$$
is a trivial fibration. The right lifting property for $G_{*}$ with respect to 
this last morphism then precisely says that  
$Hom(G_{m},H') \longrightarrow Hom(G_{m},H)$ is surjective.
\hfill $\Box$ \\

The previous lemma implies that 
for any linear representation $V$ of $G_{m}$ we have 
$H^{2}(G_{m},V)=0$, as this group classifies extensions
$$\xymatrix{1 \ar[r] & V \ar[r] & J \ar[r] & G_{m} \ar[r] & 1.}$$
This in turns implies that $H^{i}(G_{m},V)$ for any $i>1$, and thus that the natural
morphism
$$\pi_{0}(\mathbb{R}\underline{Hom}_{*}(Bh_{G_{m}},
Bh_{K(K,V,n)})) \longrightarrow Hom(G_{m},K)$$
is injective, as its fibers are in bijections with $H^{n}(G_{m},V)$. 
For the other morphism, we consider the morphism of affine group schemes
$$Map(\Delta^{1},K(K,V,n)) \longrightarrow 
Map(\Delta^{1},K) \times_{Map(\partial\Delta^{1},K)}Map(\partial \Delta^{1},K(K,V,n)).$$
As the morphism $K(K,V,n) \longrightarrow K$ is a fibration (because of lemma \ref{ll0} $(1)$)
and is relatively $(1)$-connected (i.e. its fibers are $1$-connected), the above morphism
is faithfully flat. The right lifting property of $G_{m}$ with respect to this morphism (which 
is insured by the sub-lemma \ref{sl})
implies that the morphism
$$\pi_{0}(\underline{Hom}(G_{m},K(K,V,n))) \longrightarrow Hom(G_{m},K)$$
is injective. This finishes the proof of the proposition \ref{ll4}.
\hfill $\Box$ \\

From the proof of the proposition \ref{ll4} we also extract the following important 
corollary that will be used in the sequel.

\begin{cor}\label{cll4}
Let $G_ {*}\in \sff{sGAff}_{\mathbb{U}}$, $K$ be a affine group scheme of finite type, 
$V$ a finite dimensional linear representation of $K$, and $n\geq 1$. Then the natural
morphism
$$\mathbb{R}\underline{Hom}(G_{*},K(K,V,n)) \longrightarrow 
\mathbb{R}\underline{Hom}_{\sff{SGp(k)}}(h_{G_{*}},h_{K(K,V,n)}) \simeq
\mathbb{R}\underline{Hom}_{\sff{SPr_{*}(k)}}(Bh_{G_{*}},Bh_{K(K,V,n)})$$
is an isomorphism in $\opi{Ho}(\sff{SSet})$.
\end{cor}

\subsection{Cosimplicial Hopf algebras and schematic homotopty types}

We are now ready to explain how cosimplicial Hopf algebras are
models for schematic homotopy types. 

\bigskip

We now consider $\sff{SGp}(k)$, the category
of presheaves of $\mathbb{V}$-simplicial groups on $(\sff{Aff},\op{fpqc})$.
It will be endowed with the model category structure for which
equivalences and fibrations are defined via the forgetful
functor $\sff{SGp}(k) \longrightarrow \sff{SPr}(k)$: a morphism in 
$\sff{SGp}(k)$ is a fibration and/or an equivalence if and only if
it is so as a morphism in $\sff{SPr}(k)$ (by forgetting the group structure).
We consider the Yoneda functor
$$h_{-} : \sff{sGAff}_{\mathbb{U}} \longrightarrow \sff{SGp}(k).$$

The functor
$h$ sends equivalences to local equivalences of simplicial presheaves and thus 
defines above possesses a functor
$$h : \opi{Ho}(\sff{sGAff}_{\mathbb{U}}) \longrightarrow
\opi{Ho}(\sff(SGp)(k)).$$

Consider the classifying space functor
\[
\bB : \sff{SGp}(k) \longrightarrow \sff{SPr}(k)_{*},
\]
from the category of presheaves of simplicial groups to the category of
pointed simplicial presheaves. It is well known (see for example
\cite[Theorem~1.4.3]{t1} and \cite[Proposition~1.5]{t2}) that this
functor preserves equivalences and induces a fully faithful functor on
the homotopy categories 
\[
\bB : \op{Ho}(\sff{SGp}(k)) \longrightarrow \op{Ho}(\sff{SPr}_{*}(k)).
\]
Composing with the  functor
\[
h\, :
\op{Ho}(\sff{sGAff}_{\mathbb{U}}) \longrightarrow 
\op{Ho}(\sff{SGp}(k))
\]
one gets a functor
\[
\bB h\, :
\op{Ho}(\sff{sGAff}_{\mathbb{U}}) \longrightarrow
\op{Ho}(\sff{SPr}_{*}(k)). 
\]

With this notations we now have:

\begin{thm}\label{t1}
The functor
\[
\bB h\, :
\op{Ho}(\sff{sGAff}_{\mathbb{U}}) \longrightarrow
\op{Ho}(\sff{SPr}_{*}(k)) 
\]
is fully faithful when restricted to the full sub-category 
$\op{Ho}^{P}(\sff{sGAff}_{\mathbb{U}}) \subset 
\op{Ho}(\sff{sGAff}_{\mathbb{U}})$ consisting of
P-local objects. 
Its essential image consists precisely of all psht.
\end{thm}

{\bf Proof:} We first analyze the full faithfulness properties of
$\bB h$. 
Let $G_{*}$ and $G_{*}'$ be two P-local and cofibrant simplicial affine group schemes.
We need to prove that the natural morphism of
simplicial sets
\[
\underline{\op{Hom}}(G_{*},G'_{*}) \to
\underline{\op{Hom}}_{\sff{SPr}_{*}(k)}(\bB h_{G'},\bB h_{G'_{*}})\to 
\mathbb{R}\underline{\op{Hom}}_{\sff{SPr}_{*}(k)}(\bB h_{G'},\bB h_{G'_{*}})
\]
is a weak equivalence. As the functor $B$ is fully faithful, it is enough to show that 
the induced morphism
\[
\underline{\op{Hom}}(G_{*},G'_{*}) \to
\mathbb{R}\underline{\op{Hom}}_{\sff{SGp}(k)}(h_{G'},h_{G'_{*}})
\]
is an equivalence.
As $G_{*}'$ is P-local, it can be written
as a transfinite composition of homotopy pull-backs of objects of the form $K(K,V,n)$
(see \cite[Theorem 5.1.5]{hir}, and also the proof of proposition \ref{pp2}), for some
affine group scheme of finite type $K$, some
finite dimensional linear representation $V$ of $K$ and some integer $n\geq 1$. 

\begin{lem}\label{ll5}
The functor
$$h : \sff{sGAff}_{\mathbb{U}} \longrightarrow \sff{SGp}(k)$$
commutes with homotopy limits of P-local objects.
\end{lem}

\textit{Proof of the lemma \ref{ll5}:} It is enough to show that 
$h$ preserves homotopy pull-backs and homotopy products (possibly infinite) of P-local objects. 
The case of homotopy pull-backs follows from \cite[Lemma 1.15]{ja} as we have already seen
that the functor $h$ sends fibrations to local fibrations of simplicial presheaves. 
The case of infinite homotopy product of P-local objects would follow from the fact that 
$h$ commutes with products and sends P-local objects to fibrant objects in 
$\sff{SGp}(k)$. As the P-local objects are obtained by  transfinite composition of
homotopy pull-backs of objects of the form $K(K,V,n)$ it is enough to check that 
$h_{K(K,V,n)}$ is a fibrant simplicial presheaf. But, as a simplicial presheaf it is
isomorphic to $h_{K}\times K(\mathbb{G}_{a},n)^{d}$, where $d$ is the dimension
of $V$, which is known to be fibrant as both $h_{K}$ and 
$K(\mathbb{G}_{a},n)$ are fibrant (see e.g. \cite[Lemma 1.1.2]{t1}). 
\hfill $\Box$ \\

Using the lemma \ref{ll5} we are reduced to prove that the natural morphism
$$
\underline{\op{Hom}}(G_{*},K(K,V,n)) \longrightarrow 
\mathbb{R}\underline{\op{Hom}}_{\sff{SPr}_{*}(k)}(\bB h_{G'},\bB h_{K(K,V,n)})$$
is an equivalence. But this is something we have already seen during the
proof of proposition \ref{ll4} (see corollary \ref{cll4}). 

It remains now to prove that the image of the functor
$\bB h$ consists precisely of all psht. The fact that
the image of $\bB h$ is contained in the category of
psht follows from \cite[Corollary~3.2.7]{t1} and the fact that 
the P-local affine group schemes are generated by 
homotopy limits by the objects $K(K,V,n)$ (note that 
$\bB h_{K(K,V,n)}$ is a psht with $\pi_{1}=K$,
$\pi_{n}=V$ and $\pi_{i}=0$ for $i\neq n$). Now suppose that $F$ is
a psht. We will show that $F$ is $\boldsymbol{P}$-equivalent to an
object of the form $Bh\, G_{*}$ for some P-local object $G_{*}\in
\sff{sGAff}_{\mathbb{U}}$ (see \cite[Definition~3.1.1]{t1} for
the notion of $\boldsymbol{P}$-local equivalences). Since by
definition a psht is a $\boldsymbol{P}$-local object, this implies
that $F$ is isomorphic to $Bh\, G_{*}$ as required.
For this, we assume that $F$ is cofibrant as an object in
$\sff{SPr}_{*}(k)$, and we construct a $\boldsymbol{P}$-local model of
$F$ by the small object argument. Let $\mycal{K}$ be a
$\mathbb{U}$-small set of representatives of all psht of the form
$\bB h_{K(K,V,n)}$ for some affine group scheme $K$ which is of finite
type as a $k$-algebra, some finite dimensional linear representation $V$ of $K$,  and
for some $n\geq 1$.  We consider the $\mathbb{U}$-small set $K$ of all
morphisms of the form
\[
G^{\Delta^{m}} \longrightarrow G^{\partial \Delta^{m}},
\]
for $G \in \mycal{K}$ and $m\geq 0$. Finally, we let $\Lambda(K)$ be the
$\mathbb{U}$-small set of morphisms in $\sff{SPr}_{*}(k)$ which are
fibrant approximations of the morphisms of $K$.
We construct a tower of cofibrant objects in $F/\sff{SPr}_{*}(k)$
\[
\xymatrix{F \ar[r] &  \dots F_{i} \ar[r] & F_{i-1} \ar[r] & \dots
  \ar[r] & F_{0}=*,} 
\]
defined inductively in the following way. We let $I_{i}$ be the
set of all commutative squares in $\sff{SPr}_{*}(k)$
\[
\xymatrix{
F \ar[r] \ar[d] & G_{1} \ar[d]^-{u} \\
F_{i-1} \ar[r] & G_{2},}
\]
where $u \in \Lambda(K)$. Let $J_{i}$ be a sub-set in $I_{i}$ of
representative for the isomorphism classes of objects in the homotopy
category of commutative squares in $\sff{SPr}_{*}(k)$ (i.e. $J_{i}$ is
a subset of representative for the equivalent classes of squares in
$I_{i}$). Note that the set $J_{i}$ is $\mathbb{U}$-small since all
stacks in the previous diagrams are affine $\infty$-gerbes.

We now define an object $F_{i}'$ as the pull-back square
\[
\xymatrix{
F_{i}' \ar[r] \ar[d] & \prod_{j\in J_{i}}G_{1} \ar[d] \\
F_{i-1} \ar[r] & \prod_{j\in J_{i}}G_{2},}
\]
and finally
$F\rightarrow F_{i}=Q(F_{i}')$ as a cofibrant replacement of $F
\rightarrow F_{i}$.
This defines the $i$-step of the tower
from its $(i-1)$-step. Finally, we consider the morphism
\[
\alpha : F \longrightarrow \widetilde{F}:= \op{Lim}_{i}F_{i}.
\]
We first claim that the object $\widetilde{F}\in
\op{Ho}(\sff{SPr}_{*}(k))$ lies in the 
essential image of the functor $\bB\op{Spec}$. 
Indeed, it is a $\mathbb{U}$-small homotopy limit
of objects belonging to this essential image, and as the functor
$\bB h$ is fully faithful and commutes with homotopy
limits of P-local objects (see \ref{ll5}) we see 
that $\widetilde{F}$ stays in the essential image of
$\bB h$.  It remains to
see that the morphism $\alpha$ is a $\boldsymbol{P}$-equivalence. For this, let
$F_{0}$ be a psht of the form $\bB h_{K(K,V,n)}$. 
 Using that $\widetilde{F}$ and
$F_{0}$ are both in the image of $\bB h$, 
we see that the morphism
\[
\alpha^{*} : \mathbb{R}\underline{\op{Hom}}(\widetilde{F},F_{0}) 
\longrightarrow
\mathbb{R}\underline{\op{Hom}}(F,F_{0})
\]
is isomorphic in $\op{Ho}(\sff{SSet})$ to the natural morphism
\[
\alpha^{*} : \op{Colim}_{i}\mathbb{R}\underline{\op{Hom}}(F_{i},F_{0})
\longrightarrow 
\mathbb{R}\underline{\op{Hom}}(F,F_{0}).
\]
Thus, by inductive construction of the tower, it is then clear that
for any $m$, 
the morphism 
\[
\mathbb{R}\underline{\op{Hom}}(\widetilde{F},F_{0})^{\partial
  \Delta^{m}} 
\longrightarrow
\mathbb{R}\underline{\op{Hom}}(\widetilde{F},F_{0})^{\Delta^{m}}
\times^{h}_{\mathbb{R}\underline{\op{Hom}} 
(F,F_{0})^{\Delta^{m}}}
\mathbb{R}\underline{\op{Hom}}(F,F_{0})^{\partial \Delta^{m}}
\] 
is
surjective on connected components. This implies that $\alpha^{*}$
is an isomorphism, and therefore that $F\longrightarrow \widetilde{F}$
is a $\boldsymbol{P}$-equivalence. \ \hfill $\Box$

\

\bigskip

\noindent
The following corollary completes the characterization of psht
given in \cite[Theorem~3.2.4,Proposition~3.2.9]{t1}.

\begin{cor}\label{c2}
An object $F \in \op{Ho}(\sff{SPr}_{*}(k))$ is a psht if and only if
$\pi_{1}(F,*)$ is an affine group scheme and $\pi_{i}(F,*)$
is an affine unipotent group scheme.
\end{cor}
{\bf Proof:} This follows from \cite[Cor. 3.2.7]{t1} and the essential
surjectivity part 
of theorem \ref{t1}. \ \hfill $\Box$ 

\ 

\bigskip

We are now ready to explain the relations between 
cosimplicial Hopf algebras and schematic homotopy types.
The first step is to identify the P-equivalences in 
$\sff{Hopf}^{\Delta}_{\mathbb{U}}$. 

\begin{cor}\label{cc2'}
A morphism of cosimplicial Hopf algebras 
$$f : B_{*} \longrightarrow B_{*}'$$ 
is a P-equivalence if and only if the induced morphism on 
the total complexes
$$Tot(B_{*}) \longrightarrow Tot(B'_{*})$$
is a quasi-isomorphism.
\end{cor}

\textit{Proof:} We let $G_{*}\in \sff{sGAff}_{\mathbb{U}}$ and
$G_{*} \longrightarrow G_{*}^{P}$ be a P-localization. Using corollary \ref{cll4} we see that the induced morphism
$$\bB h_{G_{*}} \longrightarrow \bB h_{G_{*}^{P}}$$
is a P-equivalence of simplicial presheaves. In the same way, it has been 
proved in \cite[Corollary 3.2.7]{t1} that the natural morphism
$$\bB h_{G_{*}} \longrightarrow \bB \mathbb{R}Spec\, \mathcal{O}(G_{*})$$
is a P-equivalence of simplicial presheaves. As both 
$\bB h_{G_{*}^{P}}$ and $\bB \mathbb{R}Spec\, \mathcal{O}(G_{*})$
are P-local simplicial presheaves, we deduce that they are naturally equivalent as pointed simplicial 
presheaves. 

From this and theorem \ref{t1}, we deduce that 
a morphism $G_{*} \longrightarrow G_{*}'$ of simplicial affine group schemes is 
a P-equivalence if and only if the induced morphism of simplicial presheaves
$$\bB \mathbb{R}Spec\, \mathcal{O}(G_{*}) \longrightarrow \bB \mathbb{R}Spec\, \mathcal{O}(G'_{*})$$
is an equivalence. This is also equivalent to say that 
the induced morphism
$$\mathbb{R}Spec\, \mathcal{O}(G_{*}) \longrightarrow \mathbb{R}Spec\, \mathcal{O}(G_{*}')$$
is an equivalence. By \cite[Corollary 2.2.3]{t1} this is equivalent to say that 
$$\mathcal{O}(G_{*}') \longrightarrow \mathcal{O}(G_{*})$$
is an equivalence of cosimplicial algebras, or in other words that 
$$Tot(\mathcal{O}(G_{*}')) \longrightarrow Tot(\mathcal{O}(G_{*}))$$
is a quasi-isomorphism. 
\hfill $\Box$ \\

If we put \ref{cc2'} and \ref{t1} together we find the following result, explaining in 
which sense cosimplicial Hopf algebras are algebraic models for
pointed schematic homotopy types.

\begin{cor}\label{cc2t1}
Let $B_{*} \mapsto B_{*}^{P}$ be a P-localization in $\sff{Hopf}^{\Delta}_{\mathbb{U}}$. Then, 
the composite functor
$$B_{*} \mapsto B^{P}_{*} \mapsto \bB h_{Spec\, B^{P}_{*}}$$
induces an equivalence between the
$\opi{Ho}^{P}(\sff{Hopf}^{\Delta}_{\mathbb{U}})^{op}$, the localized category 
of cosimplicial Hopf algebras along quasi-isomorphisms, and the homotopy category 
of pointed schematic homotopy types. 
\end{cor}

\noindent
Using cosimplicial Hopf algebras as models for psht we can describe
the schematization functor $X \mapsto \sch{X}{k}$ of \cite{t1} in the
following explicit way.

Let $X$ be a pointed and connected simplicial set in $\mathbb{U}$. As
$X$ is connected, one can choose an equivalent $X$ which is a reduced
simplicial set (i.e. a simplicial set with when $X_{0}=*$). We apply
Kan's loop group construction associating to $X$ a simplicial group
$GX$. Kan's loop group functor is left adjoint to the classifying
space functor associating to each simplicial group its classical
classifying space. Explicitly (see \cite[\S V]{gj}) $GX$ is the
simplicial group whose group of $n$-simplices is the free group
generated by the set $X_{n+1}-s_{0}(X_{n})$. Applying the
pro-algebraic completion functor lewelwise we turn $GX$ into a
simplicial affine group scheme $GX^{\op{alg}}$ in $\mathbb{U}$.  The
object $GX^{\op{alg}}$ can be thought of as an \textit{algebraic loop
space for $X$}.  The corresponding cosimplicial algebra $B_{*} :=
\mathcal{O}(GX^{\op{alg}})$ of $k$-valued regular functions on
$GX^{\op{alg}}$ is such that the level $n$ component $B_{n}$ of
$B_{*}$ is the sub-Hopf algebra of $k^{GX_{n}}=\op{Hom}(GX_{n},k)$ of
functions spanning a finite dimensional sub-$GX_{n}$-module. 

\begin{cor}\label{c3'}
With the above notations, there exists a natural 
isomorphism in $\op{Ho}(\sff{SPr}_{*}(k))$
\[
\sch{X}{k}\simeq \bB  \op{Spec}\, B^{P}_{*},
\]
where $B_{*}^{P}$ is a P-local model for $B_{*}$. 
\end{cor}

\textit{Proof:} As shown in the proof of \cite[Theorem 3.3.4]{t1}, the
schematization $\sch{X}{k}$ is given by 
$\bB \mathbb{R}Spec\, B_{*}$. An we have already seen during the
proof of corollary \ref{cc2'} that there exists a natural equivalence
$$\bB \mathbb{R}Spec\, B_{*} \simeq \bB  \op{Spec}\, B^{P}_{*}.$$
\hfill $\Box$ \\

\subsection{Equivariant cosimplicial algebras}

In this section we provide another algebraic model for psht. This
model is based on equivariant cosimplicial algebras and is quite similar 
to the approach of \cite{bs}. The main difference between the two 
approaches  is that we will work equivariantly over 
a affine group scheme, and not over a discrete group as this is
done in  \cite{bs}.

\bigskip

\begin{center} \textit{Equivariant stacks} \end{center}

We will start by some reminders on the notion of equivariant stacks. \\

For the duration of this section we fix a presheaf of groups $G$ on
$(\sff{Aff})_{\opi{ffqc}}$, which will be considered as a group object
in $\sff{SPr}(k)$. We will make the assumption that $G$ is
cofibrant as an object of $\sff{SPr}(k)$.
For example, $G$ could be representable (i.e. an affine group scheme),
or a constant presheaf associated to
a group in $\mathbb{U}$.

Because the direct product makes the category $\sff{SPr}(k)$ into a
cofibrantly generated symmetric monoidal
closed model category, for which the monoid axiom of
\cite[Section~3]{ss} is
satisfied, the category $\opi{G-}\sff{SPr}(k)$ of
simplicial presheaves equipped with a left action of $G$ is
again a closed model category \cite{ss}. Recall that the fibrations
(respectively equivalences) in $\opi{G-}\sff{SPr}(k)$ are defined
to be the morphisms inducing fibrations (respectively equivalences)
between the underlying
simplicial presheaves. The model category $\opi{G-}\sff{SPr}(k)$
will be called the {\em model category of $G$-equivariant simplicial
presheaves}, and the objects in $\op{Ho}(\opi{G-}\sff{SPr}(k))$
will be called {\em $G$-equivariant stacks}. For any $G$-equivariant
stacks $F$ and $F'$ we will denote by
$\underline{\op{Hom}}_{G}(F,F')$ the simplicial set of morphisms in
$\op{\opi{G-}SPr}(k)$, and  by
$\mathbb{R}\underline{\op{Hom}}_{G}(F,F')$ its derived version.

On the other hand, to any group $G$ one can associate  its classifying
simplicial presheaf $\opi{BG} \in \sff{SPr}_{*}(k)$
(see \cite[$\S 1.3$]{t1}). The
object $\opi{BG} \in \op{Ho}(\sff{SPr}_{*}(k))$ is uniquely characterized
 up to
a unique isomorphism by the properties
\[
\pi_{0}(\opi{BG})\simeq * \qquad \pi_{1}(\opi{BG},*)\simeq G \qquad
\pi_{i}(\opi{BG},*)=0 \; for \; i>1.
\]
Consider the comma category $\sff{SPr}(k)/\opi{BG}$, of
objects over the classifying simplicial presheaf $\opi{BG}$, endowed
with its natural simplicial closed model structure (see
\cite{ho}).
Note that since
the model structure on $\sff{SPr}(k)$ is
proper, the model category $\sff{SPr}(k)/\opi{BG}$ is an invariant,
up to a Quillen equivalence, of the isomorphism class
of $\opi{BG}$ in $\op{Ho}(\sff{SPr}_{*}(k))$. This implies in particular
that one is free to choose any model for $\opi{BG} \in
\op{Ho}(\sff{SPr}(k))$
when dealing with the homotopy category
$\op{Ho}(\sff{SPr}(k)/\opi{BG})$. We will set
$\opi{BG}:=\opi{EG}/G$, where $\opi{EG}$ is a cofibrant model of $*$ in
$\opi{G-}\sff{SPr}(k)$ which is fixed once for all.

We now define a pair of adjoint functors
\[
\xymatrix{
\opi{G-}\sff{SPr}(k) \ar@<1ex>[r]^{\opi{De}} &
\sff{SPr}(k)/\opi{BG},
\ar@<1ex>[l]^{\opi{Mo}}
}
\]
where $\opi{De}$ stands for \textit{descent} and $\opi{Mo}$ for
\textit{monodromy}.
If $F$ is a $G$-equivariant simplicial presheaf, then
$\opi{De}(F)$ is defined to be $(\opi{EG}\times F)/G$, where $G$ acts
diagonally on $\opi{EG}\times F$.
Note that there is a natural projection
$\opi{De}(F) \longrightarrow \opi{EG}/G=\opi{BG}$, and so
$\opi{De}(F)$ is naturally
an object in $\sff{SPr}(k)/\opi{BG}$.
This adjunction is easily seen to be a Quillen adjunction. Furthermore,
using the reasoning of \cite{t2}, one
can also  show that this Quillen adjunction is actually a Quillen
equivalence.  For future reference we state this as a lemma:

\begin{lem}\label{leqstack}
The Quillen adjunction $(\opi{De},\opi{Mo})$ is a Quillen equivalence.
\end{lem}
{\bf Proof:} The proof is essentially the same as the proof of
\cite[$2.22$]{t2} and is left to the reader. \ \hfill $\Box$

\

\bigskip

\noindent
The previous lemma implies that the derived Quillen adjunction induces
an equivalence of categories
\[
\op{Ho}(\opi{G-}\sff{SPr}(k)) \simeq \op{Ho}(\sff{SPr}(k)/\opi{BG}).
\]
For any $G$-equivariant stack $F \in \op{Ho}(\opi{G-}\sff{SPr}(\mathbb{C}))$,
we define the {\em quotient stack $[F/G]$} of $F$ by $G$ as the
object $\mathbb{L}\opi{De}(F) \in
\op{Ho}(\sff{SPr}(k)/\opi{BG})$ corresponding to $F \in
\op{Ho}(\opi{G-}\sff{SPr}(k))$.
By construction, the homotopy fiber (taken at the distinguished point
of $BG$) of the natural
projection
\[
p : [F/G] \longrightarrow \opi{BG}
\]
is canonically isomorphic to the
underlying stack of the $G$-equivariant stack $F$. 

\bigskip

\begin{center} \textit{Equivariant cosimplicial algebras and
    equivariant affine stacks} \end{center}  

Suppose that $G$ is an affine group scheme and consider the category
of $k$-linear representations of $G$.  This category will be denoted by
$\sff{Rep}(G)$. Note that it is an abelian $k$-linear tensor
category, which admits all $\mathbb{V}$-limits and
$\mathbb{V}$-colimits.  The category of cosimplicial $G$-modules is
defined to be the category $\sff{Rep}(G)^{\Delta}$, of cosimplicial
objects in $\sff{Rep}(G)$.  For any object $E \in
\sff{Rep}(G)^{\Delta}$ one can construct the normalized cochain
complex $\No(E)$ associated to $E$, which is a cochain complex in
$\sff{Rep}(G)$.  Its cohomology representations $H^{i}(\No(E)) \in
\sff{Rep}(G)$ will simply be denoted by $H^{i}(E)$. This construction
is obviously functorial and gives rise to various cohomology functors $H^{i} :
\sff{Rep}(G)^{\Delta} \longrightarrow \sff{Rep}(G)$.  As the category
$\sff{Rep}(G)$ is $\mathbb{V}$-complete and $\mathbb{V}$-co-complete,
its category of cosimplicial objects $\sff{Rep}(G)^{\Delta}$ has a
natural structure of a simplicial category (\cite[II,Example
$2.8$]{gj}).

Following the argument of \cite[II.4]{q2}, it can be seen 
that there exists a simplicial finitely generated closed model
structure on the category $\sff{Rep}(G)^{\Delta}$ with the properties:
\begin{itemize}
\item A morphism $f : E \longrightarrow E'$ is an equivalence if and only if
for any $i$, the induced morphism $H^{i}(f) : H^{i}(E) \longrightarrow
H^{i}(E')$ is an isomorphism.
\item A morphism $f : E \longrightarrow E'$ is a cofibration if and
only if for any $n>0$,
the induced morphism $f_{n} : E_{n} \longrightarrow E'_{n}$ is a monomorphism.
\item A morphism $f : E \longrightarrow E'$ is a fibration if and only
if it is an epimorphism whose
kernel $K$ is such that for any $n\geq 0$, $K_{n}$ is an injective
object in $\sff{Rep}(G)$.
\end{itemize}

\bigskip

The category $\sff{Rep}(G)^{\Delta}$ is endowed with a symmetric monoidal
structure, given by the tensor product of cosimplicial
$G$-modules (defined levelwise). In particular we can consider the
category $\opi{G-}\sff{Alg}^{\Delta}$ of commutative unital monoids in
$\sff{Rep}(G)^{\Delta}$. It is reasonable to view the objects in
$\opi{G-}\sff{Alg}^{\Delta}$ as cosimplicial
algebras equipped with an action of the group scheme $G$. Motivated by
this remark we will refer to the category $\opi{G-}\sff{Alg}^{\Delta}$ as
the {\em category of $G$-equivariant cosimplicial algebras}. From
another point of view, the category $\opi{G-}\sff{Alg}^{\Delta}$
is also the category of simplicial affine schemes in $\mathbb{V}$
equipped with an action of $G$.

Every $G$-equivariant cosimplicial algebra $A$ has an underlying
cosimplicial $G$-module again denoted by $A \in
\sff{Rep}(G)^{\Delta}$. This defines a forgetful functor
$$\opi{G-}\sff{Alg}^{\Delta} \longrightarrow \sff{Rep}(G)^{\Delta}$$
which has a left adjoint $L$, given by the free commutative
monoid construction.

\begin{prop}\label{peqmod}
There exists a simplicial cofibrantly generated closed model structure
on the category $\opi{G-}\sff{Alg}^{\Delta}$ satisfying
\begin{itemize}
\item A morphism $f : A \longrightarrow A'$ is an equivalence if and
only if the induced morphism in $\sff{Rep}(G)^{\Delta}$
is an equivalence.
\item A morphism $f : A \longrightarrow A'$ is a fibration if and only
if the induced morphism in
$\sff{Rep}(G)^{\Delta}$ is a fibration.
\end{itemize}
\end{prop}
{\bf Proof:} This again an application of the small 
object argument and more precisely of theorem \cite[Theorem 2.1.19]{ho}. 

Let $I$ and $J$ be sets of generating
cofibrations and trivial cofibrations in $\sff{Rep}(G)^{\Delta}$. That
is, $I$ is the set of monomorphisms between finite dimensional
cosimplicial $G$-modules, and $J$ the set of
trivial cofibrations between finite dimensional $G$-modules.
Consider the
forgetful functor $\opi{G-}\sff{Alg}^{\Delta} \longrightarrow
\sff{Rep}(G)^{\Delta}$, and its left adjoint
$L : \sff{Rep}(G)^{\Delta} \longrightarrow \opi{G-}\sff{Alg}^{\Delta}$. The
functor $L$ sends a cosimplicial representation $V$ of $G$ to the
free commutative $G$-equivariant cosimplicial algebra generated by
$V$. We will
apply the small object argument to the sets $L(I)$ and $L(J)$. 

By construction, the morphisms in $L(I)-inj$ are precisely the morphisms
inducing surjective quasi-isomorphisms on the associated normalized
cochain complexes. In the same way, the morphisms in 
$L(J)-inj$ are precisely the morphisms
inducing surjections on the associated normalized
cochain complexes. From this, it is easy to see that amoung the conditions $(1)$ to $(6)$
only the inclusion $J-cell  \subset W$ requires a proof as well as $(6)$. 
However, $(6)$ can also be replaced by 
$W\cap I-inj\subset I-inj$ which is easily seen by the previous descriptions. 
Therefore, it only remains
to check that for any morphism $A
\longrightarrow B$ in $J$, and
any morphism $L(A) \longrightarrow C$, the induced morphism
$C\longrightarrow C\coprod_{L(A)}L(B)$
is an equivalence. But since the
forgetful functor
$\opi{G-}\sff{Alg}^{\Delta} \longrightarrow Alg^{\Delta}$ is a left adjoint,
it commutes with
colimits and so $C\longrightarrow C\coprod_{L(A)}L(B)$ will be an
equivalence due to the fact that 
$\opi{Alg}^{\Delta}$ is endowed with a model category structure
for which $L(A) \longrightarrow L(B)$ is a trivial cofibration (see \cite[Theorem 2.1.2]{t1}).
\hfill $\Box$

\

\bigskip

\

\noindent
We will also need the following result:

\begin{lem}\label{leq1}
The forgetful functor $\opi{G-}\sff{Alg}^{\Delta} \longrightarrow
\sff{Alg}^{\Delta}$ is a left Quillen functor.
\end{lem}
{\bf Proof:} The forgetful functor has a right adjoint 
$$
F : \sff{Alg}^{\Delta} \longrightarrow \opi{G-}\sff{Alg}^{\Delta},
$$ which to a cosimplicial algebra $A \in \sff{Alg}^{\Delta}$ assigns
the $G$-equivariant algebra $F(A) := \mathcal{O}(G)\otimes A$, where
$G$ acts on $\mathcal{O}(G)$ by left translation.

By adjunction, it is enough to show that $F$
is right Quillen. By definition, the
functor $F$ preserves
equivalences. Moreover, if $A \longrightarrow B$ is a fibration of
cosimplicial algebras (i.e. an epimorphism),  then
the map $\mathcal{O}(G)\otimes A \longrightarrow \mathcal{O}(G)\otimes
B$ is again an epimorphism. The kernel of the latter
morphism is then isomorphic to $\mathcal{O}(G)\otimes K$, where
$K$ is the kernel of $A \longrightarrow B$.
But, for any vector space $V$, $\mathcal{O}(G)\otimes V$ is always an
injective object in $\sff{Rep}(G)$. This implies
that $\mathcal{O}(G)\otimes K_{n}$ is an injective object in
$\sff{Rep}(G)^{\Delta}$, and therefore $F(A) \longrightarrow F(B)$
is a fibration in $\opi{G-}\sff{Alg}^{\Delta}$ (see Proposition
\ref{peqmod}). \ \hfill $\Box$

\

\bigskip

For any $G$-equivariant cosimplicial algebra $A$, one can define its
(geometric) spectrum \linebreak
$\op{Spec}_{G}\, A \in \opi{G-}\sff{SPr}(k)$, by
taking the usual spectrum of its underlying cosimplicial algebra and
keeping track of the $G$-action. Explicitly,
if $A$ is given by a morphism of cosimplicial algebras $A
\longrightarrow A\otimes \mathcal{O}(G)$, one finds
a morphism of simplicial schemes
\[
G\times \op{Spec}\, A\simeq \op{Spec}\, (A\otimes \mathcal{O}(G))
\longrightarrow \op{Spec}\, A,
\]
which induces a well defined $G$-action on the simplicial scheme
$\op{Spec}\, A$. Hence, by passing to the simplicial presheaves
represented by $G$ and $\op{Spec}\, A$, one gets the $G$-equivariant
simplicial presheaf $\op{Spec}_{G}(A)$. This procedure defines
a functor
\[
\op{Spec}_{G}\, : (\opi{G-}\sff{Alg}^{\Delta})^{\op{op}} \longrightarrow
\opi{G-}\sff{SPr}(k),
\]
and we have the following

\begin{cor}\label{ceq1}
The functor $\op{Spec}_{G}$ is right Quillen.
\end{cor}

{\bf Proof:} Clearly $\op{Spec}_{G}$ commutes with
$\mathbb{V}$-limits. Furthermore, the category
$\opi{G-}\sff{Alg}^{\Delta}$ possesses a small set of small generators in
$\mathbb{V}$. For example, one can take a set of representatives
of $G$-equivariant cosimplicial algebras $A$ with $A_{n}$ of finite
type for any $n\geq 0$. This implies that
$(\opi{G-}\sff{Alg}^{\Delta})^{\op{op}}$ possesses a small set of small co-generators in
$\mathbb{V}$. The existence of the left adjoint to
$\op{Spec}_{G}$ then follows from the special adjoint theorem of \cite[\S V.8]{mac}.

To prove that $\op{Spec}_{G}$ is right Quillen it remains to prove
that it preserves fibrations and trivial fibrations.
In other words, one needs to show that if $A \longrightarrow A'$ is a
(trivial) cofibration of $G$-equivariant
cosimplicial algebras, then $\op{Spec}_{G}(A') \longrightarrow
\op{Spec}_{G}(A)$ is a (trivial) fibration in
$\opi{G-}\sff{SPr}(k)$. As fibrations and equivalences in
$\opi{G-}\sff{SPr}(k)$ are defined on the underlying object, this
follows immediately from the fact that the non-equivariant
$\op{Spec}$ is right Quillen and from Lemma \ref{leq1}. \hfill $\Box$

\

The left adjoint of $\op{Spec}_{G}$ will be denoted by
$\mathcal{O}_{G} : \opi{G-}\sff{SPr}(\mathbb{C}) \longrightarrow
\opi{G-}\sff{Alg}^{\Delta}$.

\

\bigskip

The previous corollary allows one to form the right derived functor of
$\op{Spec}_{G}$:
\[
\mathbb{R}\op{Spec}_{G} : \op{Ho}(\opi{G-}\sff{Alg}^{\Delta})^{\op{op}}
\longrightarrow \op{Ho}(\opi{G-}\sff{SPr}(k)),
\]
which possesses a left adjoint $\mathbb{L}\mathcal{O}_{G}$.
One can then compose this functor with the quotient stack functor
$[-/G]$, and obtain a functor
\[
[\mathbb{R}\op{Spec}_{G}(-)/G] : \op{Ho}(\opi{G-}\sff{Alg}^{\Delta})^{\op{op}}
\longrightarrow \op{Ho}(\sff{SPr}(k)/\opi{BG}),
\]
which still possesses a left adjoint due to the fact that
$[-/G]$ is an equivalence
of categories. We will denote this left adjoint again by
\[
\mathbb{L}\mathcal{O}_{G} : \op{Ho}(\sff{SPr}(k)/\opi{BG})
\longrightarrow \op{Ho}(\opi{G-}\sff{Alg}^{\Delta})^{\op{op}}.
\]

\

\begin{prop}\label{peq2}
If $A \in \op{Ho}(\opi{G-}\sff{Alg}^{\Delta})^{\op{op}}$ is isomorphic to some
$G$-equivariant cosimplicial algebra in
$\mathbb{U}$, then the adjunction morphism
\[
A \longrightarrow
\mathbb{L}\mathcal{O}_{G}(\mathbb{R}\op{Spec}_{G}\, A)
\]
is an isomorphism.

In particular, the functors $\mathbb{R}\op{Spec}_{G}$ and
$[\mathbb{R}\op{Spec}_{G}(-)/G]$ become fully faithful when
restricted to the full sub-category of
$\op{Ho}(\opi{G-}\sff{Alg}^{\Delta})$ consisting of $G$-equivariant cosimplicial
algebras isomorphic to some object in $\mathbb{U}$.
\end{prop}
{\bf Proof:} Let $A$ be a cofibrant $G$-equivariant cosimplicial
algebra in $\mathbb{U}$, and let $F:=\op{Spec}_{G}\, A$.
Since $F$ and 
$\op{hocolim}_{n\in \Delta^{\op{op}}}\op{Spec}_{G}\, A_{n}$ are weakly
equivalent in 
$\opi{G-}\sff{SPr}(k)$, and since $\mathcal{O}_{G}$ is a left
Quillen functor we get
\[
\mathbb{L}\mathcal{O}_{G}(F)\simeq \op{holim}_{n \in
\Delta}\mathbb{L}\mathcal{O}_{G}(\op{Spec}_{G}\, A_{n})
\simeq \op{holim}_{n \in \Delta}A_{n}\simeq A.
\]
The proposition is proven. \  \hfill $\Box$

\begin{df}\label{deqstack}
An equivariant stack $F \in \op{Ho}(\opi{G-}\sff{SPr}(k))$ is a
$G$-equivariant affine stack if it is isomorphic to some
$\mathbb{R}\op{Spec}_{G}(A)$, with $A$ a $G$-equivariant
cosimplicial algebra in $\mathbb{U}$.
\end{df}

\

\bigskip

We conclude this section with a proposition showing
that stacks of the form $[\mathbb{R}\op{Spec}_{G}\, (A)/G]$
are often pointed schematic homotopy types.

\begin{prop}\label{peqpsht}
Let $A \in \opi{G-}\sff{Alg}^{\Delta}$ be a $G$-equivariant cosimpicial
algebra in $\mathbb{U}$, such that the underlying algebra
of $A$ has an  augmentation $x : A \longrightarrow k$. Assume also
that $A$  is 
cohomologically connected, i.e.  that $H^{0}(A)\simeq k$. 
Then, the
quotient stack $[\mathbb{R}\op{Spec}_{G}\, (A)/G]$ is a pointed
schematic homotopy type.
\end{prop}
{\bf Proof:} The augmentation map $x : A
\longrightarrow k$ induces a morphism
$* \longrightarrow \mathbb{R}\op{Spec}\, A$, and therefore gives rise
to a well defined point
\[
* \longrightarrow \mathbb{R}\op{Spec}\, A \longrightarrow
[\mathbb{R}\op{Spec}_{G}\, (A)/G].
\]
Consider the natural morphism in $\op{Ho}(\sff{SPr}_{*}(k))$
\[
[\mathbb{R}\op{Spec}_{G}\, (A)/G] \longrightarrow
\opi{BG}.
\]
Its homotopy fiber is isomorphic to $\mathbb{R}\opi{Spec}\, A$.
Using the general theorems \cite[Corollary~2.4.10, Theorem~3.2.4]{t1}
and the long exact sequence on homotopy sheaves, it is then
enough to
show that $\mathbb{R}\opi{Spec}\, A$ is a connected affine
stack. Equivalently, we need to show that
$A$ is cohomologically connected, which is true by hypothesis.
The proposition is proven. \hfill $\Box$

\bigskip

We now define a category
$\op{Ho}(\sff{EqAlg}^{\Delta,*}_{\mathbb{U}})$, of 
\emph{equivariant augmented and 1-connected cosimplicial
algebras} in the following way. The objects of
$\op{Ho}(\sff{EqAlg}^{\Delta,*}_{\mathbb{U}})$ 
are triplets $(G,A,u)$, consisting of

\begin{itemize}
\item an affine group scheme $G$ in $\mathbb{U}$

\item a $G$-equivariant cosimplicial algebra $A$ in $\mathbb{U}$ such that
$H^{0}(A)\simeq k$ and $H^{1}(A)=0$

\item  a morphism of $G$-equivariant cosimplicial algebras $u : A
  \rightarrow \mathcal{O}(G)$ 
(where $\mathcal{O}(G)$ is the algebra of functions on $G$ considered as
a $G$-equivariant cosimplicial algebra via the regular representation
  of $G$).  
\end{itemize}

A morphism $f : (G,A,u) \longrightarrow (H,B,v)$ in
$\op{Ho}(\sff{EqAlg}^{\Delta,*}_{\mathbb{U}})$ 
is a pair $(\phi,a)$, consisting of

\begin{itemize}
\item a morphism of affine group schemes $\phi : G \longrightarrow H$ 
\item a morphism in $\op{Ho}(\opi{G-}\sff{Alg}^{\Delta}/\mathcal{O}(G))$
  (the homotopy 
  category of 
the model category of objects over $\mathcal{O}(G)$)
$a : B \longrightarrow A$. Here, $B$ is considered as an object in
  $\opi{G-}\sff{Alg}^{\Delta}/\mathcal{O}(G)$ 
by composing the action and the augmentation with
the morphisms $G \rightarrow H$ and $\mathcal{O}(H) \rightarrow
  \mathcal{O}(G)$. 
\end{itemize}

This defines the category
$\op{Ho}(\sff{EqAlg}^{\Delta,*}_{\mathbb{U}})$. Note however that this
category is not (a priori) the homotopy category of a model category
$\sff{EqAlg}^{\Delta,*}_{\mathbb{U}}$ (at least we did not define any
such model category), and is defined in an ad-hoc way.

Let $(G,A,u) \in \op{Ho}(\sff{EqAlg}^{\Delta,*}_{\mathbb{U}})$. We
consider the associated 
$G$-equivariant stack $\mathbb{R}\op{Spec}_{G}\, A$. The morphism
$u : A \longrightarrow \mathcal{O}(G)$ induces a natural morphism
\[
\mathbb{R}\op{Spec}_{G}\, \mathcal{O}(G)\simeq G \longrightarrow
\mathbb{R}\op{Spec}_{G}\, A,
\]
where $G$ acts on itself by left translations. Therefore,
$\mathbb{R}\op{Spec}_{G}\, A$ 
can be considered in a natural way as an object in
$\op{Ho}(G/\opi{G-}\sff{SPr}(k))$, 
the homotopy category 
of $G$-equivariant stacks under $G$. Applying the
quotient stack 
functor of \cite[Definition~1.2.2]{kpt} we get a well defined object
\[
\left([G/G]\simeq * \longrightarrow [\mathbb{R}\op{Spec}_{G}\, A/G]
\right) \in \op{Ho}([G/G]/\sff{SPr}(k))\simeq 
\op{Ho}(\sff{SPr}_{*}(k)).
\]
Tracing carefully through the definitions, it is straightforward to
 check
 that the previous construction defines a functor
\begin{equation} \label{eq:eqalg}
\xymatrix@R-2pc{
\op{Ho}(\sff{EqAlg}^{\Delta,*}_{\mathbb{U}})^{\op{op}} \ar[r] &
\op{Ho}(\sff{SPr}_{*}(k))  \\
 (G,A,u) \ar@{|->}[r] & [\mathbb{R}\op{Spec}_{G}\, A/G].
}
\end{equation}

\begin{prop}\label{peq3}
The functor \eqref{eq:eqalg} is fully faithful and its essential image
consists exactly of 
psht.
\end{prop}
{\bf Proof:} The full faithfulness follows easily from
Proposition~\ref{peq2}. Also, 
Proposition~\ref{peqpsht} implies that this functor takes values in
the sub-category of psht. It remains to prove that any psht is
in the essential image of \eqref{eq:eqalg}.

Let $F$ be a psht. By Theorem~\ref{t1} we can write $F$ as $\bB G_{*}$
where $G_{*}$ is a fibrant object in $\sff{sGAff}_{\mathbb{U}}$. We consider
the projection $F \longrightarrow BG$, where $G:=\pi_{0}(G_{*})\simeq
\pi_{1}(F)$, as well as the cartesian square
\[
\xymatrix{
F \ar[r] & BG \\
F^{0} \ar[r] \ar[u] & EG, \ar[u]}
\]
where as usual $EG$ is the simplicial presheaf having $EG_{m} := G^{m+1}$
and the face and degeneracy maps are given by the projections and
diagonal embeddings (see \cite[\S 1.3]{t1} for details). The
simplicial presheaf $F^{0}$ is a pointed affine scheme in $\mathbb{U}$
equipped with a natural action of $G$. Taking its cosimplicial algebra
of functions one gets a $G$-equivariant cosimplicial algebra
$A:=\mathcal{O}(F^{0})$. As the fiber of $F^{0} \longrightarrow F$ is
naturally isomorphic to $G$ one gets furthermore a $G$-equivariant
morphism $G \longrightarrow F^{0}$, giving rise to an $G$-equivariant
morphism $u : A \longrightarrow \mathcal{O}(G)$. One can then check
that $F^{0}$ is naturally isomorphic in $\op{Ho}([G/G]-\sff{SPr}(k))$
to the image of the object $(G,A,u) \in
\op{Ho}(\opi{G-}\sff{Alg}^{\Delta}/\mathcal{O}(G))^{\op{op}}$ by the
functor $\mathbb{R}\op{Spec}_{G}$ (this follows from the fact that $F^{0}$
is an affine stack). Finally, one has isomorphisms in
$\op{Ho}(\sff{SPr}_{*}(k))$
\[
\left(* \rightarrow F\right) \simeq \left([G/G] \rightarrow [F^{0}/G]\right)
\simeq \left([\mathbb{R}\op{Spec}_{G}\, \mathcal{O}(G)/G] \rightarrow
       [\mathbb{R}\op{Spec}_{G}\, A/G]\right),
\]
showing that $F$ belongs to the essential image of the required
functor. \hfill $\Box$ 

\

\medskip 

\begin{cor}\label{ccomp}
The categories $\op{Ho}(\sff{Hopf}^{\Delta}_{\mathbb{U}})$ and
$\op{Ho}(\sff{EqAlg}_{\mathbb{U}}^{\Delta,*})$ 
are equivalent.
\end{cor}
{\bf Proof:} Indeed by Theorem  \ref{t1} and Proposition \ref{peq3}
they are both equivalent to the 
full sub-category of $\op{Ho}(\sff{SPr}_{*}(k))$ consisting of
psht. \hfill $\Box$  

\begin{rmk}
\emph{The previous corollary is a generalization of the equivalence between
reduced nilpotent Hopf dg-algebras and $1$-connected reduced dga, given
by the bar and cobar constructions (see \cite[\S 0]{ta}). It could be
interesting to produce 
explicit functors between $\op{Ho}(\sff{Hopf}^{\Delta}_{\mathbb{U}})$ and
$\op{Ho}(\sff{EqAlg}_{\mathbb{U}}^{\Delta,*})$ without passing
through the category of psht.} 
\end{rmk}

\begin{center} \textit{An explicit model for $\sch{X}{k}$}
\end{center}  

Let $X$ be a pointed and connected simplicial set in $\mathbb{U}$. In
this paragraph we give an explicit
model for $\sch{X}{k}$ which is based on the notion of
equivariant affine stacks we just introduced.

The main idea of the construction is the following observation.  Let
$G:=\pi_{1}(X,x)^{\op{alg}}$ be the pro-algebraic completion of
$\pi_{1}(X,x)$ over $k$.  By Lemma \ref{leqstack} and Corollary \ref{c2},
the natural morphism $\sch{X}{k} \longrightarrow \opi{BG}$ corresponds
to a $G$-equivariant affine stack. Furthermore, the universal property
of the schematization suggests that the corresponding $G$-equivariant
cosimplicial algebra is the cosimplicial algebra of cochains of $X$
with coefficients in the local system $\mathcal{O}(G)$. We will show
that this guess is actually correct.

Let $\pi_{1}(X,x) \longrightarrow G$
be the universal morphism, and let $X \longrightarrow
B(G(k))$ be the corresponding morphism of simplicial sets.
This latter morphism is well defined up to homotopy, and we choose once
and for all a representative for it.
Let $p : P \longrightarrow X$ be the corresponding $G$-torsor in
$\sff{SPr}(k)$. More precisely,
$P$ is the simplicial presheaf sending an affine scheme $\op{Spec}
A \in \sff{Aff}$ to the simplicial set
$P(A):=(\opi{EG}(A)\times_{\opi{BG}(A)} X)$, where $\opi{EG}(A)
\longrightarrow \opi{BG}(A)$ is the natural projection.
The morphism $p : P \longrightarrow X$ is then a well defined morphism
in $\op{Ho}(\opi{G-}\sff{SPr}(k))$. Here the group
$G$ is acting on $P=(\opi{EG}\times_{\opi{BG}}X)$ by its action on
$\opi{EG}$,  and trivially on $X$.
Alternatively we can describe $P$ by the formula
\[
P\simeq (\widetilde{X}\times G)/\pi_{1}(X,x),
\]
where $\widetilde{X}$ is the universal covering of $X$, and
$\pi_{1}(X,x)$ acts on $\widetilde{X}\times G$ by the diagonal
action (our convention here is that $\pi_{1}(X,x)$ acts on $G$ by left
translation). We assume at this point that $\widetilde{X}$ is chosen to
be  cofibrant in the model
category of $\pi_{1}(X,x)$-equivariant simplicial sets. For example,
we may assume that $\widetilde{X}$ is
a $\pi_{1}(X,x)$-equivariant cell complex.

We consider now the $G$-equivariant affine stack
$\mathbb{R}\op{Spec}_{G}\, \mathcal{O}_{G}(P) \in
\op{Ho}(\opi{G-}\sff{SPr}(k))$,
which comes naturally equipped with its adjunction morphism $P
\longrightarrow \mathbb{R}\op{Spec}_{G}\,
\mathbb{L}\mathcal{O}_{G}(P)$.
This induces a well defined morphism in $\op{Ho}(\sff{SPr}(k))$:
\[
X \simeq [P/G] \longrightarrow [\mathbb{R}\op{Spec}_{G}\,
\mathbb{L}\mathcal{O}_{G}(P)/G].
\]
Furthermore, since $X$ is pointed, this morphism induces a natural morphism in
$\op{Ho}(\sff{SPr}_{*}(k))$
\[
f : X \longrightarrow [\mathbb{R}\op{Spec}_{G}\,
\mathbb{L}\mathcal{O}_{G}(P)/G].
\]
With this notation we now have the following important

\begin{thm}\label{teq1}
The natural morphism $f : X \longrightarrow [\mathbb{R}\op{Spec}_{G}\,
\mathbb{L}\mathcal{O}_{G}(P)/G]$
is a model for the schematization of $X$.
\end{thm}
{\bf Proof:} Since $P\simeq (\widetilde{X}\times
G)/\pi_{1}(X,x)$, the algebra $\mathbb{L}\mathcal{O}_{G}(P)$ can be
identified with the cosimplicial algebra of cochains on $X$ with
coefficients in the local system of algebras $\mathcal{O}(G)$. More
precisely, $P$ is equivalent to the homotopy colimit of
$\widetilde{X}\times G$ viewed as a $\pi_{1}(X,x)$-diagram in
$\opi{G-}\sff{SPr}(k)$. As $\mathcal{O}_{G}$ is left Quillen, one has
equivalences
\[
\mathbb{L}\mathcal{O}_{G}(P)\simeq
\op{holim}_{\pi_{1}(X,x)}\mathbb{L}\mathcal{O}_{G}(\widetilde{X}\times G)
\simeq
\op{holim}_{\pi_{1}(X,x)}\mathbb{L}\mathcal{O}_{G}(G)^{\widetilde{X}},
\]
where $(-)^{\widetilde{X}}$
is the exponential functor (which is part of the simplicial structure
on $\opi{G-}\sff{Alg}^{\Delta}$). In particular we have an isomorphism
\[
\mathbb{L}\mathcal{O}_{G}(P)\simeq
(\mathcal{O}(G)^{\widetilde{X}})^{\pi_{1}(X,x)},
\]
of $\mathbb{L}\mathcal{O}_{G}(P)$ with
the $\pi_{1}(X,x)$-invariant $G$-equivariant cosimplicial algebra of
$\mathcal{O}(G)^{\widetilde{X}}$. Note that the identification
\[
\op{holim}_{\pi_{1}(X,x)}(\mathcal{O}(G)^{\widetilde{X}}) \simeq
\op{lim}_{\pi_{1}(X,x)}(\mathcal{O}(G)^{\widetilde{X}})
=(\mathcal{O}(G)^{\widetilde{X}})^{\pi_{1}(X,x)}
\]
uses the fact that $\widetilde{X}$ is cofibrant as a
$\pi_{1}(X,x)$-simplicial set.  The underlying cosimplicial algebra
of $\mathbb{L}\mathcal{O}_{G}(P)$ is therefore augmented,
cohomologically connected and belongs to $\mathbb{U}$. Therefore, by
Proposition \ref{peqpsht} we conclude that
$[\mathbb{R}\opi{Spec}_{G}\, \mathbb{L}\mathcal{O}_{G}(P)/G]$ is a
pointed schematic homotopy type.

In order to finish the proof of the theorem it remains to show that
the morphism $f$ is a $\boldsymbol{P}$-equivalence. For this we start
with an algebraic group $H$ and a finite dimensional linear
representation $V$ of $H$.  Let $F$ denote the pointed schematic
homotopy type $K(H,V,n)$, and let $F \longrightarrow BH$ be the
natural projection. It is instructive to observe that $F$ is naturally
isomorphic to $[\mathbb{R}\opi{Spec}_{H}\, B(V,n)/H]$, where $B(V,n)$
is the cosimplicial cochain algebra of $K(V,n)$ taken together with
the natural action of $H$.

There is a commutative diagram
$$\xymatrix{
\mathbb{R}\underline{\op{Hom}}_{*}([\mathbb{R}\op{Spec}_{G}\,
\mathbb{L}\mathcal{O}_{G}(P)/G],F) \ar[r]^-{f^{*}} \ar[d]_{p} &
\mathbb{R}\underline{\op{Hom}}_{*}(X,F) \ar[d]^{q} \\
\op{Hom}_{\sff{GAff}}(G,H) \ar[r] &
\op{Hom}_{\sff{Gp}}(\pi_{1}(X,x),H),}$$ in which the horizontal
morphism at the bottom is an isomorphism due to the universal
property of the map $\pi_{1}(X,x) \longrightarrow G$.
Therefore it suffices
 to check that $f^{*}$ induces an equivalence on the homotopy
fibres of the two projections $p$ and $q$.
Let $\rho : G \longrightarrow H$ be a morphism of groups, and
consider the homotopy fiber of $F \longrightarrow BH$,
together with the $G$-action induced from $\rho$. This is a
$G$-equivariant affine stack that will be denoted by $F_{G}$, and
whose underlying stack is isomorphic to $K(V,n)$.

The homotopy fibers of $p$ and $q$ at $\rho$ are isomorphic to
$\mathbb{R}\underline{\op{Hom}}_{G}
(\mathbb{R}\op{Spec}_{G}\,\mathbb{L}\mathcal{O}_{G}(P),F_{G})$  and
$\mathbb{R}\underline{\op{Hom}}_{G}(P,F_{G})$
respectively.
But since
$F_{G}$ is a $G$-equivariant affine stack Prop. \ref{peq2} implies
that the natural morphism
\[
\mathbb{R}\underline{\op{Hom}}_{G}(
\mathbb{R}\op{Spec}_{G}\,\mathbb{L}\mathcal{O}_{G}(P),F_{G})
\longrightarrow \mathbb{R}\underline{\op{Hom}}_{G}(P,F_{G})
\]
is an equivalence. \hfill $\Box$

\

\bigskip

Consider now $\mathcal{O}(G)$ as a locally constant sheaf of algebras
on $X$ via the natural action of $\pi_{1}(X,x)$, and let
\[
C^{\bullet}(X,\mathcal{O}(G)):=
(\mathcal{O}(G)^{\widetilde{X}})^{\pi_{1}(X,x)}
\]
be the cosimplicial algebra of cochains on $X$ with coefficients in
$\mathcal{O}(G)$. This cosimplicial algebra is equipped with a
natural $G$-action, induced by the regular representation of $G$. One
can thus consider $C^{\bullet}(X,\mathcal{O}(G))$ as a object in
$\opi{G-}\sff{Alg}^{\Delta}$.  Thus the previous theorem, and the natural
identification $\mathbb{L}\mathcal{O}_{G}(P)\simeq
C^{\bullet}(X,\mathcal{O}(G))$ as objects in
$\op{Ho}(\opi{G-}\sff{Alg}^{\Delta})$ immediately yield the following:

\begin{cor}\label{cteq1}
With the previous notations, one has
\[
\sch{X}{k}\simeq [\mathbb{R}\op{Spec}_{G}
C^{\bullet}(X,\mathcal{O}(G))/G].
\]
\end{cor}

\

\begin{center} \textit{The case of characteristic zero} \end{center}

Suppose that $k$ is of characteristic zero. In this case
Corollary~\ref{cteq1} can be reformulated in terms of the
pro-reductive completion of the fundamental group. 

Let $x\in X$ be a pointed $\mathbb{U}$-small simplicial set, and
$G^{\op{red}}$ be the pro-reductive completion of $\pi_{1}(X,x)$ over
$k$. By definition, $G^{\op{red}}$ is the universal reductive affine
group scheme over $k$ equipped with a morphism from $\pi_{1}(X,x)$
with Zariski dense image.  We consider $\mathcal{O}(G^{\op{red}})$ as
a local system of $k$-algebras on $X$, and consider the cosimplicial
algebra $C^{*}(X,\mathcal{O}(G^{\op{red}}))$ as an object in
$G^{\op{red}}-\sff{Alg}^{\Delta}$. The
proof of Theorem \ref{teq1} can also be adapted to obtain 
the following result.

\begin{cor}\label{cteq1'}
There exists a natural isomorphism in $\op{Ho}(\sff{SPr}_{*}(k))$
\[
\sch{X}{k} \simeq
[\mathbb{R}\op{Spec}_{G^{\op{red}}}
C^{\bullet}(X,\mathcal{O}(G^{\op{red}}))/G^{\op{red}}].
\]
\end{cor}
{\bf Proof:} The natural morphism
\[
f : X \longrightarrow [\mathbb{R}\op{Spec}_{G^{\op{red}}}
C^{\bullet}(X,\mathcal{O}(G^{\op{red}}))/G^{\op{red}}]=:\mathcal{X}
\]
is defined the same way as in the proof of theorem \ref{teq1}.
Let $H$ be a reductive
linear algebraic group $H$, and $V$ be a linear representation
$V$ of $H$. We set $F:=K(H,V,n)$, and we let 
$\op{Hom}_{\sff{Gp}}^{\op{zd}}(\pi_{1}(X,x),H)$ be the
subset of morphisms $\pi_{1}(X,x) \longrightarrow H$ 
with a Zariski dense image.  We also 
let 
$\mathbb{R}\underline{\op{Hom}}^{\op{zd}}(X,F)$ be
defined by the homotopy pull-back diagram
\[
\xymatrix{
\mathbb{R}\underline{\op{Hom}}^{\op{zd}}(X,F)
\ar[r] \ar[d] & \mathbb{R}\underline{\op{Hom}}(X,F) \ar[d] \\
\op{Hom}_{\sff{Gp}}^{\op{zd}}(\pi_{1}(X,x),H) \ar[r] & 
\op{Hom}_{\sff{Gp}}(\pi_{1}(X,x),H).}
\]
In the same way, we define
$\op{Hom}_{\sff{Gp}}^{\op{zd}}(\pi_{1}(\mathcal{X},x),H)$ to be the
subset of morphisms $\pi_{1}(\mathcal{X},x) \longrightarrow H$ 
with a Zariski dense image. Notice that these are precisely
the morphisms $\pi_{1}(\mathcal{X},x) \longrightarrow H$ 
that factor through the natural quotient
$\pi_{1}(\mathcal{X},x) \longrightarrow G^{\op{red}}$, 
that is
\[
\op{Hom}_{\sff{Gp}}^{\op{zd}}(\pi_{1}(\mathcal{X},x),H)=
\op{Hom}_{\sff{Gp}}(G^{\op{red}},H).
\]
Finally, we define
$\mathbb{R}\underline{\op{Hom}}^{\op{red}}(\mathcal{X},F)$ by the
homotopy pull-back diagram
\[
\xymatrix{
\mathbb{R}\underline{\op{Hom}}^{\op{red}}_{*}(\mathcal{X},F)
\ar[r] \ar[d] & \mathbb{R}\underline{\op{Hom}}_{*}(\mathcal{X},F) \ar[d] \\
\op{Hom}_{\sff{Gp}}(G^{\op{red}},H) \ar[r] & 
\op{Hom}_{\sff{Gp}}(\pi_{1}(\mathcal{X},x),H).}
\]
We have a 
commutative diagram
\[
\xymatrix{
\mathbb{R}\underline{\op{Hom}}_{*}^{\op{red}}(\mathcal{X},F) 
\ar[r]^-{f^{*}} \ar[d]_{p} &
\mathbb{R}\underline{\op{Hom}}_{*}^{\op{zd}}(X,F) \ar[d]^{q} \\
\op{Hom}_{\sff{GAff}}(G^{\op{red}},H) \ar[r] &
\op{Hom}_{\sff{Gp}}^{\op{zd}}(\pi_{1}(X,x),H).}
\]
As in the proof
of theorem \ref{teq1} the induced morphisms on the homotopy 
fibers of the vertical morphisms $p$ and $q$ are all weak equivalences.
Since the bottom horizontal morphism is bijective, this 
shows that the top horizontal morphism
\[
f^{*} : \mathbb{R}\underline{\op{Hom}}_{*}^{\op{red}}(\mathcal{X},F)
\longrightarrow 
\mathbb{R}\underline{\op{Hom}}_{*}^{\op{zd}}(X,F)
\]
is a weak equivalence. 

As this is true for any $H$, $V$ and $n$, we find that
the natural morphism 
\[
\pi_{1}(X,x) \longrightarrow
\pi_{1}(\mathcal{X},x)
\]
induces an equivalence of the categories
of reductive linear representations, and that for any such 
reductive representation $V$, the induced morphism
\[
f^{*} : H^{*}(\mathcal{X},V)\longrightarrow
H^{*}(X,V)
\]
is an isomorphism. This in turn implies that 
the morphism $f$ is a $\boldsymbol{P}$-equivalence. \ \hfill $\Box$ 

\section{Some properties of the schematization functor}
\label{sec:properties} 

In this section we have collected some basic properties of the
schematization functor. All these are purely topological but are
geared toward the schematization of algebraic varieties.

For a topological space $X$, we will write $\sch{X}{k}$ instead
of $(S(X)\otimes k)^{sch}$, where $S(X)$ is the singular simplicial
set of $X$. 

\subsection{Schematization of smooth manifolds}

The purpose of this section is to give a formula for the
schematization of a smooth manifold in terms of differential forms.

\

\smallskip

Let $X$ be a topological space, which is assumed to be locally
contractible and for which any open subset is paracompact.  For
example, $X$ could be any $CW$-complex, including that way any
covering of a compact and smooth manifold. Let $G$ be an affine group
scheme, and let us consider $\sff{Rep}(G)(X)$, the abelian category of
sheaves on $X$ with values in the category $\sff{Rep}(G)$ of linear
representations of $G$. The category $\sff{Rep}(G)(X)$ is in
particular the category of sheaves on a (small) site with values in a
Grothendieck abelian category. Thus, $\sff{Rep}(G)(X)$ has enough
injectives.  This implies that the opposite category possesses enough
projectives, and therefore \cite[II.4.11,Remark~5]{q2} can be applied
to endow the category $C^{+}(\sff{Rep}(G)(X))$ of positively graded
complexes in $\sff{Rep}(G)(X)$ with a model structure. Recall that the
equivalences in $C^{+}(\sff{Rep}(G)(X))$ are morphisms inducing
isomorphisms on cohomology sheaves (i.e. quasi-isomorphisms of
complexes of sheaves), and the fibrations are epimorphisms whose
kernel $K$ is such that each $K_{n}$ is an injective object in
$\sff{Rep}(G)(X)$.  One checks immediately that
$C^{+}(\sff{Rep}(G)(X))$ is a cofibrantly generated model category.

The Dold-Kan correspondence yields an
equivalence of categories
\[
D : C^{+}(\sff{Rep}(G)(X)) \longrightarrow \sff{Rep}(G)(X)^{\Delta} \qquad
C^{+}(\sff{Rep}(G)(X)) \longleftarrow \sff{Rep}(G)(X)^{\Delta} : N,
\]
where $N$ is the normalization functor and $D$ the denormalization
functor. Through this equivalence we can transplant the model
structure of $C^{+}(\sff{Rep}(G)(X))$ to a model structure on
$\sff{Rep}(G)(X)^{\Delta}$.  As explained in \cite[II,Example
2.8]{gj}, the category $\sff{Rep}(G)(X)^{\Delta}$ possesses a natural
simplicial structure, and it is easy to check that this simplicial
structure is compatible with the model structure. Therefore,
$\sff{Rep}(G)(X)^{\Delta}$ is a simplicial and cofibrantly generated
model category.

The two categories $\sff{Rep}(G)(X)^{\Delta}$ and
$C^{+}(\sff{Rep}(G)(X))$ have natural symmetric monoidal structures
induced by the usual tensor product on $\sff{Rep}(G)$. The
monoidal structures turn these categories into symmetric monoidal
model categories. However, the functors $D$ and $N$ are not monoidal
functors, but are related to the monoidal structure via the usual
Alexander-Whitney and shuffle products
\[
\op{aw}_{X,Y} : N(X)\otimes N(Y) \longrightarrow N(X \otimes Y) \qquad
\op{sp}_{X,Y} :
D(X)\otimes D(Y) \longrightarrow D(X\otimes Y).
\]
The morphisms $\op{aw}_{X,Y}$ are unital and associative, and
$\op{sp}_{X,Y}$ are unital, associative and commutative.

We will denote the categories of commutative monoids in the categories
$\sff{Rep}(G)(X)^{\Delta}$ and
$C^{+}(\sff{Rep}(G)(X))$ respectively by
\[
\opi{G-}\sff{Alg}^{\Delta}(X), \text{ and } \opi{G-}\sff{CDGA}(X),
\]
and call them the categories of $G$-equivariant cosimplicial algebras
on $X$, and of $G$-equivariant commutative differential algebras on
$X$.

The following proposition is standard

\begin{prop} \label{pa1}
\begin{enumerate}
\item There exist a unique simplicial model structure on
$\opi{G-}\sff{Alg}^{\Delta}$ such that a morphism is an equivalence
(respectively a fibration) if and only if the induced morphism in
$\sff{Rep}(G)(X)^{\Delta}$ is an equivalence, i.e. induces a
quasi-isomorphism on the normalized complexes of sheaves (respectively
a fibration, i.e. induces an epimorphism on the normalized complexes
of sheaves whose kernel is levelwise injective).
\item There exist a unique simplicial model structure on
$\opi{\opi{G-}\sff{CDGA}}(X)$ such that
a morphism is an equivalence (respectively a fibration) if and only if
the induced morphism
in $C^{+}(\sff{Rep}(G)(X))$ is an equivalence, i.e. a
quasi-isomorphism of complexes of sheaves
(respectively a fibration, i.e. a epimorphism of complexes of sheaves
whose kernel is 
levelwise injective).
\end{enumerate}
\end{prop}

\

\smallskip

\noindent
By naturality, the functor $\op{Th}$ of Thom-Sullivan cochains
(see \cite[4.1]{his}) extends to a functor
\[
\op{Th} : \op{Ho}(\opi{G-}\sff{CDGA}(X)) \longrightarrow
\op{Ho}(\opi{G-}\sff{Alg}^{\Delta}(X)).
\]
This functor is an equivalence, and its inverse is the denormalization
functor
\[
D : \op{Ho}(\opi{G-}\sff{Alg}^{\Delta}(X)) \longrightarrow
\op{Ho}(\opi{G-}\sff{CDGA}(X)).
\]
Let $f : X \longrightarrow Y$ be a continuous map of topological spaces
(again paracompact
and loaclly contractible). The inverse and direct image functors of
sheaves induce  Quillen adjunctions
\[
\xymatrix@1{
\opi{G-}\sff{CDGA}(Y) \ar@<1ex>[r]^-{f^{*}} & \opi{G-}\sff{CDGA}(X)
\ar@<1ex>[l]^-{f_{*}}
}
\quad \text{ and } \quad
\xymatrix@1{
\opi{G-}\sff{Alg}^{\Delta}(Y) \ar@<1ex>[r]^-{f^{*}} &
\opi{G-}\sff{Alg}^{\Delta}(X). \ar@<1ex>[l]^-{f_{*}}}
\]
Furthermore, one checks immediately that the following diagram 
\[
\xymatrix{\op{Ho}(\opi{G-}\sff{CDGA}(X)) \ar[r]^-{D} \ar[d]_-{\mathbb{R}f_{*}}
& \op{Ho}(\opi{G-}\sff{Alg}^{\Delta}(X)) \ar[d]^-{\mathbb{R}f_{*}} \\
\op{Ho}(\opi{G-}\sff{CDGA}(Y)) \ar[r]_-{D} &
\op{Ho}(\opi{G-}\sff{Alg}^{\Delta}(Y))}
\]
commutes. 
As usual, in the case when $Y=*$, the functor $f_{*}$ will be
denoted simply  by $\Gamma(X,-)$. The previous
diagram should be understood  as a functorial isomorphism
$D\mathbb{R}\Gamma(X,A)\simeq \mathbb{R}\Gamma(X,D(A))$, for any $A
\in \op{Ho}(\opi{G-}\sff{CDGA}(X))$.

The category $\sff{Rep}(G)(X)^{\Delta}$ is naturally enriched over the category
of sheaves of simplicial sets on $X$. Indeed, for $F \in \sff{SSh}(X)$
a simplicial sheaf, and $E \in \sff{Rep}(G)(X)^{\Delta}$, one can
define $F\otimes E$ to be the sheaf associated to the presheaf defined
by the following formula
\[
\xymatrix@R=1pt@C=9pt{
\sff{Open}(X)^{\op{op}} \ar[r] & \sff{Rep}(G)^{\Delta} \\
U \ar@{|->}[r] & F(U)\otimes E(U),
}
\]
where $\sff{Open}(X)$ denotes the category of open sets in $X$ and
$F(U)\otimes E(U)$ is viewed as an object in $\sff{Rep}(G)^{\Delta}$
via the natural simplicial structure on the model category
$\sff{Rep}(G)^{\Delta}$. It is straightforward to check that this
structure makes the model category $\sff{Rep}(G)(X)^{\Delta}$ into a
$\sff{SSh}(X)$-model category in the sense of \cite[\S 4.2]{ho},
where $\sff{SSh}(X)$ is taken with the injective model structure
defined in \cite{ja}.  In particular, one can define functors
\[
\xymatrix@R=1pt@C=9pt{
\sff{SSh}(X)^{\op{op}}\otimes \sff{Rep}(G)(X)^{\Delta} \ar[r] &
\sff{Rep}(G)(X)^{\Delta} \\
(F,E) \ar@{|->}[r] & E^{F},
}
\]
and
\[
\xymatrix@=6pt{
(\sff{Rep}(G)(X)^{\Delta})^{\op{op}}\otimes \sff{Rep}(G)(X)^{\Delta}
\ar[r] & \sff{SSh}(X) \\
(E,E') \ar@{|->}[r] & \underline{\op{Hom}}_{\sff{SSh}(X)}(E,E'),
}
\]
which are related by the usual adjunction isomorphisms
$$
\underline{\op{Hom}}_{\sff{SSh}(X)}(F\otimes E,E')\simeq
\underline{\op{Hom}}_{\sff{SSh}(X)}(E,(E')^{F})\simeq
\underline{\op{Hom}}_{\sff{SSh}(X)}(F,
\underline{\op{Hom}}_{\sff{SSh}(X)}(E,E')),$$
where $\underline{\op{Hom}}_{\sff{SSh}(X)}$ denotes the internal
$\op{Hom}$ in $\sff{SSh}(X)$.

\begin{lem}\label{la1}
Let $A \in \sff{Rep}(G)(X)$ be an injective object. Then, the presheaf
of abelian groups
on $X$ underlying $A$ is acyclic.
\end{lem}
{\bf Proof:} As the space $X$ is paracompact its sheaf cohomology
coincides with its \v{C}ech cohomology. Thus it is enough to show that for
any open covering $\{U_{i}\}_{i \in I}$, the \v{C}ech complex
\[
A(X) \longrightarrow \prod_{i\in I}A(U_{i}) \longrightarrow \ldots
\longrightarrow
\prod_{(i_{0},\ldots,i_{p})\in I^{p+1}}A(U_{i_{0}}\cap \ldots \cap
U_{i_{p}}) \longrightarrow \ldots
\]
is exact.

Since $A$ is an injective object, it is fibrant as a constant
cosimplicial object $A \in \sff{Rep}(G)(X)^{\Delta}$.
Let $\{U_{i}\}_{i\in I}$ be an open covering of $X$, and let
$N(U/X)$ be its nerve.
This is the simplicial sheaf on $X$ whose sheaf of $n$-simplices is
defined by the following formula
\[
\xymatrix@=6pt{
N(U/X)_{n} \; : & \sff{Open}(X)^{\op{op}} \ar[r] & \sff{SSet} \\
& V \ar@{|->}[r] & \prod_{(i_{0},\ldots,i_{n}) \in
I^{n+1}}\op{Hom}_{X}(V,U_{i_{0}}\cap \ldots \cap U_{i_{n}}).
}
\]
The simplicial sheaf $N(U/X)$ is contractible (i.e. equivalent to $*$
in $\sff{SSh}(X)$) as can be easily seen on stalks. Furthermore, since
every object in $\sff{SSh}(X)$ is cofibrant, and $A \in
\sff{Rep}(G)(X)^{\Delta}$ is fibrant, the natural morphism
\[
A^{\bullet}\simeq A \longrightarrow A^{N(U/X)}
\]
is an equivalence of fibrant objects in $\sff{Rep}(G)(X)^{\Delta}$
(here $A^{N(U/X)}$ is part of the $\sff{SSh}(X)$-module structure
on $\sff{Rep}(G)(X)^{\Delta})$). Therefore, since  the global
sections functor
\[
\Gamma(X,-) : \sff{Rep}(G)(X)^{\Delta} \longrightarrow
\sff{Rep}(G)^{\Delta}
\]
is right Quillen, we conclude
that the induced morphism $\Gamma(X,A)\longrightarrow
\Gamma(X,A^{N(U/X)})$ is a quasi-isomorphism.
Since $\Gamma(X,A^{N(U/X)})$ is just the \v{C}ech complex
of $A$ for the covering $\{U_{i}\}_{i \in I}$, this implies that $A$
is acyclic on $X$. \hfill $\Box$

\

\begin{lem}\label{la2}
Let $A \in \opi{G-}\sff{CDGA}(X)$ be a $G$-equivariant commutative
differential graded algebra over $X$ such that
for all $n\geq 0$, the sheaf of abelian groups $A_{n}$ is an acyclic
sheaf on $X$. Then, the natural morphism
in $\op{Ho}(\opi{G-}\sff{CDGA})$
\[
\Gamma(X,A) \longrightarrow \mathbb{R}\Gamma(X,A)
\]
is an isomorphism.
\end{lem}
{\bf Proof:} Let $A \longrightarrow RA$ be a fibrant model of $A$
in $\opi{G-}\sff{CDGA}(X)$. Then, by the definition
of a fibration in $\opi{G-}\sff{CDAG}(X)$, each $A_{n}$ is an injective
object in $\sff{Rep}(G)(X)$, and is therefore
acyclic by lemma \ref{la1}. The morphism $A \longrightarrow RA$ is
thus a quasi-isomorphism of
complexes of acyclic sheaves of abelian groups on $X$, which implies
that the induced morphism
on global sections
\[
\Gamma(X,A) \longrightarrow \Gamma(X,RA)
\]
is a quasi-isomorphism of complexes. By the definition of an
equivalence in
$\opi{G-}\sff{CDGA}$ this implies that
the morphism
\[
\Gamma(X,A) \longrightarrow \Gamma(X,RA)=:\mathbb{R}\Gamma(X,A)
\]
is actually an isomorphism in $\op{Ho}(\opi{G-}\sff{CDGA})$. \hfill $\Box$

\

\medskip

Now, let $(X,x)$ be
a pointed connected compact smooth manifold and  let $X^{\op{top}}$
denote the underlying
topological space of $X$. Let $L_{B}(X)$ be the
category of semi-simple local systems of
finite dimensional $\mathbb{C}$-vector spaces on $X^{\op{top}}$. It is
a rigid $\mathbb{C}$-linear tensor category which is naturally equivalent
to the category of finite dimensional semi-simple
representations of the fundamental group $\pi_{1}(X^{\op{top}},x)$.

The category of semi-simple $C^{\infty}$ complex vector bundles with flat
connections on $X$ will be denoted as before by $L_{DR}(X)$. Recall
that the category $L_{DR}(X)$ is a rigid $\mathbb{C}$-linear tensor
category, and the functor which maps a flat bundle to its monodromy
representations at $x$, induces an equivalence
of tensor categories $L_{B}(X)\simeq L_{DR}(X)$ (this is again the
Riemann-Hilbert correspondence).

Let $G_{X}:=\pi_{1}(X^{\op{top}},x)^{\op{red}}$ be the pro-reductive completion
of the group $\pi_{1}(X^{\op{top}},x)$. Note that
it is the Tannaka dual of the category $L_{B}(X)$.
The algebra $\mathcal{O}(G_{X})$ can be viewed as the
left regular representation of $G_{X}$. Through the
universal morphism $\pi_{1}(X^{\op{top}},x) \longrightarrow G_{X}$, we can
also consider ${\mathcal O}(G_{X})$ as
a linear representation of $\pi_{1}(X^{\op{top}},x)$. This linear
representation is not finite dimensional, but it is
admissible in the sense that it equals the union of its finite
dimensional sub-representations. Therefore,
the algebra $\mathcal{O}(G_{X})$ corresponds to an object in the
$\mathbb{C}$-linear tensor category $T_{B}(X)$, of
$\op{Ind}$-local systems on $X^{\op{top}}$. By convention all of
our $\op{Ind}$-objects are labelled by
$\mathbb{U}$-small index categories.

Furthermore, the algebra structure on $\mathcal{O}(G_{X})$, gives rise
to a morphism
\[
\mu : \mathcal{O}(G_{X}) \otimes \mathcal{O}(G_{X}) \longrightarrow
\mathcal{O}(G_{X}),
\]
which is easily checked to be a morphism in
$T_{B}(X)$. This means that if
$\mathcal{O}(G_{X})$ is written as
the colimit of local systems $\{V_{i}\}_{i \in
I}$, then the product $\mu$ is given
by a compatible system of morphisms in $L_{B}(X)$
\[
\mu_{i,k} : V_{i} \otimes V_{i} \longrightarrow V_{k},
\]
for some index $k \in I$ with $i\leq k \in I$.  The morphism $\mu = \{
\mu_{i,k} \}_{i,k \in I}$, endows the object $\mathcal{O}(G_{X}) \in
T_{B}(X)$ with a structure of a commutative unital monoid.  Through
the Riemann-Hilbert correspondence $T_{B}(X)\simeq T_{DR}(X)$, the
algebra $\mathcal{O}(G_{X})$ can also be considered as a commutative
monoid in the tensor category $T_{DR}(X)$ of $\op{Ind}$-objects in
$L_{DR}(X)$.

Let $\{(V_{i},\nabla_{i})\}_{i \in I} \in T_{DR}(X)$ be the
object corresponding to $\mathcal{O}(G_{X})$.
For any $i \in I$, one can form the de Rham complex of
$C^{\infty}$-differential forms
$$(A^{\bullet}_{DR}(V_{i}),\nabla_{i}):=\xymatrix@1{A^{0}(V_{i})
\ar[r]^-{\nabla_{i}} & A^{1}(V_{i})
\ar[r]^-{\nabla_{i}} & \ldots  \ar[r]^-{\nabla_{i}}&
A^{n}(V_{i}) \ar[r]^-{D_{i}} &  \ldots}.$$
In this way we obtain an inductive system of complexes
$\{(A^{\bullet}_{DR}(V_{i}),D_{i})\}_{i \in I}$ whose colimit
complex is defined to be the de Rham complex of the local system 
$\mathcal{O}(G_{X})$ on $X$:
\[
(A^{\bullet}_{DR}(\mathcal{O}(G_{X})),\nabla) :=
\op{colim}_{i \in I}(A^{\bullet}_{DR}(V_{i}),\nabla_{i}).
\]
The complex $(A^{\bullet}_{DR}(\mathcal{O}(G_{X})),\nabla)$ has a
natural structure of a commutative differential graded algebra, coming
from the commutative monoid structure on $\{(V_{i},\nabla_{i})\}_{i
\in I} \in T_{DR}(X)$. Using wedge products of differential forms, and
the monoidal structure, we obtain in the usual fashion morphisms of
complexes
$$(A^{\bullet}_{DR}(V_{i}),\nabla_{i}) \otimes
(A^{\bullet}_{DR}(V_{j}),\nabla_{j}) \longrightarrow
(A^{\bullet}_{DR}(V_{k}),\nabla_{k})$$
which, after passing to the colimit along $I$, turn
$(A^{\bullet}_{DR}(\mathcal{O}(G_{X})),\nabla)$ into a commutative
differential graded algebra.

The affine group scheme $G_{X}$ acts via the right
regular representation on the $\op{Ind}$-local
system $\mathcal{O}(G_{X})$, and this
action is compatible with the algebra structure.
By functoriality,
this action gives rise to
an action of $G_{X}$ on the corresponding objects in
$T_{DR}(X)$. Furthermore, if
$G_{X}$ acts on an inductive system of flat bundles $(V_{i},\nabla_{i})$,
then it acts
naturally on its de Rham complex
$\op{colim}_{i \in I}(A^{\bullet}_{DR}(V_{i}),\nabla_{i})$, by
acting on the
spaces of differential forms with coefficients in the various
$V_{i}$. Indeed, if
the action of $G_{X}$ is given by a co-module structure
\[
\{(V_{i},\nabla_{i})\}_{i \in I}\longrightarrow \{\mathcal{O}(G_{X})
\otimes (V_{i},\nabla_{i})\}_{i \in I},
\]
then one obtains a morphism of $\op{Ind-}{C}^{\infty}$-bundles by tensoring
with the sheaf $A^{n}$ of differential
forms on $X$
\[
\{V_{i}\otimes A^{n}\}_{i \in I}\longrightarrow
\{\mathcal{O}(G_{X})\otimes (V_{i}\otimes A^{n})\}_{i \in I}.
\]
Taking global sections on $X$, one has a morphism
\[
\op{colim}_{i \in I}A^{n}(V_{i}) \longrightarrow \op{colim}_{i \in I}
A^{n}(V_{i})\otimes \mathcal{O}(G_{X}),
\]
which defines an action of $G_{X}$ on the space of differential forms
with values in the $\op{Ind-}{C}^{\infty}$-bundle $\{V_{i}\}_{i \in
I}$. Since this action is compatible with the differentials
$\nabla_{i}$, one obtains an action of $G_{X}$ on the de Rham complex
$\op{colim}_{i \in I
}(A^{\bullet}_{DR}(V_{i}),\nabla_{i})$. Furthermore since the action
is compatible with the algebra structure on $\mathcal{O}(G_{X})$ it
follows that $G_{X}$ acts on $\op{colim}_{i \in
I}(A^{\bullet}_{DR}(V_{i}),\nabla_{i})$ by algebra automorphisms.
Thus, the group scheme $G_{X}$ acts in a natural way on the complex
$(A^{\bullet}_{DR}(\mathcal{O}(G_{X})),\nabla)$, turning it into a
well defined $G_{X}$-equivariant commutative differential graded
algebra.

Applying the denormalization functor $D : \op{Ho}(\opi{G_{X}-}\sff{CDGA})
\longrightarrow 
\op{Ho}(\opi{G_{X}-}\sff{Alg}^{\Delta})$,
we obtain a well defined $G_{X}$-equivariant cosimplicial algebra denoted by
\[
C^{\bullet}_{DR}(X,\mathcal{O}(G_{X})):=
D(A^{\bullet}_{DR}(\mathcal{O}(G_{X})),\nabla)
\in \op{Ho}(\opi{G_{X}-}\sff{Alg}^{\Delta}).
\]
\

\noindent
To summarize: for any $(X,x)$ - a pointed connected smooth manifold we let
$G_{X}:=\pi_{1}(X,x)^{\op{red}}$ to be the pro-reductive
completion of its fundamental group. The $G_{X}$-equivariant
commutative differential graded algebra
of de Rham cohomology of $X$ with coefficients in $\mathcal{O}(G_{X})$
is denoted by
\[
(A^{\bullet}_{DR}(\mathcal{O}(G_{X})),\nabla) \in
\op{Ho}(\opi{G_{X}-}\sff{CDGA}).
\]
Its denormalization is denoted by
\[
C^{\bullet}_{DR}(X,\mathcal{O}(G_{X})):=
D(A^{\bullet}_{DR}(\mathcal{O}(G_{X})),\nabla) \in
\op{Ho}(\opi{G_{X}-}{Alg}^{\Delta}).
\]
\

\smallskip

Any smooth map $f : (Y,y) \longrightarrow
(X,x)$ of pointed connected
smooth manifolds induces a morphism
$G_{Y}:=\pi_{1}(Y,y)^{\op{red}} \longrightarrow
G_{X}:=\pi_{1}(X,x)^{\op{red}}$, and 
therefore a well defined functor
\[
f^{*} : \op{Ho}(\opi{G_{X}-}\sff{Alg}^{\Delta}) \longrightarrow
\op{Ho}(\opi{G_{Y}-}\sff{Alg}^{\Delta}).
\]
It is not difficult to check that the pull-back of differential forms
via $f$ induces
a well defined morphism in $\op{Ho}(\opi{G_{Y}-}\sff{Alg}^{\Delta})$
\[
f^{*} : f^{*}C^{\bullet}_{DR}(X,\mathcal{O}(G_{X}))
\longrightarrow C^{\bullet}_{DR}(Y,\mathcal{O}(G_{Y})).
\]
Furthermore, this morphism depends functorially (in an obvious
fashion) on the
morphism $f$.

The last discussion, allows us to construct a functor
\[
(X,x) \mapsto [\mathbb{R}\op{Spec}_{G_{X}}\,
C^{\bullet}_{DR}(X,\mathcal{O}(G_{X}))/G_{X}],
\]
from the category of pointed connected smooth manifolds to the
category of pointed schematic
homotopy types. To simplify notation, we will denote this functor by
$(X,x) \mapsto (X\otimes \mathbb{C})^{\op{diff}}$. 

The following corollary is a generalization of Poincar\'e lemma.

\begin{cor}\label{ca1}
Let $X$ be a pointed connected compact and smooth manifold. Then,
there exist an isomorphism
in $\op{Ho}(\opi{G-}\sff{CDGA})$
\[
(A_{DR}(\mathcal{O}(G_{X})),D) \simeq
\mathbb{R}\Gamma(X,\mathcal{O}(G_{X})),
\]
which is functorial in $X$.
\end{cor}
{\bf Proof:} This follows from lemma \ref{la2}, and the fact that the
sheaves $A^{n}(\mathcal{O}(G_{X}))$ of differential forms with
coefficients in the $\op{Ind}$-flat bundle associated to
$\mathcal{O}(G_{X})$, are filtered colimits of soft sheaves on $X$. In
particular they are acyclic since $X$ is compact.

The functoriality is straighforward.  \hfill $\Box$

\

\medskip

As a consequence of the Corollary \ref{ca1} we obtain the following:

\begin{cor}\label{ca2}
Let $X$ be a pointed connected compact and smooth manifold. Then
there exist an isomorphism
in $\op{Ho}(\opi{G-}\sff{Alg}^{\Delta})$
\[
D(A_{DR}(\mathcal{O}(G_{X})),D) \simeq
\mathbb{R}\Gamma(X,\mathcal{O}(G_{X})),
\]
where $D$ is the denormalization functor, and $\mathcal{O}(G_{X})$ is
considered
as an object in $\opi{G-}\sff{Alg}^{\Delta}(X)$. This isomorphism is
furthermore 
functorial in $X$.
\end{cor}

Let $(S,s)$ be a pointed connected simplicial set, let $\pi:=\pi_{1}(S,s)$
be its fundamental group, and let
$G_{S}$ be the pro-reductive completion of $\pi$. The conjugation
action of $\pi$ on $G_{S}$ is an action
by group scheme automorphism,
and therefore gives rise to a natural action on the model category
$\opi{G_{S}-}\sff{Alg}^{\Delta}$. More precisely, if $\gamma \in \pi$
and if $A \in \opi{G_{S}-}\sff{Alg}^{\Delta}$ corresponds to a cosimplicial
$\mathcal{O}(G_{S})$-co-module $E$, one defines
$\gamma\cdot E$ to be the co-module structure
\[
\xymatrix{E \ar[r] & E\otimes \mathcal{O}(G_{S})
\ar[r]^-{\opi{Id}\otimes
\gamma} & E\otimes \mathcal{O}(G_{S})}.
\]
We will be concerned with the fixed-point model category of
$\opi{G_{S}-}\sff{Alg}^{\Delta}$ under the group $\pi$,
denoted by $(\opi{G_{S}-}\sff{Alg}^{\Delta})^{\pi}$, and described in
\cite{kpt}.

The category $(\opi{G_{S}-}\sff{Alg}^{\Delta})^{\pi}$ is naturally
enriched over the category $\opi{\pi-}\sff{SSet}$, of
$\pi$-equivariant simplicial sets. Indeed, for $X \in
\opi{\pi-}\sff{SSet}$ and $A \in
(\opi{G_{S}-}\sff{Alg}^{\Delta})^{\pi}$ one can define $X\otimes A$
whose underlying $G_{S}$-equivariant cosimplicial algebra is
$X^{for}\otimes A^{\op{for}}$ (where $X^{\op{for}}$ and $A^{\op{for}}$
are the underlying simplicial set and $G_{S}$-equivariant
cosimplicial algebra of $X$ and $A$ respectively, and
$X^{\op{for}}\otimes A^{\op{for}}$ uses the simplicial structure of
$\opi{G_{S}-}\sff{Alg}^{\Delta}$). The action of $\pi$ on
$X^{\op{for}}\otimes A^{\op{for}}$ is then defined diagonally. In
particular, we 
can use the exponential product
\[
A^{X}\in (\opi{G_{S}-}\sff{Alg}^{\Delta})^{\pi}, \text{ where } X \in
\opi{\pi-}\sff{SSet} \text{ and } A \in
(\opi{G_{S}-}\sff{Alg}^{\Delta})^{\pi}.
\]
The functor
\[
\xymatrix@R=1pt@C=10pt{
(\opi{\pi-}\sff{SSet})^{\op{op}}\times (\opi{G_{S}-}\sff{Alg}^{\Delta})^{\pi} 
\ar[r] &
(\opi{G_{S}-}\sff{Alg}^{\Delta})^{\pi} \\
(X,A) \ar@{|->}[r] & A^{X}
}
\]
is a bi-Quillen functor (see \cite[$\S 4$]{ho}), and can be derived
into a functor
\[
\xymatrix@R=1pt@C=10pt{
\op{Ho}(\opi{\pi-}\sff{SSet})^{\op{op}}\times
\op{Ho}((\opi{G_{S}-}\sff{Alg}^{\Delta})^{\pi})  
\ar[r] & \op{Ho}((\opi{G_{S}-}\sff{Alg}^{\Delta})^{\pi}) \\
(X,A) \ar@{|->}[r] & A^{\mathbb{R}X}.
}
\]
Recall that by definition
\[
A^{\mathbb{R}X}:=(RA)^{QX},
\]
where $RA$ is a fibrant model for $A$ in
$(\opi{G_{S}-}\sff{Alg}^{\Delta})^{\pi}$, and $QX$ is a cofibrant model
for $X$ in $\opi{\pi-}\sff{SSet}$.

In the following defintion, $\mathcal{O}(G_{S})$ is considered
together with its $\pi$ and $G_{S}$ actions, i.e. is viewed
as an object in $(\opi{G_{S}-}\sff{Alg}^{\Delta})^{\pi}$.

\begin{df}\label{da2}
The $G_{S}$-equivariant cosimplicial algebra of cochains of $S$ with
coefficients in the local
system $\mathcal{O}(G_{S})$ is defined to be
\[
C^{\bullet}(S,\mathcal{O}(G_{S})):=
\mathcal{O}(G_{S})^{\mathbb{R}\widetilde{S}}\in
\op{Ho}(\opi{G_{S}-}\sff{Alg}^{\Delta}),
\]
where $\widetilde{S} \in \op{Ho}(\opi{\pi-}\sff{SSet})$ is the
universal covering of 
$S$.
\end{df}

Since $\mathcal{O}(G_{S})$ is always a fibrant object in
$(\opi{G_{S}-}\sff{Alg}^{\Delta})^{\pi}$, one has
\[
C^{\bullet}(S,\mathcal{O}(G_{S}))\simeq
\mathcal{O}(G_{S})^{\widetilde{S}}\in
\op{Ho}(\opi{G_{S}-}{Alg}^{\Delta}),
\]
as soon as $\widetilde{S}$ is chosen to be cofibrant in
$\opi{\pi-}\sff{SSet}$ (e.g. is chosen to be a $\pi$-equivariant cell
complex).

Going  back to our space $X$, which is assumed to be
pointed and connected, let
$S:=S(X)$ be its singular simplicial set, naturally pointed by the
image $s \in S$ of $x \in X$. Observe that in this
case one has a  natural isomorphism $G_{S}\simeq G_{X}$, induced by the
natural isomorphism
$\pi_{1}(X,x)\simeq \pi_{1}(S,s)$.

\begin{prop}\label{la3}
There exist an isomorphism in $\op{Ho}(\opi{G_{S}-}\sff{Alg}^{\Delta})$
\[
C^{\bullet}(S,\mathcal{O}(G_{S})) \simeq
\mathbb{R}\Gamma(X,\mathcal{O}(G_{S})),
\]
where $\mathcal{O}(G_{S})$ is considered as an object in
$\opi{G_{S}-}\sff{Alg}^{\Delta}(X)$. Furthermore,
this isomorphism is functorial in $X$.
\end{prop}
{\bf Proof:} We will prove the existence of the
isomorphism. The functoriality statement is straightforward and is
left to the reader.

Let $\pi = \pi_{1}(X,x)\simeq \pi_{1}(S,s)$
Let $\widetilde{X} \longrightarrow X$ be the universal cover of
$X$. Let $\widetilde{S} \longrightarrow
S(\widetilde{X})$ be a cofibrant replacement of $S(\widetilde{X}) \in
\opi{\pi-}\sff{SSet}$ and let  $p : \widetilde{S} \longrightarrow
S(X)$ be the natural
$\pi$-equivariant projection. For any open subset $U \subset X$, we
will denote by $\widetilde{S}_{U}$ the fiber product
\[
\widetilde{S}_{U}:=\widetilde{S}\times_{S(X)}S(U) \in
\opi{\pi-}\sff{SSet}.
\]
Consider the presheaf $C^{\bullet}(-,\mathcal{O}(G_{S}))$ of
$G_{S}$-equivariant cosimplicial algebras on $X$, defined by
\[
\xymatrix@=6pt{
C^{\bullet}(-,\mathcal{O}(G_{S})) \; : &
\sff{Open}(X)^{\op{op}} \ar[r] & \opi{G_{S}-}\sff{Alg}^{\Delta} \\
& U \ar@{|->}[r] & C^{\bullet}(U,\mathcal{O}(G_{S}))
:=\mathcal{O}(G_{S})^{\widetilde{S}_{U}},
}
\]
where $\mathcal{O}(G_{S})^{\widetilde{S}_{U}}$ is the exponentiation
of $\mathcal{O}(G_{S}) \in (\opi{G_{S}-}\sff{Alg}^{\Delta})^{\pi}$
by $\widetilde{S}_{U} \in \opi{\pi-SSet}$. We denote by
$aC^{\bullet}(-,\mathcal{O}(G_{S})) \in \opi{G_{S}-}\sff{Alg}^{\Delta}(X)$
the associated sheaf.

We have a natural morphism in $\opi{G_{S}-}\sff{Alg}^{\Delta}$
\[
C^{\bullet}(X,\mathcal{O}(G_{S})) \longrightarrow
\Gamma(X,aC^{\bullet}(-,\mathcal{O}(G_{S})))
\longrightarrow \Gamma(X,RaC^{\bullet}(-,\mathcal{O}(G_{S})))
\longleftarrow
\Gamma(X,R\mathcal{O}(G)_{S}),
\]
where $RaC^{\bullet}(-,\mathcal{O}(G_{S})))$ is a fibrant replacement
of $aC^{\bullet}(-,\mathcal{O}(G_{S})))$, and
$R\mathcal{O}(G)_{S}$ is a fibrant replacement of $\mathcal{O}(G)_{S}$.
Since  $X$  is assumed to be locally contractible, the natural
morphism in $\opi{G_{S}-}\sff{Alg}^{\Delta}(X)$
\[
\mathcal{O}(G_{S}) \longrightarrow C^{\bullet}(-,\mathcal{O}(G_{S})),
\]
induced over every open subset $U \subset X$ by the projection
$\widetilde{S}_{U} \longrightarrow *$,
is an equivalence. Therefore, it only remains to show that the morphism
\[
C^{\bullet}(X,\mathcal{O}(G_{S})) \longrightarrow
\Gamma(X,RaC^{\bullet}(-,\mathcal{O}(G_{S}))
\]
is an equivalence in $\opi{G_{S}-Alg}^{\Delta}$.

For this, let $U_{\bullet} \longrightarrow X$ be an open hyper-cover of $X$
such that each $U_{n}$ is the disjoint union
of contractible open subset of $X$ (such a hyper-cover exists again
due to local contractibility of  $X$). The
simplicial sheaf represented by $U_{\bullet}$ is equivalent to $*$
in $\sff{SSh}(X)$. As
$RaC^{\bullet}(-,\mathcal{O}(G_{S}))$ is fibrant we obtain a natural
equivalence in $\opi{G_{S}-}\sff{Alg}^{\Delta}$
\[
\Gamma(X,RaC^{\bullet}(-,\mathcal{O}(G_{S}))\simeq
\Gamma(X,RaC^{\bullet}(-,\mathcal{O}(G_{S}))^{U_{\bullet}})\simeq
\op{holim}_{[n]\in
  \Delta}\Gamma(U_{n},RaC^{\bullet}(-,\mathcal{O}(G_{S}))).
\]
Furthermore, it is shown is \cite[Lemma~2.10]{t2} that the natural
morphism
\[
\op{hocolim}_{[n]\in \Delta^{\op{op}}}U_{\bullet} \longrightarrow X
\]
is a weak equivalence, and therefore
\[
\op{hocolim}_{[n]\in \Delta^{\op{op}}}\widetilde{S}_{U_{\bullet}}
\longrightarrow \widetilde{S}_{X}=\widetilde{S}
\]
is an equivalence in $\opi{\pi-}\sff{SSet}$. This implies that the natural
morphism in $\opi{G_{S}-}\sff{Alg}^{\Delta}$
\[
C^{\bullet}(X,\mathcal{O}(G_{S})) \longrightarrow \op{holim}_{[n]
\in \Delta}
C^{\bullet}(U_{n},\mathcal{O}(G_{S}))
\]
is an equivalence. Moreover, there exist a commutative diagram
\[
\xymatrix{
C^{\bullet}(X,\mathcal{O}(G_{S})) \ar[r] \ar[d] & \opi{holim}_{[n] \in
\Delta}
C^{\bullet}(U_{n},\mathcal{O}(G_{S})) \ar[d] \\
\Gamma(X,RaC^{\bullet}(-,\mathcal{O}(G_{S})) \ar[r] &
\op{holim}_{[n]\in
\Delta}\Gamma(U_{n},RaC^{\bullet}(-,\mathcal{O}(G_{S}))),}
\]
i.e. without a loss of generality we may assume
that $X$ is contractible. Under this assumption it only remains to
check that the natural morphism
\[
\mathcal{O}(G_{S}) \longrightarrow
\mathbb{R}\Gamma(X,\mathcal{O}(G_{S}))
\]
is an isomorphism in $\op{Ho}(\opi{G_{S}-}\sff{Alg}^{\Delta})$. But,
by lemma 
\ref{la1} and the description of fibrant objects
in $\opi{G_{S}-}\sff{Alg}^{\Delta}(X)$, one has
\[
H^{n}(\mathbb{R}\Gamma(X,\mathcal{O}(G_{S})))\simeq
H^{n}(X,\mathcal{O}(G_{S})).
\]
Since the space $X$ is paracompact and contractible, its sheaf
cohomology with coefficients in the
constant local system $\mathcal{O}(G_{S})$ is trivial, which completes
the proof of the proposition. \ \hfill $\Box$

\

\medskip

As an immediate corollary of Corollary\ref{ca2} and Proposition
\ref{la3} one
obtains the following description of the schematization of
a smooth manifold in terms of its differential forms.

\begin{prop}\label{pdiff}
Let $X$ be a pointed connected compact and smooth manifold. Then,
there exist an 
isomorphism 
in $\op{Ho}(\opi{G-}\sff{Alg}^{\Delta})$
\[
D(A_{DR}(\mathcal{O}(G_{X})),D) \simeq
C^{\bullet}(S(X),\mathcal{O}(G_{X})),
\]
where $D$ is the denormalization functor, and $\mathcal{O}(G_{X})$ is
considered
as an object in $\opi{G-}\sff{Alg}^{\Delta}(X)$. In other words, one has a
natural isomorphism 
\[
\sch{X}{k}\simeq
[\mathbb{R}\op{Spec}_{G_{X}}D(A_{DR}(\mathcal{O}(G_{X})),D)/G_{X}].
\] 
\end{prop}
\

\subsection{A schematic Van Kampen theorem}

Let $X$ be a pointed and connected ($\mathbb{U}$-small) topological
space, and $\{U_{i}\}_{i \in I}$ be a finite open covering. We will
assume that each $U_{i}$ contains the base point, and that each of the
$p+1$-uple intersection $U_{i_{0},\dots,i_{p+1}}:=U_{i_{0}}\cap \dots
\cap U_{i_{p}}$ is connected. Here we have in mind the example of a smooth
projective complex variety $X$ covered by Zariski open subsets
containing the base point. We form the poset $N(U)$ (the nerve of
$\{U_{i}\}_{i\in I}$), whose objects are strings of indices
$(i_{0},\dots,i_{p}) \in I^{p+1}$ (for various lenghts $p\geq 0$), and
where $(i_{0},\dots,i_{p})\leq (j_{0},\dots,j_{p})$ if and only if
$U_{i_{0},\dots,i_{p+1}} \subset U_{j_{0},\dots,j_{p+1}}$. There
exists a natural functor
\[
\xymatrix@R-2pc{
N(U)  \ar[r] & \sff{Top}^{\op{con}}_{*} \\
(i_{0},\dots,i_{p}) \ar@{|->}[r] & U_{i_{0},\dots,i_{p}},
}
\]
where $\sff{Top}_{*}^{\op{con}}$ is the category of pointed and
connected $\mathbb{U}$-topological spaces.  This diagram is
furthermore augmented to the constant diagram $X$. Applying the
schematization functor, we get a diagram of pointed simplicial
presheaves
\[
\xymatrix@R-2pc{
N(U) \ar[r] & \sff{SPr}_{*}(k) \\
(i_{0},\dots,i_{p}) \ar@{|->}[r] & \sch{U_{i_{0},\dots,i_{p}}}{k}
}
\]
which is naturally augmented towards $\sch{X}{k}$ (here we use a
version of the schematization functor which is defined on the level of
simplicial sets and not only on its homotopy category. For example,
one can use the functor $Z \mapsto \bB\op{Spec}\, \mathcal{O}(GZ^{\op{alg}})^{P}$
where $(-)^{P}$ is a cofibrant replacement functor for the P-local model structure on 
$\sff{Hopf}^{\Delta}$, see 
 Corollary~\ref{c3'}).

\begin{prop}\label{p3}
For any psht $F$, the natural morphism
\[
\mathbb{R}\underline{\op{Hom}}(\sch{X}{k},F) \longrightarrow
\op{Holim}_{\alpha \in
  N(U)}\mathbb{R}\underline{\op{Hom}}(\sch{U_{\alpha}}{k},F),
\]
 induced by the augmentation, is an equivalence of simplicial sets.
\end{prop}
{\bf Proof:} Using the universal property of the schematization
functor it is enough to prove that 
the natural morphism
\[
\op{Hocolim}_{\alpha \in N(U)}U_{\alpha} \longrightarrow X
\]
is a weak equivalence of topological spaces. When $X$ is a CW complex this is
something well known (see for example \cite[Lemma~2.10]{t2}). In the
general case it is enough to 
consider the commutative square
\[
\xymatrix{
\op{Hocolim}_{\alpha \in N(U)}|S(U_{\alpha})| \ar[r] \ar[d] & |S(X)| 
\ar[d] \\
\op{Hocolim}_{\alpha \in N(U)}U_{\alpha} \ar[r] &  X,}
\]
where $|S(-)|$ is the geometric realization of the singular
functor. Since the functors $|-|$ and $S(-)$ form 
a Quillen equivalence, the vertical morphisms are both weak
equivalences. The top 
horizontal morphism being an equivalence this finishes the proof of
the proposition. \hfill $\Box$  

\begin{rmk}
Another interpretation of Proposition \ref{p3} is to say that
$\sch{X}{k}$ is the homotopy colimit, {\em in the category of psht},
of the diagram $\alpha \mapsto U_{\alpha}$. This point of view however
is somewhat trickier since, when apropriately defined, the homotopy
colimits of psht are not the same as homotopy colimits of pointed
simplicial presheaves. These subtleties go beyond the scope of the
present paper and we will not discuss them here.
\end{rmk}

Proposition \ref{p3} is a Van Kampen type result, as it states that
the schematization of a space $X$ is uniquely determined by the
diagram $\alpha \mapsto \sch{U_{\alpha}}{k}$.  This property
will be useful only if the space $X$ behaves well
locally with respect to the schematization functor. This is for
example the case when $X$ is a smooth projective complex
variety. Indeed, locally for the Zariski topology $X$ is a $K(\pi,1)$,
where $\pi$ is a successive extension of free groups of finite
type. In  Theorem~\ref{thm-good-extension} below, we will show that 
such groups $\pi$
are \textit{algebraically good}, and therefore $\sch{K(\pi,1)}{k}$
is itself $1$-truncated. So in a way the schematization of a
smooth projective complex variety is well understood locally for the
Zariski topology, and proposition \ref{p3} tells us that globally the
schematization is obtained by the homotopy colimits of the
schematization of a covering.

\subsection{Good groups}

Recall that a discrete group $\Gamma$ in $\mathbb{U}$ is algebraically
good (relative to the field $k$) if the natural morphism
\[
\sch{K(\Gamma,1)}{k} \longrightarrow K(\Gamma^{\op{alg}},1)
\]
is an isomorphism. The terminology \textit{algebraically good} mimicks
the corresponding pro-finite notion introduced by J.P.Serre in
\cite{se}. It is justified by
the following lemma. Let $H^{\bullet}_{H}(\Gamma^{\op{alg}},V)$ denote the
Hochschild cohomology of
the affine group scheme $\Gamma^{\op{alg}}$ with coefficients in a
linear representation $V$ (as defined in \cite[I]{sga3}).

\begin{lem}\label{lb1}
Let $\Gamma$ be a group in $\mathbb{U}$ and $\Gamma^{\op{alg}}$ be
its pro-algebraic completion. Then,
$\Gamma$ is an algebraically good group if and only if for every
finite dimensional  linear
representation $V$ of $\Gamma^{\op{alg}}$, the natural morphism
\[
H^{\bullet}_{H}(\Gamma^{\op{alg}},V) \longrightarrow
H^{\bullet}(\Gamma,V)
\]
is an isomorphism.
\end{lem}
{\bf Proof:} Using \cite[$\S 1.5$]{t1} and \cite[Corollary
$3.3.3$]{t1}, the fact that $H^{\bullet}_{H}(\Gamma^{\op{alg}},V)
\cong H^{\bullet}(\Gamma,V)$ for all $V$ is just a reformulation of
the fact that $K(\Gamma,1) \longrightarrow K(\Gamma^{\op{alg}},1)$
is a $\boldsymbol{P}$-equivalence. Since $K(\Gamma^{\op{alg}},1)$ is a
pointed 
schematic homotopy type, this 
implies the lemma. \hfill $\Box$

\

\medskip

\noindent
We will also use the following very general fact.

\begin{prop}\label{prop-gen}
Let $\Gamma$ be a group in $\mathbb{U}$, $\Gamma^{\op{alg}}$ be
its pro-algebraic completion, and $n>1$ an integer. The following are
equivalent.
\begin{enumerate}
\item[(1)] For any linear representation $V$ of $\Gamma^{\op{alg}}$
the induced morphism
\[
H^{i}_{H}(\Gamma^{\op{alg}},V) \longrightarrow H^{i}(\Gamma,V)
\]
is an isomorphism for $i<n$ and injective for $i=n$.
\item[(2)] For any linear representation $V$ of $\Gamma^{\op{alg}}$ the
induced morphism
\[
H^{i}_{H}(\Gamma^{\op{alg}},V) \longrightarrow H^{i}(\Gamma,V)
\]
is surjective for $i<n$.
\end{enumerate}
\end{prop}
{\bf Proof:} Let $\sff{Rep}(\Gamma)$ be the category of linear complex
representations of $\Gamma$ (possibly infinite
dimensional), and $\sff{Rep}(\Gamma^{\op{alg}})$ the category of linear
representation of the affine group scheme $\Gamma^{\op{alg}}$ (again
maybe of infinite dimension). The proposition follows from \cite[Lemma
$11$]{ad} applied to the inclusion functor
$\sff{Rep}(\Gamma^{\op{alg}}) \longrightarrow \sff{Rep}(\Gamma)$.  \ \hfill
$\Box$

\

\bigskip

As an immediate corollary we obtain that finite groups are
algebraically good over $k$. Furthermore for any finitely
generated group $\Gamma$ the universal property of the pro-algebraic
completion implies that the natural map $\Gamma \to \Gamma^{\op{alg}}$
induces an isomorphism on cohomology in degrees $0$ and $1$. Therefore
the previous
proposition implies that the natural map
$H^{2}_{H}(\Gamma^{\op{alg}},V)\longrightarrow H^{2}(\Gamma,V)$ is always
injective. In
particular any free group of finite type will be algebraically good
over $k$. In the same vein we have:

\begin{prop}\label{prop-good}
The following groups are algebraically good groups.
\begin{itemize}
\item[(1)] Abelian groups of finite type.
\item[(2)] Fundamental groups of a compact Riemann surface
when $k=\mathbb{C}$.
\end{itemize}
\end{prop}
{\bf Proof:} $(1)$ For any abelian group of finite type $M$, there
exist a morphism of exact sequences
\[
\xymatrix{
0\ar[r] \ar[d] & H \ar[r] \ar[d] & M \ar[r] \ar[d] & M_{0} \ar[r]
\ar[d] & 0 \ar[d] \\
0\ar[r]  & H \ar[r]  & M^{\op{alg}} \ar[r] & M_{0}^{\op{alg}} \ar[r]
& 0,}
\]
where $H$ is a finite group, and $M_{0}$ is torsion free. The
comparison of the associated Leray spectral
sequences shows that without a loss of generality one may assume $M$
to be torsion free. Then, again comparing the two Leray spectral
sequences for the rows of the diagram
\[
\xymatrix{
0\ar[r] \ar[d] & \mathbb{Z}^{n-1} \ar[r] \ar[d] & \mathbb{Z}^{n}
\ar[r] \ar[d] & \mathbb{Z} \ar[r] \ar[d] & 0 \ar[d] \\
0\ar[r]  & (\mathbb{Z}^{n-1})^{\op{alg}} \ar[r]  & (\mathbb{Z}^{n})^{\op{alg}}
\ar[r] & \mathbb{Z}^{\op{alg}} \ar[r]  & 0,}
\]
and proceeding by induction on the rank of $M$, one reduces the proof
to the case $M=\mathbb{Z}$. But,
as $\mathbb{Z}$ is a free group it is algebraically good.

\medskip

\noindent
$(2)$ Note that by $(1)$ one can suppose that $g>2$.
As $\Gamma_{g}$ is the fundamental group
of a compact Riemann surface of genus $g$, one has
$H^{i}(\Gamma_{g},V)=0$ for any $i>2$ and
any linear representation $V$. Applying \ref{prop-gen}
it is enough to prove that for any finite dimensional linear
representation $V$ of $\Gamma$, the 
natural morphism
$H^{2}_{H}(\Gamma_{g}^{\op{alg}},V) \longrightarrow
H^{2}(\Gamma_{g},V)$ is 
an isomorphism. 

\begin{lem}
Let $V$ be any linear representation of finite dimension of
$\Gamma_{g}$, and
$x \in H^{2}(\Gamma_{g},V)$. Then, there exist a linear representation
of finite dimension $V'$, and
an injective morphism \linebreak 
$j : V \hookrightarrow V'$, such that $j(x)=0\in
H^{2}(\Gamma_{g},V')$.
\end{lem}
{\bf Proof:} Using Poincare duality the assertion of the lemma can
be restated as follows. Given any finite dimensional linear
representation $E$ of $\Gamma_{g}$ and any invariant vector $e \in
H^{0}(\Gamma_{g}, E) = E^{\Gamma_{g}}$ show
that there exists a finite dimensional linear representation $E^{'}$
and a $\Gamma_{g}$-equivariant surjection $q : E'\to
E$ so that $e \not\in \op{im}[H^{0}(\Gamma_{g}, E') \to
H^{0}(\Gamma_{g}, E)]$. Let $\mathbb{I}$ denote the one dimensional
trivial representation of $\Gamma_{g}$. Then the element $e \in
H^{0}(\Gamma_{g}, E)$ can be viewed as an injective homomorphism $e :
\mathbb{I} \to E$ of finite dimensional $\Gamma_{g}$-modules. Let $F :=
E/{\mathbb I}$ be the corresponding quotient module and let $\Sigma$
be an irreducible $r$-dimensional representation of $\Gamma_{g}$ for
some $r \geq 1$. The short exact sequence of $\Gamma_{g}$-modules
\[
0 \to {\mathbb I} \stackrel{e}{\to}  E \to F \to 0
\]
induces a long exact sequence of $\op{Ext}$'s in the category of
$\Gamma_{g}$-modules:
\[
\xymatrix@1{
\ldots \ar[r] & \op{Ext}^{1}(F,\Sigma) \ar[r] &  \op{Ext}^{1}(E,\Sigma)
\ar[r]^-{\opi{ev}_{e}} & \op{Ext}^{1}({\mathbb I},\Sigma) \ar[r] &
\op{Ext}^{2}(F,\Sigma) \ar[r] & \ldots,
}
\]
Now observe that $\dim \op{Ext}^{1}({\mathbb I},\Sigma)
= \dim H^{1}(\Gamma_{g},\Sigma) \geq
\chi(\Gamma_{g},\Sigma) = 2r(g-1)$. Furthermore by Poincare duality we
have $\op{Ext}^{2}(F,\Sigma) \cong \op{Hom}(\Sigma,F)^{\vee}$ and
since $\Sigma$ contains only finitely many irreducible
subrepresentations of $\Gamma_{g}$ it follows that for the generic
choice of $\Sigma$, the natural map
\[
\xymatrix@1{\op{Ext}^{1}(E,\Sigma)
\ar[r]^-{\opi{ev}_{e}} & \op{Ext}^{1}({\mathbb I},\Sigma)}
\]
will be surjective. Choose such a $\Sigma$ and let $\epsilon \in
\op{Ext}^{1}(E,\Sigma)$ be an element such that $\opi{ev}_{e}(\epsilon)
\neq 0 \in \op{Ext}^{1}({\mathbb I},\Sigma)$. Let now
\[
0 \to \Sigma \to E' \stackrel{q}{\to} E \to 0
\]
be the extension of $\Gamma_{g}$ modules corresponding to
$\epsilon$. Then the image of $e \in H^{0}(\Gamma_{g},E)$ under the
first edge homomorphism $H^{0}(\Gamma_{g},E) \to
H^{1}(\Gamma_{g},\Sigma) \cong \op{Ext}^{1}({\mathbb I},\Sigma)$ is
precisely $\opi{ev}_{e}(\epsilon) \neq 0$. The lemma is proven. \
\hfill $\Box$

\

\medskip

\noindent
Now, let $V$ by any linear representation of finite dimension of
$\Gamma_{g}$, $x \in H^{2}(\Gamma_{g},V)$,
and let $V \hookrightarrow V'$ as in the previous lemma. Then, the
morphism of long exact sequences
$$\xymatrix{
H^{1}_{H}(\Gamma_{g}^{\op{alg}},V'/V) \ar[r] \ar[d] &
H^{2}_{H}(\Gamma_{g}^{\op{alg}},V) \ar[r] \ar[d] &
H^{2}_{H}(\Gamma_{g}^{\op{alg}},V') \ar[r] \ar[d] & \dots \\
H^{1}(\Gamma_{g},V'/V) \ar[r] & H^{2}(\Gamma_{g},V) \ar[r] &
H^{2}(\Gamma_{g},V') \ar[r] & \dots }$$
shows that $x$ can be lifted to $H^{2}(\Gamma_{g}^{\op{alg}},V)$. This implies that
$H^{2}_{H}(\Gamma_{g}^{\op{alg}},V) \longrightarrow H^{2}(\Gamma_{g},V)$ is
surjective, and finishes the proof
of the proposition. \hfill $\Box$

\

\bigskip

Before we construct more examples of algebraically good groups we
describe a convenient criterion (which was already used implicitly in
the proof of proposition \ref{prop-good}) for checking if a given group
is good.  For a finitely generated group $\Gamma$ we denote by
$\sff{Rep}(\Gamma)$ the category of all representations of $\Gamma$ in
$k$-vector spaces (possibly infinite dimensional) and by
$\sff{Rep}(\Gamma^{\op{alg}})$ the category of all algebraic
representations of the affine pro-algebraic group $\Gamma^{\op{alg}}$.
Recall also that a group $\Gamma$ is of type $(\boldsymbol{F})$ (over
$k$) if: (a) for every $n$ and every finite dimensional complex
representation $V$ of $\Gamma$ the group $H^{n}(\Gamma,V)$ is finite
dimensional, and (b) $H^{\bullet}(\Gamma,-)$ commutes with inductive
limits of finite dimensional complex representations of $\Gamma$.
With this notation we now have:

\begin{lem}\label{lem-good-criterion} Let $\Gamma$ be a finitely
generated group of type $(\boldsymbol{F})$. The following properties
of $\Gamma$ are equivalent:
\begin{itemize}
\item[{\bf (a)}] $\Gamma$ is algebraically good over $k$.
\item[{\bf (b)}] For every positive integer $n$, every finite
dimensional representation $L \in \sff{Rep}(\Gamma)$ and every class
$\alpha \in H^{n}(\Gamma,L)$, there exists an injection
$\imath : L \hookrightarrow W_{\alpha}$
in some finite dimensional $\Gamma$-module $W_{\alpha} \in
\sff{Rep}(\Gamma)$ so that the induced map
$\imath^{*} : H^{n}(\Gamma,L) \to H^{n}(\Gamma,W_{\alpha})$,
annihilates $\alpha$, i.e. $\imath^{*}(\alpha) = 0$.
\item[{\bf (c)}] For every positive integer $n$, every finite
dimensional representation $L \in \sff{Rep}(\Gamma)$,  there exists an
injection $\imath : L \hookrightarrow W_{\alpha}$
in some finite dimensional $\Gamma$-module $W_{\alpha} \in
\sff{Rep}(\Gamma)$ so that the induced map $\imath^{*} :
H^{n}(\Gamma,L) \to H^{n}(\Gamma,W_{\alpha})$,
is identically zero, i.e. $\imath \equiv 0$.
\end{itemize}
\end{lem}
{\bf Proof:} First we show that {\bf (a)} $\Rightarrow$ {\bf
(b)}. Suppose that $\Gamma$ is good over $k$. Fix an integer
$n > 0$, a representation $L \in \sff{Rep}(\Gamma)$ with $\dim(L) <
+\infty$ and some class $\alpha \in H^{n}(\Gamma,L)$.

Consider the
regular representation of $\Gamma^{\op{alg}}$ on the algebra
$\mathcal{O}(\Gamma^{\op{alg}})$ of $k$-valued regular functions
on the affine group $\Gamma^{\op{alg}}$. The natural map $s : \Gamma \to
\Gamma^{\op{alg}}$ allows us to view $\mathcal{O}(\Gamma^{\op{alg}})$ as a
$\Gamma$-module and so the tensor product
\[
L' := L\otimes_{{\mathbb C}} \mathcal{O}(\Gamma^{\op{alg}})
\]
can be interpreted both as a $\Gamma$-module and as a
$\Gamma^{\op{alg}}$-module.

Note that when viewed as an object in $\sff{Rep}(\Gamma^{\op{alg}})$, the
module $L'$ is injective. In particular
\[
H^{i}_{H}(\Gamma^{\op{alg}},L') = 0, \quad \text{ for all } i > 0.
\]
Now write $L'$ as an inductive limit ${\displaystyle  L'=\op{colim}
L'_{i}}$, with finite dimensional $\Gamma$-modules $L'_{i}$.
We consider the natural inclcusion $L \hookrightarrow L'$, which induces
inclusions $L' \hookrightarrow L_{i}'$ for all $i \gg 0$.
Since
$\Gamma$ is algebraically good we have
\[
H^{n}(\Gamma,L) \cong H^{n}_{H}(\Gamma^{\op{alg}},L),
\]
and so we may view $\alpha$ as an element in
$H^{n}_{H}(\Gamma^{\op{alg}},L)$. Furthermore since $\Gamma$ is of type
$(\boldsymbol{F})$ we get
\[
H_{H}^{n}(\Gamma^{\op{alg}},L') =
H_{H}^{n}(\Gamma^{\op{alg}},\op{colim}_{\to}L_{i}') 
= \op{colim}_{\to} H_{H}^{n}(\Gamma^{\op{alg}},L_{i}'),
\]
and since $L'$ was chosen so that $H^{n}_{H}(\Gamma^{\op{alg}},L') = 0$, it
follows that
\[
\alpha \to 0 \in H^{n}_{H}(\Gamma^{\op{alg}},L_{i}')
\]
for all sufficiently big $i$. Combined with the fact that $L
\hookrightarrow 
L_{i}'$ for all $i \gg 0$ this yields the implication {\bf (a)}
$\Rightarrow$ {\bf (b)}.

\medskip

\

\noindent
We will prove the implication {\bf (b)} $\Rightarrow$ {\bf (a)} by
induction on $n$. More precisely, for every integer $n >0$  consider
the condition
\[
(*_{n}) \;  \begin{minipage}[t]{6in}
for every finite dimensional represenation $V \in \sff{Rep}(\Gamma)$
the natural map $s : \Gamma \to \Gamma^{\op{alg}}$ induces an
isomorphism 
\[
s^{*} : H_{H}^{k}(\Gamma^{\op{alg}},V) \stackrel{\cong}{\to}
H^{k}(\Gamma,V) 
\]
for all $k < n$.
\end{minipage}
\]
We need to show that $(*_{n})$ holds for all integers $n >0$.
By the universal property of pro-algebraic completions we know that
$(*_{1})$ holds. This provides the base of the induction.
Assume next that $(*_{n})$ holds. This automatically implies that
\[
s^{*} : H^{n}_{H}(\Gamma^{\op{alg}},L) \to H^{n}(\Gamma,L)
\]
is injective, and so we only have to show that $s^{*} :
H^{n}_{H}(\Gamma^{\op{alg}},L) \to H^{n}(\Gamma,L)$ is also surjective.

Fix $\alpha \in H^{n}(\Gamma,L)$ and let $W \in \sff{Rep}(\Gamma)$ be
a finite dimensional representation for which we can find an injection
$\imath : L \hookrightarrow W$ so that $\imath(\alpha) = 0$. In
particular we can find a class $\beta \in H^{n-1}(\Gamma,W/L)$ which
is mapped to $\alpha$ by the edge homomorphism of the
long exact sequence in cohomology associated to the sequence of
$\Gamma$-modules
\[
0 \to L \stackrel{\imath}{\to} W \to W/L \to 0.
\]
However, by the inductive hypothesis $(*_{n})$ we have an isomorphism
\[
H^{n-1}(\Gamma,W/L) \cong H^{n-1}_{H}(\Gamma^{\op{alg}},W/L),
\]
and so from the commutative diagram
\[
\xymatrix{H_{H}^{n}(\Gamma^{\op{alg}},L) \ar[r]^-{s^{*}} & H^{n}(\Gamma,L) \\
H_{H}^{n-1}(\Gamma^{\op{alg}},W/L) \ar[r]^-{\cong}_-{s^{*}} \ar[u] &
H^{n-1}(\Gamma,W/L) \ar[u]
}
\]
it follows that $\alpha$ comes from $H_{H}^{n}(\Gamma^{\op{alg}},L)$.

\

\medskip

\noindent
The implication {\bf (c)} $\Rightarrow$ {\bf (b)} is obvious. For the
implication {\bf (b)} $\Rightarrow$ {\bf (c)} we need to construct a
finite dimensional representation $W \in \sff{Rep}(L)$ and a
monomorphism $\imath : L \to W$ so that $\imath(H^{n}(\Gamma,L)) = \{ 0
\} \subset H^{n}(\Gamma,W)$. Choose a basis $e_{1}, e_{2}, \ldots,
e_{m}$ of $H^{n}(\Gamma,L)$. By {\bf (b)} we can find monomorphisms
$\imath_{1} : L \hookrightarrow W_{1}$,
$\imath_{2} : L \hookrightarrow W_{2}$, \ldots, $\imath_{m} : L
\hookrightarrow W_{m}$, so that $\imath_{i}(e_{i}) = 0 \in
H^{n}(\Gamma,W_{i})$ for $i = 1, \ldots, m$. Consider now the sequence
of finite dimensional representations $V_{k} \in \sff{Rep}(\Gamma)$
and monomorphisms $\jmath_{k} : L \hookrightarrow V_{k}$ constructed
inductively as follows:
\begin{itemize}
\item $V_{1} := W_{1}$, $\jmath_{1} := \imath_{1}$.
\item Assuming that $\jmath_{k-1} : L \hookrightarrow V_{k-1}$ has
already been constructed, define $V_{k}$ as the pushout:
\[
\xymatrix@M=8pt{
L \; \ar@{^{(}->}[r]^-{-\imath_{k}} \ar@{_{(}->}[d]_-{\jmath_{k-1}} &
W_{k} \ar[d] \\
V_{k-1} \ar[r] & V_{k}
}
\]
in $\sff{Rep}(\Gamma)$. Explicitly $V_{k} = (V_{k-1}\oplus
W_{k})/\operatorname{im}[\xymatrix@1@C=50pt{L \ar[r]^-{\jmath_{k-1}\times
(-\imath_{k})} & V_{k-1}\oplus W_{k}}]$. Moreover the natural map
$\jmath_{k-1}\times \imath_{k} : L \rightarrow V_{k-1}\oplus
W_{k}$ induces a monomorphism $\jmath_{k} : L \hookrightarrow V_{k}$
which completes the step of the induction.
\end{itemize}
Let now $W := V_{m}$ and $\imath := \jmath_{m}$. By construction we
have inclusions $W_{k} \subset W$ for all $k = 1, \ldots, m$ and for
each $k$ the map $\jmath$ factors as
\[
\xymatrix@M=8pt{
L \ar@{^{(}->}[rr]^-{\imath} \ar@{^{(}->}[dr]_-{\imath_{k}} & & W \\
& W_{k}  \ar@{^{(}->}[ur] &
}
\]
In particular $\imath(L) \subset \cap_{k = 1}^{m} W_{k} \subset W$ and
so $\imath(H^{n}(\Gamma,L)) = 0$. The lemma is proven.
\hfill $\Box$

\

\bigskip

\noindent
Using the basic good groups (e.g. free, finite, abelian, surface
groups) as building blocks and the criterion from
Lemma~\ref{lem-good-criterion} we can construct more good groups as
follows:

\begin{thm} \label{thm-good-extension} Suppose that
\begin{equation} \label{good-seq}
1 \to F \to \Gamma \to \Pi \to 1
\end{equation}
is a short exact sequence of finitely generated groups of type
$(\boldsymbol{F})$,  such that
$\Pi$ is algebraically good over $k$ and $F$ is
free. 
$\Gamma$ is algebraically good over $k$.
\end{thm}
{\bf Proof:} The proof is essentially contained in an argument of
Beilinson which appears in \cite[Lemma~$2.2.1$ and $2.2.2$]{bei} in a
slightly different guise. Since the statement of the theorem is an
important ingredient in the localization technique for computing
schematic homotopy types,  we decided to write up the proof
in detail in our context.

Let $p : X \to S$ be the Serre fibration corresponding to the short
exact sequence \eqref{good-seq} and let $Y$ denote the homotopy fiber
of $p$. Note that $X = K(\Gamma,1)$, $S = K(\Pi,1)$ and $Y = K(F,1)$.
Given a representation $L \in \sff{Rep}(\Gamma)$ (respectively in
$\sff{Rep}(\Pi)$ or $\sff{Rep}(F)$) we write $\mathbb{L}$ for the
corresponding local system of $k$-vector spaces on $X$ (respectively
$S$ or $Y$).

To prove that $\Gamma$ is good it suffices (see
Lemma~\ref{lem-good-criterion}) to show that for every positive
integer $n$ and any finite dimensional representation $L \in
\sff{Rep}(\Gamma)$ we can find an injection $\imath : L \to W$ into a
finite dimensional $W \in \sff{Rep}(\Gamma)$ so that the induced map
$\imath : H^{n}(X,\mathbb{L}) \to H^{n}(X,{\mathbb W})$ is identically
zero.

The representation $W$ will be constructed in three steps:

\

\medskip
\noindent
{\bf Step 1.} \quad Fix $n$ and $L$ as above. Then we can find a
finite dimensional $W \in \sff{Rep}(\Gamma)$ and an injection $\imath
: L \hookrightarrow W$ such that the induced map
\[
H^{n}(S,p_{*}\mathbb{L}) \to H^{n}(S,p_{*}\mathbb{W})
\]
is identically zero.

Indeed, since $\Pi$ is assumed to be algebraically good and since
$p_{*}\mathbb{L}$ is a finite dimensional local system on $S$ we can
find (see Lemma~\ref{lem-good-criterion}) a finite dimensional local
system $\mathbb{V}$ on $S$ and an injection $p_{*}\mathbb{L}
\hookrightarrow \mathbb{V}$ which induces the zero map
\[
H^{n}(S,p_{*}\mathbb{L}) \stackrel{0}{\longrightarrow} H^{n}(S,\mathbb{V})
\]
on cohomology.

Now define a local system $\mathbb{W} \to X$ as the pushout
\[
\xymatrix@M=8pt{p^{*}p_{*}\mathbb{L} \ar@{^{(}->}[r] \ar[d] &
p^{*}\mathbb{V} \ar[d] \\
\mathbb{L} \ar@{^{(}->}[r] & \mathbb{W}
}
\]
i.e. $\mathbb{W} := (p^{*}\mathbb{V}\oplus
\mathbb{L})/p^{*}p_{*}\mathbb{L}$. By pushing this cocartesian square
down to $S$ and using the projection formula $p_{*}p^{*}p_{*}\mathbb{L}
\cong p_{*}\mathbb{L}$, we get a commutative diagram of local systems
on $S$:
\[
\xymatrix{
& p_{*}p^{*}\mathbb{V} = \mathbb{V} \ar[dd] \\
p_{*}\mathbb{L} \ar[ur] \ar[dr] & \\
& p_{*}\mathbb{W}.
}
\]
In particular the natural map
\[
H^{n}(S,p_{*}{\mathbb L}) \to H^{n}(S,p_{*}{\mathbb W})
\]
factors through
\[
H^{n}(S,p_{*}{\mathbb L}) \stackrel{0}{\longrightarrow}
H^{n}(S,{\mathbb V})
\]
and so is identically zero. This finishes the proof of Step
1. In fact, the same reasoning can be used to prove the following
enhancement of Step 1.

\

\medskip

\noindent
{\bf Step 2.} \quad Let $n$ be a positive integer and let $L, P \in
\sff{Rep}(\Gamma)$ be some
finite dimensional representations. Then we can find
a finite dimensional representation $W \in \sff{Rep}(\Gamma)$ and an
injection of $\Gamma$-modules $L \hookrightarrow W$ so that the
induced map $H^{2}(S,p_{*}(\mathbb{P}^{\vee}\otimes \mathbb{L})) \to
H^{2}(S,p_{*}(\mathbb{P}^{\vee}\otimes \mathbb{W}))$ vanishes
identically.

Indeed, the goodness of $\Pi$ together with the fact that
$p_{*}({\mathbb P}^{\vee}\otimes \mathbb{L})$ is finite dimensional
again implies the
existence of an injection
$p_{*}({\mathbb P}^{\vee}\otimes {\mathbb L}) \hookrightarrow \mathbb{M}$
into a finite dimensional local system ${\mathbb M}$ on $S$, so that the
induced map
\[
H^{2}(S,p_{*}(\mathbb{P}^{\vee}\otimes\mathbb{L})) \to
H^{2}(S,p_{*}\mathbb{M})
\]
is identically zero. Now define a local system ${\mathbb W} \to X$
as the pushout
\[
\xymatrix@M+8pt{p^{*}p_{*}({\mathbb P}^{\vee}\otimes
\mathbb{L})\otimes {\mathbb P} \ar[d] \ar@{^{(}->}[r] &
p^{*}\mathbb{M}\otimes \mathbb{P} \ar[d] \\
\mathbb{L} \ar@{^{(}->}[r] & \mathbb{W}
}
\]
where $p^{*}p_{*}({\mathbb P}^{\vee}\otimes
\mathbb{L})\otimes {\mathbb P} \to \mathbb{L}$ is the natural morphism
corresponding to $\operatorname{id}_{p_{*}(\mathbb{P}^{\vee}\otimes
\mathbb{L})}$ under the identifications
\[
\opi{Hom}(p_{*}(\mathbb{P}^{\vee}\otimes
\mathbb{L}),p_{*}(\mathbb{P}^{\vee}\otimes \mathbb{L})) =
\opi{Hom}(p^{*}p_{*}(\mathbb{P}^{\vee}\otimes
\mathbb{L}),\mathbb{P}^{\vee}\otimes \mathbb{L}) =
\opi{Hom}(p^{*}p_{*}(\mathbb{P}^{\vee}\otimes
\mathbb{L})\otimes \mathbb{P},\mathbb{L}).
\]
Note that the definition of $\mathbb{W}$ implies that the natural map
$p^{*}p_{*}(\mathbb{P}^{\vee}\otimes \mathbb{L}) \hookrightarrow
\mathbb{P}^{\vee}\otimes \mathbb{L}$ factors through $p^{*}\mathbb{M}$
and so by pushing forward to $S$ and using the projection formula we
get a commutative diagram of finite dimensional local systems on $S$:
\[
\xymatrix@M+8pt{p_{*}(\mathbb{P}^{\vee}\otimes \mathbb{L})
\ar@{^{(}->}[r] \ar[dr] & \mathbb{M} \ar[d] \\
& p_{*}(\mathbb{P}^{\vee}\otimes \mathbb{W}).
}
\]
Thus the inclusion $\mathbb{L} \hookrightarrow \mathbb{W}$ will induce
the zero map
\[
H^{2}(S,p_{*}(\mathbb{P}^{\vee}\otimes
\mathbb{L})) \stackrel{0}{\longrightarrow}
H^{2}(S,p_{*}(\mathbb{P}^{\vee}\otimes
\mathbb{W}))
\]
as claimed. Step 2 is proven.

\

\medskip

\noindent
{\bf Step 3.} \quad Let $L \in \sff{Rep}(\Gamma)$ be a finite
dimensional representation and let
\[
\xymatrix@M=8pt{
Y \ar@{^{(}->}[r] & X \ar[d]^-{p} \\
& S
}
\]
be a fibration corresponding to the sequence \eqref{good-seq}. Then
there exists an injection $L \hookrightarrow W$ into a finite
dimensional $W \in \sff{Rep}(\Gamma)$ inducing the zero map
\[
H^{1}(Y,\mathbb{L}_{|Y}) \stackrel{0}{\longrightarrow}
 H^{1}(Y,\mathbb{W}_{|Y})
\]
on the fiberwise cohomology.

In order to construct $\mathbb{W}$ we begin by choosing an injection
$\imath_{Y} : \mathbb{L}_{|Y} \hookrightarrow \mathbb{E}_{Y}$ into a suitable
local system $\mathbb{E}_{Y} \to Y$, so that the natural map
\[
H^{1}(Y,\mathbb{L}_{|Y}) \to H^{1}(Y,\mathbb{E}_{Y}),
\]
that $\imath_{Y}$ induces on cohomology is identically zero.

We will choose $\mathbb{E}_{Y}$ as follows. Start with the trivial
rank one local system $\mathbb{I}_{Y} \to Y$ and consider the trivial
local system
\[
\mathbb{P}_{Y} := H^{1}(Y,\mathbb{L}_{|Y})\otimes_{\mathbb C}
\mathbb{I}_{Y}
\]
on $Y$. Let
\[
(e) \; 0 \longrightarrow \mathbb{L}_{|Y}
\stackrel{\imath_{Y}}{\longrightarrow}  \mathbb{E}_{Y} \longrightarrow
\mathbb{P}_{Y}
\longrightarrow 0
\]
be the associated tautological extension of local systems on $Y$ (the
universal extension of $\mathbb{P}_{Y}$ by $\mathbb{L}_{|Y}$).
The definition of $\mathbb{E}_{Y}$ ensures that the pushforward of any
extension $0 \to \mathbb{L}_{|Y} \to \mathbb{L}' \to \mathbb{I}_{Y}
\to 0$ via the map $\imath_{Y} : \mathbb{L}_{|Y} \to \mathbb{E}_{Y}$ will be
split. In particular the map $H^{1}(Y,\mathbb{L}_{|Y}) \to
H^{1}(Y,\mathbb{E}_{Y})$ induced by $\imath_{Y}$ is identically zero.

Observe also that we can view $\mathbb{P}_{Y}$ as the restriction
$\mathbb{P}_{Y} = \mathbb{P}_{|Y}$ of the global local system
$\mathbb{P} := p^{*}(R^{1}p_{*}\mathbb{L})$.
If it happens that the extension $(e)$ is also a restriction from a
global extension of local systems on $X$, we can take as $W$ the
monodromy representation
of this global local system and this will complete the proof of
Step 2.

By universality we know that $e \in
\opi{Ext}^{1}_{Y}(\mathbb{P}_{Y},\mathbb{L}_{|Y})$ will be fixed by
the monodromy action of $\Pi$
and so
\[
e \in \opi{Ext}^{1}_{Y}(\mathbb{P}_{|Y},\mathbb{L}_{|Y})^{\Pi} =
H^{0}(S,R^{1}p_{*}(\mathbb{P}^{\vee}\otimes \mathbb{L})).
\]
The group $H^{0}(S,R^{1}p_{*}(\mathbb{P}^{\vee}\otimes \mathbb{L}))$
fits in the following chunk of the Leray spectral sequence for the
map $p : X \to S$:
\[
0 \to H^{1}(S,p_{*}(\mathbb{P}^{\vee}\otimes \mathbb{L})) \to
H^{1}(X,\mathbb{P}^{\vee}\otimes \mathbb{L}) \to
H^{0}(S,R^{1}p_{*}(\mathbb{P}^{\vee}\otimes \mathbb{L}))
\stackrel{\delta}{\to} H^{2}(S,p_{*}(\mathbb{P}^{\vee}\otimes
\mathbb{L})).
\]
The obstruction for lifting $\mathbb{E}_{Y}$ to a global local system
on $X$ is precisely the class
\[
\delta(e) \in
H^{2}(S,p_{*}(\mathbb{P}^{\vee}\otimes \mathbb{L})).
\]
We can kill this
obstruction by using Step 2. Indeed, by the
Step 2 we can find an injection $L
\hookrightarrow V$ of finite dimensional $\Gamma$
modules inducing the zero map
\[
H^{2}(S,p_{*}(\mathbb{P}^{\vee}\otimes \mathbb{L}))
\stackrel{0}{\longrightarrow} H^{2}(S,p_{*}(\mathbb{P}^{\vee}\otimes
\mathbb{V})).
\]
Since the morphism $\mathbb{L} \hookrightarrow \mathbb{V}$ induces a
morphism of long exact sequences
\[
\xymatrix{\ldots
\ar[r] &
H^{1}(X,\mathbb{P}^{\vee}\otimes \mathbb{L}) \ar[r] \ar[d] &
H^{0}(S,R^{1}p_{*}(\mathbb{P}^{\vee}\otimes \mathbb{L}))
\ar[r]^-{\delta} \ar[d] &  H^{2}(S,p_{*}(\mathbb{P}^{\vee}\otimes
\mathbb{L})) \ar[d]_-{0} \ar[r] & \ldots \\
\ldots
\ar[r] &
H^{1}(X,\mathbb{P}^{\vee}\otimes \mathbb{V}) \ar[r] &
H^{0}(S,R^{1}p_{*}(\mathbb{P}^{\vee}\otimes \mathbb{V}))
\ar[r]^-{\delta} &  H^{2}(S,p_{*}(\mathbb{P}^{\vee}\otimes
\mathbb{V})) \ar[r] & \ldots
}
\]
we conclude that the  push-forward
\[
\xymatrix{ 0 \ar[r] & \mathbb{L}_{|Y} \ar[r] \ar[d] & \mathbb{E}_{Y}
\ar[r] \ar[d] & \mathbb{P}_{|Y} \ar[r] \ar@{=}[d] & 0 \\
0 \ar[r] & \mathbb{V}_{|Y} \ar[r] & {\mathbb W}_{Y} \ar[r] &
\mathbb{P}_{|Y} \ar[r] & 0
}
\]
of the extension $(e)$ via
the map ${\mathbb L}_{|Y} \hookrightarrow {\mathbb V}_{|Y}$
is the restriction to $Y$ of some global extension
\[
0 \to \mathbb{V} \to \mathbb{W} \to \mathbb{P} \to 0
\]
of finite dimensional local systems on $X$. Now by constrcution the
natural map $H^{1}(Y,\mathbb{L}_{|Y}) \to H^{1}(Y,\mathbb{W}_{|Y})$,
induced by the injection $\mathbb{L} \hookrightarrow \mathbb{W}$, will
vanish identically, which completes the proof of Step 3.

\

\medskip

\noindent
By putting together the existence results established in Steps 1-3
we can now finish the proof of the theorem. Fix a positive integer $n$
and let $L \in \sff{Rep}(\Gamma)$ be a finite dimensional
representation of $\Gamma$. By Step 3 we can find an injection $L
\hookrightarrow W$ in a finite dimensional representation $W \in
\sff{Rep}(\Gamma)$, so that the induced map $R^{1}p_{*}\mathbb{L} \to
R^{1}p_{*}\mathbb{W}$ vanishes identically. The Leray spectral
sequence for $p$ yields a
commutative diagram with exact rows:
\[
\xymatrix{\ldots \ar[r] & H^{n}(S,p_{*}\mathbb{L}) \ar[r] \ar[d] &
H^{n}(X,\mathbb{L})
\ar[r] \ar[d] & H^{n-1}(S,R^{1}p_{*}\mathbb{L}) \ar[d]_-{0} \ar[r] &
\ldots \\
\ldots \ar[r] & H^{n}(S,p_{*}\mathbb{W}) \ar[r] &
H^{n}(X,\mathbb{W})
\ar[r] & H^{n-1}(S,R^{1}p_{*}\mathbb{W})  \ar[r] &
\ldots
}
\]
and so
\[
H^{n}(X,\mathbb{L}) \to
\operatorname{im}[H^{n}(S,p_{*}\mathbb{W}) \to
H^{n}(X,\mathbb{W})] \subset H^{n}(X,\mathbb{W}).
\]
Furthermore, by Step 1 we can find a finite dimensional local system
$\mathbb{W}'$ and an injection $\mathbb{W} \hookrightarrow
\mathbb{W}'$ inducing
\[
H^{n}(S,p_{*}\mathbb{W}) \stackrel{0}{\longrightarrow}
H^{n}(S,p_{*}\mathbb{W}').
\]
By functoriality of the Leray spectral sequence this implies that
the natural map
\[
H^{n}(X,\mathbb{L}) \to H^{n}(X,\mathbb{W}),
\]
induced
by the injection $\mathbb{L} \hookrightarrow \mathbb{W} \subset
\mathbb{W}'$, is identically zero. The theorem is proven. \hfill $\Box$

\

\begin{rmk} \label{rem-extension-finite} {\bf (i)} The only property of
free groups that was used in the proof of the previous theorem is the
fact that free groups have cohomological dimension $\leq 1$. Therefore
over $\mathbb{C}$ we can use the same reasoning (in fact one only
needs Step 1) as in the proof of the theorem to argue that the
extension of a good group by a finite group will also be good.

\

\smallskip

\noindent
{\bf (ii)} By a classical result of M.Artin every point $x$ on a
complex smooth variety $X$ admits a Zariski neighborhood $x \in U
\subset X$ which is a tower of smooth affine morphisms of relative
dimension one. In particular the underlying topological space of $U$
(for the classical topology) is a $K(\pi,1)$ with $\pi$ being a
successive extension of free groups of finite type. Since $\pi$ is
manifestly of type $(\boldsymbol{F})$ the previous theorem implies that
$\pi$ is algebraically good
over ${\mathbb C}$.  In fact an easy application of the Leray spectral
sequence shows that any successive extension of free groups of finite
type will be of type $(\boldsymbol{F})$ and so the previous theorem
again implies that such successive extensions are algebraically good.
\end{rmk}

\subsection{Lefschetz exactness}

Let $f : X \longrightarrow Y$ be a morphism of connected topological
spaces, and $n\geq 1$ be an integer. We will say that $f$ is an
$n$-epimorphism if induced morphism $\pi_{i}(X) \longrightarrow
\pi_{i}(Y)$ is an isomorphism for $i<n$ and an epimorphism for
$i=n$. In the same way, one defines the notion of an $n$-epimorphism
of psht. A morphism $f : F \longrightarrow F'$ of psht is an
$n$-epimorphism if the induced morphism of sheaves $\pi_{i}(F)
\longrightarrow \pi_{i}(F')$ is an isomorphism for $i<n$ and an
epimorpism for $i=n$.

\begin{prop}\label{p4}
Let $f : X \longrightarrow Y$ be an $n$-epimorphism of pointed and
connected topological spaces (in $\mathbb{U}$), 
then the induced morphism $f : \sch{X}{k} \longrightarrow \sch{Y}{k}$
is an $n$-epimorphism of psht. 
\end{prop}
{\bf Proof:} For $n=1$ the proposition is clear, as
$\pi_{1}(\sch{X}{k})\simeq \pi_{1}(X)^{\op{alg}}$ and as the functor
$G \mapsto G^{\op{alg}}$ preserves epimorphisms. So we will assume
$n>1$.

One first observes that $f$ is an $n$-epimorphism if and only if it
satisfies the following two conditions:
\begin{itemize}
\item The induced morphism $\pi_{1}(X) \longrightarrow \pi_{1}(Y)$ is
  an isomorphism. 
\item  For any local system of abelian
groups $M$ on $Y$, the induced morphism $H^{i}(Y,M) \longrightarrow
H^{i}(X,f^{*}M)$ 
is an isomorphism for $i<n$ and a monomorphism for $i=n$.
\end{itemize}

In the same way, $f : \sch{X}{k} \longrightarrow \sch{Y}{k}$
is an $n$-epimorphism 
if and only if the two following conditions are satisfied.

\begin{itemize}
\item The induced morphism $\pi_{1}(\sch{X}{k}) \longrightarrow
  \pi_{1}(\sch{Y}{k})$ is an isomorphism. 

\item For any local system of $k$-vector spaces of finite dimension
$L$ on $Y$, the induced morphism $H^{i}(Y,L) \longrightarrow
H^{i}(X,f^{*}L)$ is an isomorphism for $i<n$ and a monomorphism for
$i=n$.
\end{itemize}

The proposition now follows immediately from the fundamental property
of the schematization \cite[Definition 3.3.1]{t1}. \hfill $\Box$  

\begin{cor}\label{c6}
Let $Y$ be a smooth projective (connected) complex variety of
dimension $n+1$, and $X \hookrightarrow Y$ a smooth hyperplane
section. Then, for any base point $x \in X$, the induced morphism
$$\pi_{i}(\sch{X}{k},x) \longrightarrow \pi_{i}(\sch{Y}{k},f(x))$$ is
an isomorphism for $i<n$ and an epimorphism for 
$i=n$.
\end{cor}
{\bf Proof:} This follows from Proposition \ref{p4} and the Lefschetz
theorem 
on hyperplane sections. 
\hfill $\Box$ 

\subsection{Homotopy fibers of schematizations}

Let $X \longrightarrow B$ be a morphism of pointed and connected
$\mathbb{U}$-simplicial sets, and suppose that its homotopy
fiber $Z$ is also connected. We want to relate the schematization of
$Z$ with the homotopy fiber $F$ of $\sch{X}{k} \longrightarrow
\sch{Y}{k}$.  The universal property of the schematization
induces a natural morphism of psht
\[
\sch{Z}{k} \longrightarrow F
\]
 which in general is far from
being an isomorphism. The following proposition gives a sufficient
condition in order for this morphism to be an isomorphism.

\begin{prop}\label{p5}
Assume that $Z$ is $1$-connected and of $\mathbb{Q}$-finite type. Then
the natural morphism $\sch{Z}{k} \longrightarrow F$ is an equivalence.
\end{prop}
{\bf Proof:} Note first that $\sch{Z}{k} \simeq (Z\otimes
k)^{\op{uni}}$, as $Z$ is simply connected 
(see \cite[Corollary 3.3.8]{t1}).

The fibration sequence $Z \rightarrow X \rightarrow B$ is classified
by a morphism (well defined in the homotopy category) $B
\longrightarrow B\mathbb{R}\underline{\op{Aut}}(Z)$, where
$\mathbb{R}\underline{\op{Aut}}(Z)$ is the space of auto-equivalences
of $Z$. We consider the morphism $B\mathbb{R}\underline{\op{Aut}}(Z)
\longrightarrow B\mathbb{R}\underline{\op{Aut}}((Z\otimes
k)^{\op{uni}})$, and we observe that
$B\mathbb{R}\underline{\op{Aut}}((Z\otimes k)^{\op{uni}})$ is the
pointed simplicial set of global sections of the pointed stack
$B\mathbb{R}\mathcal{AUT}((Z\otimes k)^{\op{uni}})$.  The finitness
assumtion on $Z$ implies that $B\mathbb{R}\mathcal{AUT}((Z\otimes
k)^{\op{uni}})$ is a psht, and therefore the morphism $B
\longrightarrow B\mathbb{R}\underline{\op{Aut}}((Z\otimes
k)^{\op{uni}})$ can be inserted in a (homotopy) commutative square of
pointed simplicial sets
\[
\xymatrix{
B \ar[d] \ar[r] & B\mathbb{R}\underline{\op{Aut}}((Z\otimes
k)^{\op{uni}}) \ar[d] \\ 
(B\otimes k)^{\op{sch}} \ar[r] & B\mathbb{R}\mathcal{AUT}((Z\otimes
k)^{\op{uni}}).} 
\]
By passing to the associated fibration sequences, one gets a morphism
of fibration sequences of pointed simplicial presheaves
\[
\xymatrix{
Z \ar[r] \ar[d] & X \ar[d] \ar[r] & B \ar[d] \\
\sch{Z}{k} \ar[r] & \widetilde{F} \ar[r] & \sch{B}{k}.}
\]
We deduce from this diagram a morphism between the Leray spectral
sequences, for any local system $L$ on $B$ of finite dimensional
$k$-vector spaces
\[
\xymatrix{ H^{p}(B,H^{q}(Z,L)) \ar[d]  \ar@{=>}[r] & H^{p+q}(X,L) \ar[d] \\
H^{p}((B\otimes k)^{sch},H^{q}((Z\otimes k)^{sch},L)) \ar@{=>}[r] &
H^{p+q}(\widetilde{F},L).} 
\]
The finiteness condition on $Z$ ensures that this morphism is an
isomorphisms on the $E_{2}$-term, showing that $X \longrightarrow
\widetilde{F}$ is a $\boldsymbol{P}$-equivalence. Therefore, the
induced morphism $\sch{X}{k} \longrightarrow \widetilde{F}$ is an
equivalence of psht, which induces the proposition. \hfill $\Box$ 

\

\medskip

The important consequence of Proposition \ref{p5} is that under these
conditions the 
schematization of the base $\sch{B}{k}$ \emph{acts} on the
schematization of the 
fiber $\sch{Z}{k}$. By this we mean that the fibration
sequence $\sch{Z}{k} \longrightarrow \sch{Z}{k} \longrightarrow \sch{B}{k}$
gives rise to a classification morphism 
\[
\sch{B}{k} \longrightarrow B\mathbb{R}\mathcal{AUT}(\sch{Z}{k}),
\]
where $\mathbb{R}\mathcal{AUT}$ denotes the stack of
aut-equivalences. Using \cite[Theorem 1.4.3]{t1} this morphism 
can also be considered as a morphism of $H_{\infty}$-stacks
\[
\Omega_{*}(B\otimes k)^{sch} \longrightarrow
\mathbb{R}\mathcal{AUT}(\sch{Z}{k}), 
\]
or in other words as an action loop stack $\Omega_{*}\sch{B}{k}$ on
$\sch{Z}{k}$. This action contains of course the monodromy action of
$\pi_{1}(B)$ on the homotopy groups of $\sch{Z}{k}$, but also
higher homotopical invariants, as for examples higher monodromy maps
\[
\pi_{i}(\sch{B}{k})\otimes \pi_{n}(\sch{Z}{k})
\longrightarrow 
\pi_{n+i-1}(\sch{Z}{k}).
\]
A typical situation of applications of Proposition \ref{p5} is the
case where 
$X \longrightarrow B$ is a smooth projective familly of simply
connected complex 
varieties over a smooth and projective base. In this situation, the
Hodge decomposition 
constructed in \cite{kpt} is compatible with the monodromy
action of 
$\sch{B}{\mathbb{C}}$ on the  schematization of the fiber. In other
words, if $Z$ is the homotopy type 
of the fiber, then $\sch{Z}{\mathbb{C}}$ has an action of
$\Omega_{*}\sch{B}{\mathbb{C}}$ which 
is $\mathbb{C}^{*}$-equivariant with respect to the Hodge
decomposition. This can 
also be interpreted by saying that $\sch{Z}{\mathbb{C}}$ form a
\textit{variation of non-abelian schematic Hodge structures on $B$}.
Of course, this variation contains more information than the
variations of Hodge structures on the rational homotopy 
groups of the fiber, since it captures higher homotopical data. An
example of such 
non-trivial 
higher invariants is given in \cite{s}.

\end{document}